\documentclass[11pt,a4paper,reqno]{article}

\usepackage{amsmath,amsthm,amssymb,mathtools}

\usepackage{parskip}  

\usepackage{fullpage}
\usepackage{enumitem}
 
\usepackage{url}
\usepackage[latin1]{inputenc}

\usepackage[dvipsnames]{xcolor}
\usepackage{hyperref,graphicx}

\usepackage{pgfplots} 
\usepackage{caption}
\usepackage{subcaption}
 
\usepackage{booktabs}

\usepackage{scalerel,stackengine}
\newcommand\pig[1]{\scalerel*[5.5pt]{\Big#1}{%
  \ensurestackMath{\addstackgap[1.5pt]{\big#1}}}}
\newcommand\pigl[1]{\mathopen{\pig{#1}}}
\newcommand\pigr[1]{\mathclose{\pig{#1}}}

\DeclareMathOperator{\Li}{Li}

\DeclareMathOperator{\sgn}{sgn}

\newtheorem{theorem}{Theorem}%[section]

\newtheorem{lemma}[theorem]{Lemma}

\newtheorem{proposition}[theorem]{Proposition} 

\newtheorem{definition}[theorem]{Definition}

\newtheorem{remark}[theorem]{Remark}
\newtheorem{remarks}[theorem]{Remarks} 
\numberwithin{theorem}{section}

\theoremstyle{remark}
\theoremstyle{definition}

 \author{Malin Pal\"o Forsstr\"om\thanks{Chalmers University of Technology and the University of Gothenburg, Sweden}}

\begin{document}

\title{Poisson representations for tree-indexed Markov chains} 
 % 60J10

\maketitle

\begin{abstract}
	In~\cite{fgs}, the class of Poisson representable processes was introduced. Several well-known processes were shown not to belong to this class, with examples including both the Curie Weiss model and the Ising model on \( \mathbb{Z}^2 \) for certain choices of parameters. Curiously,  it was also shown that all positively associated \( \{ 0,1 \}\)-valued Markov chains do belong to this class. In this paper, we interpolate between Markov chains and Ising models by considering tree-indexed Markov chains. In particular, we show that for any finite tree that is not a path, whether or not the corresponding tree-indexed Markov chain is representable always depends on the parameters. Moreover, we give an example of a family of infinite trees such that the corresponding tree-indexed Markov chains are representable for some non-trivial parameters. In addition, we give alternative proofs and arguments and also strengthen several of the results in~\cite{fgs}.
	\end{abstract}
 
\section{Introduction}

In~\cite{fgs}, Poisson representable processes were introduced. 
To define this family of processes, let \(S\) be a finite or countably infinite set.  We think of the set \( S\) as the index set of our process. Next,  let \(\nu \) be a $\sigma$-finite measure on \(\mathcal{P}(S)\backslash \{\emptyset\},\) where $\mathcal{P}(S)$ is the power set of $S$. This generates a \(\{0,1\}\)-valued process $X^\nu= (X^\nu_i)_{i\in S}$ defined as follows.
Let  $Z$ be a Poisson process on $\mathcal{P}(S)\backslash \{\emptyset\}$ with intensity measure $\nu,$ written \( Z \sim \mathrm{Poisson}(\nu). \) Then $Z $ is a random collection of nonempty subsets of $S$, and for each set \( B \in \mathcal{P}(S)\backslash \{ \emptyset \}\) independently, we have \( B \in Z \) with probability \( 1-e^{-\nu(B)}.\)  We define $(X^\nu_i)_{i\in S}$ by
$$
X^\nu_i= X^\nu(i)=  \begin{cases}
    1 &\text{if } i\in \bigcup_{B\in Z} B \cr 
    0 &\text{otherwise.}
\end{cases}
$$ 
We let \( \mathcal{R}\) denote the set of all processes $(X_i)_{i\in S}$ which are equal (in distribution) to $X^\nu$ for some measure $\nu$ on \( \mathcal{P}(S)\backslash \{ \emptyset \}.  \) In words, \( X \in \mathcal{R} \) if \( X \) is the indicator function of the union of some set-valued Poisson process on its index set.
The purpose of~\cite{fgs} was to try to understand which processes could be represented in this way, and in this paper, we deepen the understanding of this question by studying tree-indexed Markov chains in more detail.

In~\cite[Theorem 3.1]{fgs}, it was shown that all positively associated \( \{0,1\}\)-valued Markov chains belong to \( \mathcal{R} .\) In contrast, by~\cite[Theorem 6.3]{fgs}, the Ising model on \( \mathbb{Z}^d, \) \(d \geq 2 ,\) does not belong to \( \mathcal{R} \) if the coupling parameter is sufficiently small (see also~\cite{fs}). By a similar proof, it was also shown that tree-indexed Markov chains, defined below, are not always in \( \mathcal{R}. \) 
In particular, by~\cite[Theorem 6.1]{fgs}, a tree-indexed Markov chain \( X \) on a tree with a vertex \( v \) from which at least \( d \) distinct infinite paths emerge is not in \( \mathcal{R} \) whenever \( P(X_v=0) < 1/(1-r_2^{(d)}), \) where~\( r_2^{(d)} \) is an explicit constant that approaches \( 0 \) as \( d \to \infty. \) %In particular, \( r_2^{(3)}=-1, \) \(r_2^{(4)} = \sqrt{3}-2, \) and  \( r_2^{(5)} =2\sqrt{6}-5 .\) 
In this paper, we refine this picture by showing that, in fact, all finite trees exhibit a phase transition, with a subset of the parameter space corresponding to processes in \( \mathcal{R} \) and another subset corresponding to processes that do not belong to \( \mathcal{R}. \)  
In order to be able to present this result in more detail, we now define what we mean by tree-indexed Markov chains.

\begin{definition}\label{def: TIMC}
	Let \( T \) be a tree with vertex set  \( V( T) \), edge set \( E(T)\), and finite vertex degree, and let \( r,p \in [0,1].\) 
	Fix a vertex \( o \in V(T). \) For each vertex \( v \in V(T), \) let 
	\begin{equation*}
		R_v = \begin{cases}
			0 &\text{w.p. } r \cr 
			1 &\text{w.p. } 1-r
		\end{cases} 
	\end{equation*}
	be independent random variables. 
	The \emph{tree-indexed Markov chain \( X \) indexed by \( T \) with parameters \( (r,p)  \)} is defined as follows. 
	\begin{enumerate}[label=(\Alph*)]
		\item Let \( X(o) = R_o, \) and set \( V_0 \coloneqq \{ o \}. \) 
	For \( j = 1,2,\dots, \) let \( V_j \) be the set of vertices \( v \in V(T) \smallsetminus V_{j-1}\) that are adjacent to some vertex \( v^- \in V_{j-1}, \) and for each \( v \in V_j \),  let 
	\begin{equation*}
		X(v) = \begin{cases}
			X(v^-) &\text{w.p. } 1-p \cr 
			R_v & \text{w.p. }p_.
		\end{cases}
	\end{equation*}
	\end{enumerate}
	Equivalently, given \( o \) and \( (R_v)_{v \in V(T)}, \) we can construct a random variable with the same distribution as \( X \) as follows.
	\begin{enumerate}[label=(\Alph*), resume]
		\item Let \( Y \) be a random subset of \( E(T) \) obtained by removing each edge in \( E(T) \), independently, with probability \( p \). 
			Independently for each connected component \( C \subseteq V(T) \) of the resulting random subgraph of \( T\),  set \( X(v) = R_{v_0} \) for all \( v \in C, \) where \( v_0 \) is the unique vertex in \( C \) that is closest to the origin \( o .\) 
	\end{enumerate} 
\end{definition}

Since (A) and (B) above yield random processes with same distribution, we will use them both interchangeably throughout the paper.  
We note that using (A), the process \( X = (X(v))_{v \in V(T)}\) can be thought of as being sampled sequentially as a reversible and positively associated Markov chain along any path of the tree. Using (B), the tree-indexed Markov chain \( X \) is instead defined as a so-called \emph{divide-and-color process} \cite{st2019}, and this is equivalently the FK (Fortuin-Kasteleyn) representation corresponding to the tree-indexed-Markov chain.

The model described in Definition~\ref{def: TIMC} is also known as a symmetric binary channel on trees and appears, e.g., when studying the broadcasting/reconstruction problems on trees (see, e.g., ~\cite{mp2003}). Further, it can be interpreted as an Ising model on a tree with free boundary conditions, where e.g. \( r=1/2\) and \( p  = \sqrt\frac{2}{1 + e^{4J}} \) corresponds to an Ising model with coupling parameter \( J \) and no external field. Further, it is well known (see, e.g., \cite[Proposition 12.24]{gbook}), that the extremal states of the Ising model on trees are tree-indexed Markov chains.

Throughout the rest of this paper, for \( n \geq 1, \) we let \( r^{(n)} \) be the largest (strictly) negative root of the polylogarithm function \( \mathrm{Li}_{1-n}(z)  \) of index \( 1-n \), defined by the power series \(  \sum_{k=1}^\infty \frac{z^k}{k^{1-n}}.\)  We note in particular that it is is known that \( r^{(n)} \in [-1,0) \) for all \( n ,\) and that \( r^{(n)} \nearrow 0 \) as \( n \to \infty. \) %
Also, for \( n \geq 1, \) let \( \tilde B_n \) be the \( n \)th complementary Bell number, defined by \( \tilde B_n = \sum_{k=0}^n (-1)^k S(n,k), \) where \( S(n,k) \) is the Stirling number of the second kind, i.e. the number of ways to partition a set of \( n \) elements into \( k \) non-empty sets. Using the complementary Bell numbers, for integers \( k,m \geq 1, \) we define 
	\begin{equation}\label{eq: r0 and r1}
		r_0(k) \coloneqq \max_{\substack{j \leq k \colon  (-1)^{j} \tilde B_{ j}>0}}
		\frac{\bigl( (-1)^{j} \tilde B_{ j}\bigr)^{1/(j-1)}}{1+ \bigl( (-1)^{j} \tilde B_{j}\bigr)^{1/(j-1)}} \quad \text{and} \quad r_1(m) \coloneqq 1/(1-r^{(m)}).
	\end{equation}
Note that it is well known that \( \lim_{n \to \infty} r_0(n) = \lim_{n \to \infty} r_1(n) = 1. \)

	%
%We collect \( r^{(n)}\) and \( \tilde B_n \) for small values of \( n \) in Table~\ref{table: boundaries}.

\begin{table}[h]\centering
	\begin{tabular}{c c c c c }
		\( n \) &  \( r^{(n)} \) & \( \tilde B_n \) & \( r_0(n) \) & \( r_1(n) \) \\\toprule
		3 & \( -1 \)& 1 & 9 & \(1/2\)\\ 
		4 & \(-0.26795\dots \) & \( 1 \) &  \( 1/2\) & \(0.78868\dots \)\\
		5 & \( -0.10102\dots \)& \(-2\) & \( 0.54321\dots \) & \(0.90825\dots \)\\
		6 & \( -0.04310 \dots \)& \( -9 \) & \( 0.54321\dots \) & \(0.95868\dots \)\\
		7 & \( -0.01952 \dots \) & \( -9 \) & \( 0.59054\dots \)& \(0.98085 \dots \)\\
		8 & \( -0.00915 \dots \) & \( 50 \) & \( 0.63619\dots \)& \(0.99093\dots \)
	\end{tabular}
	\caption{The values of \( r^{(n)} \), \( \tilde B_n \), \( r_0(n) \) and \( r_1(n) \) for \( n = 3,4,\dots, 8. \)} \label{table: boundaries}
\end{table}

By~\cite[Theorem 3.1]{fgs}, any \( \{0,1\}\)-valued Markov chain is in \( \mathcal{R}, \) and hence, equivalently, any tree-indexed Markov chain on a finite or infinite graph with no cycles where each vertex has degree at most two is Poisson representable. Our first main result shows that this result is sharp in the sense that for any finite graph with at least one vertex of degree at least three, whether or not the corresponding Markov chain is representable always depends on the parameters.

\begin{theorem}\label{theorem: finite trees}
	Let \( T \) be a finite tree, let \( (r,p) \in (0,1) , \) and let \( X \) be a tree-indexed Markov chain on \( T \) with parameters \( (r,p). \) Let \( k \) be the size of the boundary of \( T\) (the set of vertices of degree exactly one in \( T \)) and let \( m \) be the maximal degree of any vertex in \( V(T). \) 
	Then the following holds.
	\begin{enumerate}[label=(\alph*)] 
		\item \( X \in \mathcal{R} \) if \( r > r_0(k)\) and \( p \) is sufficiently close to zero (depending on \( r \)).\label{theorem: finite trees ii}
		\item \( X \notin \mathcal{R} \) if \( r < r_0(k) \) and \( p \) is sufficiently close to zero (depending on \( r \)).\label{theorem: finite trees iii}
		\item  \( X \in \mathcal{R} \) if \( r > r_1(m) \) and  \( p \) sufficiently close to one (depending on \( r \)).\label{theorem: finite trees iv}
		\item \( X \notin \mathcal{R} \) if \( r < r_1(m) \) and \( p \) sufficiently close to one (depending on \( r \)).\label{theorem: finite trees i}
	\end{enumerate}   
\end{theorem}

A consequence of Theorem~\ref{theorem: finite trees} is that the threshold in \( r \) for small \( p \) does not depend on the whole tree but only on its local structure. In particular, the threshold is the same for a large binary tree as for a tree with only four vertices and one vertex of degree~three.

We next consider infinite trees. By extending~\cite[Lemma 2.14]{fgs} to signed measures in Lemma~\ref{lemma: finite to infinite new}, the negative results of Theorem~\ref{theorem: finite trees}\ref{theorem: finite trees i} and Theorem~\ref{theorem: finite trees}\ref{theorem: finite trees iii} immediately extend to infinite trees. However, one might ask if there is an infinite tree \( T \) and parameters \( r,p \in (0,1) \) such that the tree-indexed Markov chain on \( T \) with parameters \( (r,p) \) is in \( \mathcal{R}. \) This is the main content of our next theorem.
In this theorem, we show that for a special family of infinite trees, which we will refer to as \emph{octopus trees}, a positive result similar to parts \ref{theorem: finite trees ii} and \ref{theorem: finite trees iv} of Theorem~\ref{theorem: finite trees} still holds.

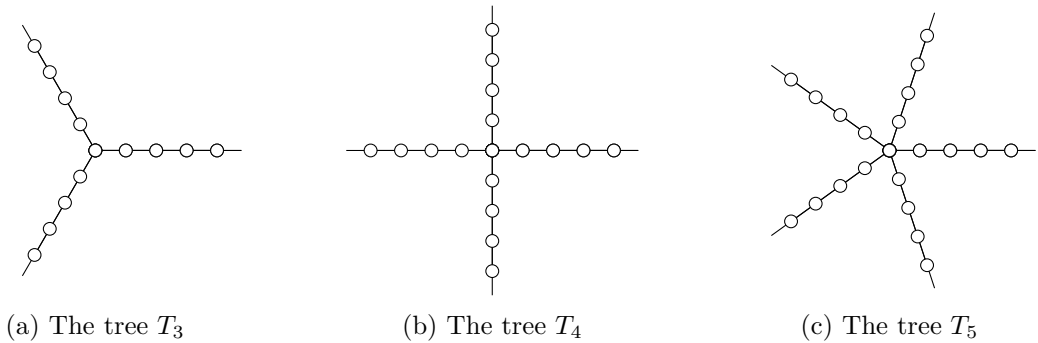
\begin{figure}[h]
    \centering
    \begin{subfigure}[t]{0.32\textwidth}
        \centering
		\begin{tikzpicture}[scale=.4] 
    	
    	\draw[white] (-4.8,-4.8) -- (4.8,4.8); 
    	
    	\foreach \i in {0,...,3}
    	{	
    		\draw[] (0,0) -- ({4.8*cos(120*\i)},{4.8*sin(120*\i)}); 
    		\draw (0,0) -- ({4*cos(120*\i)},{4*sin(120*\i)});  
    		\foreach \j in {0,...,4}
    		{
        		\filldraw[fill=white] ({\j*cos(120*\i)},{\j*sin(120*\i)}) circle (6pt);  
    		}   
    	}
     \end{tikzpicture}
		\caption{The tree \( T_3\)} 
	\end{subfigure}
    \begin{subfigure}[t]{0.32\textwidth}
    \centering
    \begin{tikzpicture}[scale=.4] 
        \centering
		\draw[white] (-4.8,-4.8) -- (4.8,4.8); 
    	
    	\foreach \i in {0,...,4}
    	{	
    		\draw[] (0,0) -- ({4.8*cos(90*\i)},{4.8*sin(90*\i)}); 
    		\draw (0,0) -- ({4*cos(90*\i)},{4*sin(90*\i)});  
    		\foreach \j in {0,...,4}
    		{
        		\filldraw[fill=white] ({\j*cos(90*\i)},{\j*sin(90*\i)}) circle (6pt);  
    		} 
    		 
    	}
     \end{tikzpicture}
	\caption{The tree \( T_4\)} 
	\end{subfigure}
    \begin{subfigure}[t]{0.32\textwidth}
        \centering
		\begin{tikzpicture}[scale=.4]

    	\draw[white] (-4.8,-4.8) -- (4.8,4.8); 
    	
    	\foreach \i in {0,...,5}
    	{	
    		\draw[] (0,0) -- ({4.8*cos(72*\i)},{4.8*sin(72*\i)}); 
    		\draw (0,0) -- ({4*cos(72*\i)},{4*sin(72*\i)});  
    		\foreach \j in {0,...,4}
    		{
        		\filldraw[fill=white] ({\j*cos(72*\i)},{\j*sin(72*\i)}) circle (6pt);  
    		}  
    		  
    	}
     \end{tikzpicture}
		\caption{The tree \( T_5\)} 
	\end{subfigure}
	\caption{The octopus trees \( T_3, \) \( T_4, \) and \( T_5. \)}
	\label{figure: trees intro}
\end{figure}

For \( m \geq 3 ,\) let \( T_m\) denote the infinite tree consisting of \( m \) disjoint semi-infinite paths whose endpoints are joined at one vertex 
%with one vertex of degree \( m\) and all other vertices having degree two 
(see Figure~\ref{figure: trees intro}). The tree \( T_m\) will be referred to as the \emph{octopus tree of degree~\( m . \)}

%Our next main result shows that for any \( m \geq 3,\) there is a non-trivial phase transition in the parameter space between a region where the corresponding tree-indexed Markov chains on \( T_m \) is in \( \mathcal{R}, \) and a region where it is not.

\begin{theorem}\label{theorem: general octopus trees}
    Let \( m \geq 3\) and \( r,p \in (0,1), \) and let \( X\) be a tree-indexed Markov chain  on \( T_m\) with parameters \( (r,p).\) Then there exists a constant \( r_2(m)  \in (0,1), \) not depending on \( p , \) such that 
    \begin{enumerate}[label=(\alph*)]
        \item \( X \notin \mathcal{R} \) if  \( r < r_1(m) ,\) where \( r_1(m) \) is defined in~\eqref{eq: r0 and r1}, and\label{item: general octopus trees a}
        \item \( X \in \mathcal{R}\) if \( r \geq r_2(m).\) \label{item: general octopus trees b}
    \end{enumerate} 
    Moreover,  
    \begin{enumerate}[label=(\alph*),resume]
    	\item \( r_2(3) = r_1(3)= 1/2. \)\label{item: general octopus trees c}
    \end{enumerate}
\end{theorem}

 The first part of this theorem, Theorem~\ref{theorem: general octopus trees}\ref{item: general octopus trees a}, is an immediate consequence of~\cite[Theorem 6.1]{fgs}, which also gives a concrete lower bound for \( r_2(m). \) The main contribution of Theorem~\ref{theorem: general octopus trees} is thus part~\ref{item: general octopus trees b}, which together with~\ref{item: general octopus trees a} implies that the set of Markov chains indexed by octopus trees has a phase transition, and part~\ref{item: general octopus trees c}, which shows that when \( m = 3, \) whether or not \( X \) is in \( \mathcal{R} \) is independent of~\( p. \) 
Moreover, we note that for finite trees, whether or not \( X\in \mathcal{R} \) typically depends on both \( r \) and \(p, \) see, e.g., \cite[Figure 1]{fgs}.

The main idea of the proof of Theorem~\ref{theorem: general octopus trees}\ref{item: general octopus trees a} is to first show that for any \( m \geq 3, \) there exists a number \(  r_1(m) \) and a finite subtree \( T \) of \( T_m ,\) such that for all fixed \( r \leq r_1(m) \) and all \( p \) sufficiently close to one, the restriction of \( X \) to \( T \) is not in  \(  \mathcal{R}. \) Using symmetries of the tree (see Proposition~\ref{proposition: scaling}), it follows that this implies that \( X \notin  R \) for all \( p \in (0,1). \) We elaborate on this and similar symmetry properties in Section~\ref{section: scaling}, but note here that the main observation of this section is that results for \( p \) close to one or \( p \) close to zero can sometimes be extended to results for general \( p \) and the same value of \( r. \)

Our final main result, Theorem~\ref{theorem: p close to zero}, gives versions of the previous results for general infinite trees.

\begin{theorem}\label{theorem: p close to zero}

	Let \( T \) be an infinite tree, let \( r,p \in (0,1), \) and let \( X \) be the tree-indexed Markov chain on \( T \) with parameters~\( (r,p). \) 	
	 Then, with \( r_0 \) and \( r_1  \) as in~\eqref{eq: r0 and r1}, the following holds.
	\begin{enumerate}[label=(\alph*)]
		\item Let \( n_1<\infty\) be smaller than or equal to the size of the boundary of \( T. \) Then \( X \notin \mathcal{R} \) if \( r < r_0(n_1) \) and \( p \) is sufficiently close to zero.
		\label{theorem: infinite trees i} 
		
		\item Let \( n_2<\infty \) be such that there is some vertex \( o \in V(T) \) with at least \( n_2 \) eventually disjoint infinite paths in \( E(T) \)  starting at \( o. \)  Then \( X \notin \mathcal{R} \) if \( r< r_1(n_2) \) and \( p \) is sufficiently close to zero. \label{theorem: infinite trees ii}  
		
		\item Let \( n_3<\infty \) be such that there is some vertex \( o \in V(T) \) with at least \( n_3 \) disjoint infinite paths in \( E(T) \)  starting at \( o. \) Then \( X \notin \mathcal{R} \) if \( r < r_1(n_3). \) \label{theorem: infinite trees iii} 
	\end{enumerate} 
	In particular, if \( T \) has an infinite boundary, then \( X \notin \mathcal{R} \)  for all  \( p \) sufficiently close to zero.
	 
\end{theorem}

%Since tree-indexed Markov chains appear when considering the Ising model on trees, it is natural to ask what consequences~Theorem~\ref{theorem: general octopus trees} and~Theorem~\ref{theorem: growing trees} has for these models. 

%Our second resolved this question completely.
%\begin{theorem}\label{theorem: ising tree new new}
%	Let \( X \) be any state of the Ising model on a regular tree. Then \( X \notin \mathcal{R}. \)
%\end{theorem}

As a part of the proofs of our main results, we extend and provide new, shorter, and more natural proofs of several of the fundamental results in~\cite{fgs}. In particular, we mention Lemma~\ref{lemma: trees}, which substantially strengthens, but at the same time greatly simplifies the proof of~\cite[Theorem 3.1]{fgs}, Proposition~\ref{prop: connected}, which extends~\cite[Proposition~3.7]{fgs} to signed measures, and~Proposition~\ref{proposition: new negative lemma for p close to one}, which extends~\cite[Theorem 5.6]{fgs} from trees with exactly one vertex of degree \( \geq 3 \) to general finite trees.

\begin{remark}
	Our main results all use the same parameters \( r = P(Y(v)=0) \) and \( p \) for all vertices and edges. However, many of the proofs throughout the paper allow the parameters to be different for different vertices and edges, and thus, many of the main results can, with small modifications, be proven in a more general setting. This is useful when, e.g., considering a Markov process \( X \) on a graph \( G, \)   since for any \( S \subseteq V(G)\), if the subgraph of \( G \) induced by \( S \) is a tree, then  \( Y = X(S)\mid \{ X(S^c) \equiv 0\} \) is a tree-indexed Markov chain where \( P(Y(v)=0) \) is typically not the same for all \( v \in S. \) 
\end{remark}

The rest of this paper will be structured as follows. 
In Section~\ref{section: preliminaries}, we review the notation and some fundamental results about Poisson representable processes from~\cite{fgs}. We also introduce the notation we will use when working with graphs and trees, give a more general definition for tree-indexed Markov chains which will simplify the proofs throughout the paper, and list a few useful properties of the polylogarithm function. Next, in Section~\ref{section: signed connected}, we extend~\cite[Proposition 3.7]{fgs} to signed measures; i.e. we show that if \( X \) satisfies a Markov property and if \( \nu \) is a signed measure corresponding to \( X ,\) then \( \nu \) has support only on connected sets.  
In~Section~\ref{section: markov formulas}, we use our results from Section~\ref{section: signed connected} to give simpler formulas for \( \nu \) in terms of \( X \) for tree-indexed Markov chains. These generalise~\cite[Theorem~3.1]{fgs} from line graphs to trees, and in addition, give more intuitive proofs also for the simpler case of a line (when the tree-indexed Markov chain is a regular Markov chain).
Next, in Section~\ref{section: scaling}, we show that on self-similar trees (see Definition~\ref{def: self-similar}), we can transfer results from one set of parameters \( (r,p) \) to other parameters \( (r,p'). \) 
In Section~\ref{section: finite and p zero}, we consider the case \( p \approx 0, \) and provide proofs of Theorem~\ref{theorem: finite trees}\ref{theorem: finite trees ii} and Theorem~\ref{theorem: finite trees}\ref{theorem: finite trees iii}. 
In Section~\ref{section: finite and p one}, we consider the case \( p \approx 1, \) and prove Theorem~\ref{theorem: finite trees}\ref{theorem: finite trees iv} and Theorem~\ref{theorem: finite trees}\ref{theorem: finite trees i}.
In Section~\ref{section: octopus trees}, we give a proof of Theorem~\ref{theorem: general octopus trees}\ref{item: general octopus trees b} and Theorem~\ref{theorem: general octopus trees}\ref{item: general octopus trees c}. 
Finally, in Section~\ref{section: general infinite trees}, we provide a proof of Theorem~\ref{theorem: p close to zero}.

\section{Preliminaries}\label{section: preliminaries}

In this section, we introduce notation and recall fundamental results from~\cite{fgs} that will be useful to us in this paper. Further, we introduce the notations we will use for graphs and trees, and also give a more general definition of tree-indexed Markov chains, which will be useful in later sections of this paper.

\subsection{Definitions and properties for Poisson representable processes}
Let a finite or countably infinite set \( S \) be given. For \( i \in S, \) we let 
$$
    \mathcal{S}_i \coloneqq \bigl\{ A \in \mathcal{P}(S)\smallsetminus \{ \emptyset \} \colon i \in A \bigr\},
$$
and for $A\subseteq S$, we let
$$
    \mathcal{S}_A^\cup \coloneqq \bigcup_{i\in A} \mathcal{S}_i 
\quad \text{and} \quad  \mathcal{S}_A^\cap \coloneqq \bigcap_{i\in A} \mathcal{S}_i .
$$
In words, \( \mathcal{S}_i \) the set of all subsets of \( S \) that contain \( i \), \( \mathcal{S}_A^\cup \) is the set of all subsets of \( S \) that have a non-empty intersection with \( A \), and \( \mathcal{S}_A^\cap \) is the set of all subsets of \( S \) that contain \( A \) as a subset.
This notation is useful to us because given a measure \( \nu \) on \( \mathcal{P}(S)\smallsetminus \{ \emptyset \}, \) it connects probabilities involving \( X^\nu \) with the measure \( \nu \) in the sense that for any set non-empty \( A \subseteq S, \) one has
\begin{equation}\label{eq: P vs nu}
	P\bigl( X^\nu(A) \equiv 0) = e^{-\nu(\mathcal{S}_A^\cup)}.
\end{equation}
(This follows immediately from the definition of \( X^\nu.\)) 
We will use this identity to go back and forth between the two random objects \( X^\nu \) and \( Y \sim \mathrm{Poisson}(\nu). \)
In particular, the identity~\eqref{eq: P vs nu} is used in the proof of the following lemma from~\cite{fgs}, which uses the M\"obius inversion theorem to obtain a formula for \( \nu.\) % in terms of the probabilities of \( X\) being identically zero on sets.
 
\begin{lemma}[Lemma 2.12 in~\cite{fgs}]\label{lemma: unique signed}
    Let \( S \) be finite, and let \( X = (X_s)_{s \in S}\) be 
    $\{0,1\}$-valued random variables such that \( P\bigl( X(I)\equiv 0)>0\) for all nonempty \( I \subseteq S.\) Then there is a unique signed measure \( \nu \) on \( \mathcal{P}(S)\smallsetminus \{ \emptyset \}\) that satisfies 
    \begin{equation}\label{eq: step 1}
        \nu\bigl(\mathcal{S}_{I}^\cup\bigr) = -\log P\bigl( X(I) \equiv 0 \bigr),\quad I \subseteq S.
    \end{equation} 
    Moreover, \( \nu \) is given by
    \begin{equation}\label{eq: Mobius inversion result}
        \begin{split}
        &\nu(K) 
        = \sum_{I \subseteq K} (-1)^{|K|-|I|} 
         \log P\bigl( X(S\smallsetminus I) \equiv 0 \bigr)  
        , \qquad \emptyset \neq K\subseteq S.
        \end{split}
    \end{equation} 
    Consequently, if \( \nu \geq 0,\) then \( X = X^\nu \in \mathcal{R}\).
\end{lemma}

Whenever \( S \) is a finite set, \( X = (X_s)_{s\in S}\) is a \(\{0,1\}\)-valued process, and \(\nu\) is given by~\eqref{eq: Mobius inversion result}, we say that \emph{\(\nu\) is the unique signed measure corresponding to \(X.\)}

A useful feature of the Poisson representation of a random process, if such a representation exists, is that it works well with various natural restrictions; a property that is in large inherited from the underlying Poisson process. This is the content of the following lemma.

\begin{lemma}[Lemma 2.14 in~\cite{fgs}]\label{lemma: finite to infinite}
	Let \( S \) be finite or countably infinite, and consider a process~$X= (X_s)_{s\in S}$. 
\begin{enumerate}[label=(\alph*)]

    \item If \( X=X^\nu \in \mathcal{R}\)  and \( B \subseteq S, \) then there is a positive measure \( \nu_B \) such that \( X|_{B} = X^{\nu_B}.\) Moreover, for any non-empty measurable subset \( \mathcal{A} \subseteq \mathcal{P}(B), \) we have \( \nu_B(\mathcal{A}) = \nu\bigl( \{ A' \in \mathcal{P}(S) \colon  A' \cap B\in \mathcal{A} \} \bigr).\)\label{item: finite to infinite i}

    \item If \( X=X^\nu \in \mathcal{R}\) and \( B \subseteq S,\) then there is a measure \( \nu_{B,0} \) on \( \mathcal{P}(B)\smallsetminus \{ \emptyset \} \) such that \( X| \bigl\{ X(B^c)\equiv 0 \bigr\} = X^{\nu_{B,0}}.\) Moreover, \( \nu_{B,0} = \nu|_{\mathcal{P}(B)}. \) \label{item: finite to infinite iii}
    
    \item\label{item: finite to infinite ii} If there exist $S_1\subseteq S_2\subseteq,\ldots$ such that $S=\bigcup_i S_i$ and  $X_{S_n} = X^{\nu_n} $ for some $\nu_n \geq 0$, then $X=X^\nu \in \mathcal{R}$ for some $\nu \geq 0$. (The projection of $\nu$ on to each $S_n$ will simply be  $\nu_n$.)
\end{enumerate} 
\end{lemma}

\subsection{Notation for graphs and trees}

In this paper we will often consider processes on graphs. When we do this, we will always assume that the graphs are connected. For graphs \( G\) we will use the following notation.
\begin{itemize}
	\item We let \( V(G) \) denote the vertex set of \( G \) and  \( E(G) \) denote the edge set of \( G. \) 
	
	\item We say that a set \( V \subseteq V(G) \) is connected if it induces a connected subgraph of \( V(G). \)
	
	\item If \( v,v' \in V(G) \) are neighbors, i.e. if \( (v,v') \in E(G)\), we write \( v \sim v'. \)
   
    \item Given a set \( V \subseteq V(G), \) we let \( E_V \) denote the set of all edges \( e \in E(G) \) with at least one endpoint in \( V. \) If \( V = \{ v \}, \) we let \( E_v \coloneqq E_{\{ v \}}. \)
     
	\item For a finite and non-empty set \( V \subseteq  V(G) \), we let 
	\[\mathcal{B}^-(V) \coloneqq  \bigl\{ v \in V\colon |\{ v' \in V \colon  v \sim v' \}| = 1 \bigr\} \qquad \text{(the inner boundary)} \]  
    %\[\mathcal{B}^-(V) \coloneqq  \bigl\{ v \in V\colon \exists  v' \in V(G) \smallsetminus V \text{ s.t. } v \sim v'  \bigr\} \qquad \text{(the inner boundary)} \]  
    and
    \[
    \partial V  \coloneqq \mathcal{B}^+(V)  \coloneqq \bigl\{ v \in V(G) \smallsetminus V\colon \exists  v' \in V \text{ s.t. } v \sim v'  \bigr\} \qquad\text{(the outer boundary).} \] 
    %and
    %\[
    %\mathcal{B}( V) \coloneqq \mathcal{B}^-(V) \cup \mathcal{B}^+(V).     \] 
    %
    For a set \( V \subseteq V(G), \) we let 
    \[
    N(V) \coloneqq \{ v' \in V(G)\smallsetminus V \colon v' \sim v \text{ for some } v \in V\} \qquad\text{(the neighborhood of \( V\))}.
    \]
    For  \( I \subseteq V, \) we let
    \[
    \partial_{V} I \coloneqq \{ v \in \mathcal{B}^+( V) \colon \exists v' \in I \text{ s.t. } v \sim v' \} = \mathcal{B}^+(V)\cap N(I).
    \]
\end{itemize}

Now assume that we are given a tree \( T \) which is rooted in some vertex \( o \in V(T). \)
\begin{itemize}
	\item  If \( V \) is connected, we let \( T_V \) be the subtree of \( T \) induced by \( V. \) 
	\item Given an edge \( e \in E(T), \) we let \( e^- \) and \( e^+ \) denote the endpoint of \(e\) that is closer resp. further away from~\(o\).
	
	\item Given a vertex \( v \in V(T)\smallsetminus \{ o \}, \) we let \( E_v^+ \) be the set of all edges \( e \in E(T) \) with \( e^- = v. \) 
\end{itemize}

	\begin{figure}[h]

		\begin{subfigure}[t]{.3\textwidth}\centering
		
		\hfil
		\begin{tikzpicture}[scale=.3] 
    	 
 		% level 6
    	\foreach \i in {0,...,3}
    		\foreach \j in {-1,1}
    			\foreach \k in {-1,1}
    				\foreach \l in {-1,1}
    					\foreach \m in {-1,1}
    						\foreach \n in {-1,1}
    						{
    							\draw[help lines] ({5*cos(120*\i+30*\j+15*\k+8*\l+4*\m)},{5*sin(120*\i+30*\j+15*\k +8*\l+4*\m)})
    							--
    							({5.8*cos(120*\i+30*\j+15*\k+8*\l+4*\m+2*\n)},{5.8*sin(120*\i+30*\j+15*\k +8*\l+4*\m+2*\n)});  
    						} 
    	
    	% level 5
    	\foreach \i in {0,...,3}
    		\foreach \j in {-1,1}
    			\foreach \k in {-1,1}
    				\foreach \l in {-1,1}
    					\foreach \m in {-1,1}
    					{
    						\draw[help lines] ({4*cos(120*\i+30*\j+15*\k+8*\l)},{4*sin(120*\i+30*\j+15*\k +8*\l)})
    						--
    						({5*cos(120*\i+30*\j+15*\k+8*\l+4*\m)},{5*sin(120*\i+30*\j+15*\k +8*\l+4*\m)});  
    						
    						\filldraw[help lines, fill=white] ({5*cos(120*\i+30*\j+15*\k+8*\l+4*\m)},{5*sin(120*\i+30*\j+15*\k +8*\l+4*\m)}) circle (4pt); 
    					}  
    	
    	% level 4
    	\foreach \i in {0,...,3}
    		\foreach \j in {-1,1}
    			\foreach \k in {-1,1}
    				\foreach \l in {-1,1}
    				{
    					\draw[help lines] ({3*cos(120*\i+30*\j+15*\k)},{3*sin(120*\i+30*\j+15*\k )}) -- ({4*cos(120*\i+30*\j+15*\k+8*\l)},{4*sin(120*\i+30*\j+15*\k +8*\l)});  
    					\filldraw[help lines, fill=white] ({4*cos(120*\i+30*\j+15*\k+8*\l)},{4*sin(120*\i+30*\j+15*\k +8*\l)}) circle (5pt);   
    				}

    	% level 3
    	\foreach \i in {0,...,3}
    		\foreach \j in {-1,1}
    			\foreach \k in {-1,1}
    			{
    				\draw[help lines] ({2*cos(120*\i+30*\j)},{2*sin(120*\i+30*\j)}) -- ({3*cos(120*\i+30*\j+15*\k)},{3*sin(120*\i+30*\j+15*\k )}); 
    				\filldraw[help lines, fill=white] ({3*cos(120*\i+30*\j+15*\k)},{3*sin(120*\i+30*\j+15*\k)}) circle (6pt);   
    			}   
    	
    	% level 2
    	\foreach \i in {0,...,3} 
    		\foreach \j in {-1,1}
    		{
    			\draw[help lines] ({cos(120*\i)},{sin(120*\i)}) -- ({2*cos(120*\i+30*\j)},{2*sin(120*\i+30*\j)}); 
    			\filldraw[help lines, fill=white] ({2*cos(120*\i+30*\j)},{2*sin(120*\i+30*\j)}) circle (7pt);   
    		}   
		
    	% level 1
    	\foreach \i in {0,...,3}
    	{	  
    		\draw[help lines] (0,0) -- ({cos(120*\i)},{sin(120*\i)});
    		\filldraw[help lines, fill=white] ({cos(120*\i)},{sin(120*\i)}) circle (8pt);  
    	}
    	
    	% level 0 
    	\filldraw[help lines, fill=white] (0,0) circle (10pt); 
    	
    	\filldraw[fill=black] (0,0) circle (10pt);  
    	
    	\foreach \i in {0,...,3}
    	{	   
    		%\draw (0,0) -- ({cos(120*\i)},{sin(120*\i)});
    		\filldraw[fill=black] ({cos(120*\i)},{sin(120*\i)}) circle (8pt);  
    	}

    	\foreach \i in {0,...,1} 
    		\foreach \j in {-1,1}
    		{
    			%\draw ({cos(120*\i)},{sin(120*\i)}) -- ({2*cos(120*\i+30*\j)},{2*sin(120*\i+30*\j)}); 
    			\filldraw[fill=black] ({2*cos(120*\i+30*\j)},{2*sin(120*\i+30*\j)}) circle (7pt);  
    		}

    	\foreach \i in {2} 
    		\foreach \j in {1}
    		{
    			%\draw ({cos(120*\i)},{sin(120*\i)}) -- ({2*cos(120*\i+30*\j)},{2*sin(120*\i+30*\j)}); 
    			\filldraw[fill=black] ({2*cos(120*\i+30*\j)},{2*sin(120*\i+30*\j)}) circle (7pt);    
    		}
    	
    	\foreach \i in {0}
    		\foreach \j in {-1,1}
    			\foreach \k in {-1,1}
    			{ 
    				%\draw ({2*cos(120*\i+30*\j)},{2*sin(120*\i+30*\j)}) -- ({3*cos(120*\i+30*\j+15*\k)},{3*sin(120*\i+30*\j+15*\k )}); 
    				\filldraw[fill=black] ({3*cos(120*\i+30*\j+15*\k)},{3*sin(120*\i+30*\j+15*\k)}) circle (6pt);   
    			}
    	
    	\foreach \i in {0}
    		\foreach \j in {1}
    			\foreach \k in {-1,1}
    			{ 
    				%\draw ({2*cos(120*\i+30*\j)},{2*sin(120*\i+30*\j)}) -- ({3*cos(120*\i+30*\j+15*\k)},{3*sin(120*\i+30*\j+15*\k )}); 
    				\filldraw[draw=black, fill=black] ({3*cos(120*\i+30*\j+15*\k)},{3*sin(120*\i+30*\j+15*\k)}) circle (6pt);  
    				
    				%\draw[red] ({3*cos(120*\i+30*\j+15*\k)},{3*sin(120*\i+30*\j+15*\k)}) circle (10pt);  
    			}
    	
    	\foreach \i in {0}
    		\foreach \j in {-1}
    			\foreach \k in {-1}
    			{ 
    				\draw ({2*cos(120*\i+30*\j)},{2*sin(120*\i+30*\j)}) -- ({3*cos(120*\i+30*\j+15*\k)},{3*sin(120*\i+30*\j+15*\k )}); 
    				\filldraw[draw=black, fill=black] ({3*cos(120*\i+30*\j+15*\k)},{3*sin(120*\i+30*\j+15*\k)}) circle (6pt);  
    				%\draw[red] ({3*cos(120*\i+30*\j+15*\k)},{3*sin(120*\i+30*\j+15*\k)}) circle (10pt);  
    			}
    	
    	\foreach \i in {1}
    		\foreach \j in {-1}
    			\foreach \k in {1}
    			{ 
    				%\draw ({2*cos(120*\i+30*\j)},{2*sin(120*\i+30*\j)}) -- ({3*cos(120*\i+30*\j+15*\k)},{3*sin(120*\i+30*\j+15*\k )}); 
    				\filldraw[draw=black, fill=black] ({3*cos(120*\i+30*\j+15*\k)},{3*sin(120*\i+30*\j+15*\k)}) circle (6pt); 
    				%\draw[red] ({3*cos(120*\i+30*\j+15*\k)},{3*sin(120*\i+30*\j+15*\k)}) circle (10pt);     
    			}  
    	
    	\foreach \i in {1}
    		\foreach \j in {1}
    			\foreach \k in {-1,1}
    			{ 
    				%\draw ({2*cos(120*\i+30*\j)},{2*sin(120*\i+30*\j)}) -- ({3*cos(120*\i+30*\j+15*\k)},{3*sin(120*\i+30*\j+15*\k )}); 
    				\filldraw[draw=black, fill=black]  ({3*cos(120*\i+30*\j+15*\k)},{3*sin(120*\i+30*\j+15*\k)}) circle (6pt);     
    			}

    	\foreach \i in {2}
    		\foreach \j in {1}
    			\foreach \k in {-1,1}
    			{ 
    				%\draw ({2*cos(120*\i+30*\j)},{2*sin(120*\i+30*\j)}) -- ({3*cos(120*\i+30*\j+15*\k)},{3*sin(120*\i+30*\j+15*\k )}); 
    				\filldraw[draw=black, fill=black] ({3*cos(120*\i+30*\j+15*\k)},{3*sin(120*\i+30*\j+15*\k)}) circle (6pt);    
    			}

    	\foreach \i in {0}
    		\foreach \j in {-1}
    			\foreach \k in {1}
    				\foreach \l in {-1,1}
    				{
    					%\draw ({3*cos(120*\i+30*\j+15*\k)},{3*sin(120*\i+30*\j+15*\k )}) -- ({4*cos(120*\i+30*\j+15*\k+8*\l)},{4*sin(120*\i+30*\j+15*\k +8*\l)});  
    					\filldraw[draw=black, fill=black] ({4*cos(120*\i+30*\j+15*\k+8*\l)},{4*sin(120*\i+30*\j+15*\k +8*\l)}) circle (5pt);

    					%\draw[red] ({4*cos(120*\i+30*\j+15*\k+8*\l)},{4*sin(120*\i+30*\j+15*\k +8*\l)}) circle (9pt);   
    				} 
	
     \end{tikzpicture}

		\caption{\( T \) in gray and \( V \) in black }
		\end{subfigure} 
		\hfil
		\begin{subfigure}[t]{.3\textwidth}\centering
		
		\hfil
		\begin{tikzpicture}[scale=.3] 
    	 
 		% level 6
    	\foreach \i in {0,...,3}
    		\foreach \j in {-1,1}
    			\foreach \k in {-1,1}
    				\foreach \l in {-1,1}
    					\foreach \m in {-1,1}
    						\foreach \n in {-1,1}
    						{
    							\draw[help lines] ({5*cos(120*\i+30*\j+15*\k+8*\l+4*\m)},{5*sin(120*\i+30*\j+15*\k +8*\l+4*\m)})
    							--
    							({5.8*cos(120*\i+30*\j+15*\k+8*\l+4*\m+2*\n)},{5.8*sin(120*\i+30*\j+15*\k +8*\l+4*\m+2*\n)});  
    						} 
    	
    	% level 5
    	\foreach \i in {0,...,3}
    		\foreach \j in {-1,1}
    			\foreach \k in {-1,1}
    				\foreach \l in {-1,1}
    					\foreach \m in {-1,1}
    					{
    						\draw[help lines] ({4*cos(120*\i+30*\j+15*\k+8*\l)},{4*sin(120*\i+30*\j+15*\k +8*\l)})
    						--
    						({5*cos(120*\i+30*\j+15*\k+8*\l+4*\m)},{5*sin(120*\i+30*\j+15*\k +8*\l+4*\m)});  
    						
    						\filldraw[help lines, fill=white] ({5*cos(120*\i+30*\j+15*\k+8*\l+4*\m)},{5*sin(120*\i+30*\j+15*\k +8*\l+4*\m)}) circle (4pt); 
    					}  
    	
    	% level 4
    	\foreach \i in {0,...,3}
    		\foreach \j in {-1,1}
    			\foreach \k in {-1,1}
    				\foreach \l in {-1,1}
    				{
    					\draw[help lines] ({3*cos(120*\i+30*\j+15*\k)},{3*sin(120*\i+30*\j+15*\k )}) -- ({4*cos(120*\i+30*\j+15*\k+8*\l)},{4*sin(120*\i+30*\j+15*\k +8*\l)});  
    					\filldraw[help lines, fill=white] ({4*cos(120*\i+30*\j+15*\k+8*\l)},{4*sin(120*\i+30*\j+15*\k +8*\l)}) circle (5pt);   
    				}

    	% level 3
    	\foreach \i in {0,...,3}
    		\foreach \j in {-1,1}
    			\foreach \k in {-1,1}
    			{
    				\draw[help lines] ({2*cos(120*\i+30*\j)},{2*sin(120*\i+30*\j)}) -- ({3*cos(120*\i+30*\j+15*\k)},{3*sin(120*\i+30*\j+15*\k )}); 
    				\filldraw[help lines, fill=white] ({3*cos(120*\i+30*\j+15*\k)},{3*sin(120*\i+30*\j+15*\k)}) circle (6pt);   
    			}   
    	
    	% level 2
    	\foreach \i in {0,...,3} 
    		\foreach \j in {-1,1}
    		{
    			\draw[help lines] ({cos(120*\i)},{sin(120*\i)}) -- ({2*cos(120*\i+30*\j)},{2*sin(120*\i+30*\j)}); 
    			\filldraw[help lines, fill=white] ({2*cos(120*\i+30*\j)},{2*sin(120*\i+30*\j)}) circle (7pt);   
    		}   
		
    	% level 1
    	\foreach \i in {0,...,3}
    	{	  
    		\draw[help lines] (0,0) -- ({cos(120*\i)},{sin(120*\i)});
    		\filldraw[help lines, fill=white] ({cos(120*\i)},{sin(120*\i)}) circle (8pt);  
    	}
    	
    	% level 0 
    	\filldraw[help lines, fill=white] (0,0) circle (10pt); 
    	
    	\foreach \i in {0,...,3}
    	{	   
    		\draw (0,0) -- ({cos(120*\i)},{sin(120*\i)}); 
    	}

    	\foreach \i in {0,...,1} 
    		\foreach \j in {-1,1}
    		{
    			\draw ({cos(120*\i)},{sin(120*\i)}) -- ({2*cos(120*\i+30*\j)},{2*sin(120*\i+30*\j)});  
    		}

    	\foreach \i in {2} 
    		\foreach \j in {1}
    		{
    			\draw ({cos(120*\i)},{sin(120*\i)}) -- ({2*cos(120*\i+30*\j)},{2*sin(120*\i+30*\j)}); 
    		}
    	
    	\foreach \i in {0}
    		\foreach \j in {-1,1}
    			\foreach \k in {-1,1}
    			{ 
    				\draw ({2*cos(120*\i+30*\j)},{2*sin(120*\i+30*\j)}) -- ({3*cos(120*\i+30*\j+15*\k)},{3*sin(120*\i+30*\j+15*\k )});  
    			}
    	
    	\foreach \i in {0}
    		\foreach \j in {1}
    			\foreach \k in {-1,1}
    			{ 
    				\draw ({2*cos(120*\i+30*\j)},{2*sin(120*\i+30*\j)}) -- ({3*cos(120*\i+30*\j+15*\k)},{3*sin(120*\i+30*\j+15*\k )}); 
    			}
    	
    	\foreach \i in {0}
    		\foreach \j in {-1}
    			\foreach \k in {-1,1}
    			{ 
    				\draw ({2*cos(120*\i+30*\j)},{2*sin(120*\i+30*\j)}) -- ({3*cos(120*\i+30*\j+15*\k)},{3*sin(120*\i+30*\j+15*\k )}); 
    			}

    	\foreach \i in {1}
    		\foreach \j in {1}
    			\foreach \k in {-1,1}
    			{ 
    				\draw ({2*cos(120*\i+30*\j)},{2*sin(120*\i+30*\j)}) -- ({3*cos(120*\i+30*\j+15*\k)},{3*sin(120*\i+30*\j+15*\k )});    
    			}

    	\foreach \i in {2}
    		\foreach \j in {1}
    			\foreach \k in {-1,1}
    			{ 
    				\draw ({2*cos(120*\i+30*\j)},{2*sin(120*\i+30*\j)}) -- ({3*cos(120*\i+30*\j+15*\k)},{3*sin(120*\i+30*\j+15*\k )});   
    			}

    	\foreach \i in {1}
    		\foreach \j in {-1}
    			\foreach \k in {1}
    			{ 
    				\draw ({2*cos(120*\i+30*\j)},{2*sin(120*\i+30*\j)}) -- ({3*cos(120*\i+30*\j+15*\k)},{3*sin(120*\i+30*\j+15*\k )});   
    			}

    	\foreach \i in {0}
    		\foreach \j in {-1}
    			\foreach \k in {1}
    				\foreach \l in {-1,1}
    				{
    					\draw ({3*cos(120*\i+30*\j+15*\k)},{3*sin(120*\i+30*\j+15*\k )}) -- ({4*cos(120*\i+30*\j+15*\k+8*\l)},{4*sin(120*\i+30*\j+15*\k +8*\l)});    
    				}

    				%-----
    				
    \filldraw[fill=black] (0,0) circle (10pt);  
    	
    	\foreach \i in {0,...,3}
    	{	   
    		\filldraw[fill=black] ({cos(120*\i)},{sin(120*\i)}) circle (8pt);  
    	}

    	\foreach \i in {0,...,1} 
    		\foreach \j in {-1,1}
    		{
    			\filldraw[fill=black] ({2*cos(120*\i+30*\j)},{2*sin(120*\i+30*\j)}) circle (7pt);  
    		}

    	\foreach \i in {2} 
    		\foreach \j in {1}
    		{
    			\filldraw[fill=black] ({2*cos(120*\i+30*\j)},{2*sin(120*\i+30*\j)}) circle (7pt);    
    		}
    	
    	\foreach \i in {0}
    		\foreach \j in {-1,1}
    			\foreach \k in {-1,1}
    			{  
    				\filldraw[fill=black] ({3*cos(120*\i+30*\j+15*\k)},{3*sin(120*\i+30*\j+15*\k)}) circle (6pt);   
    			}
    	
    	\foreach \i in {0}
    		\foreach \j in {1}
    			\foreach \k in {-1,1}
    			{  
    				\filldraw[draw=black, fill=red] ({3*cos(120*\i+30*\j+15*\k)},{3*sin(120*\i+30*\j+15*\k)}) circle (6pt);    
    			}
    	
    	\foreach \i in {0}
    		\foreach \j in {-1}
    			\foreach \k in {-1}
    			{  
    				\filldraw[draw=black, fill=red] ({3*cos(120*\i+30*\j+15*\k)},{3*sin(120*\i+30*\j+15*\k)}) circle (6pt);     
    			}
    	
    	\foreach \i in {1}
    		\foreach \j in {-1}
    			\foreach \k in {1}
    			{  
    				\filldraw[draw=black, fill=red] ({3*cos(120*\i+30*\j+15*\k)},{3*sin(120*\i+30*\j+15*\k)}) circle (6pt);      
    			}  
    	
    	\foreach \i in {1}
    		\foreach \j in {1}
    			\foreach \k in {-1,1}
    			{  
    				\filldraw[draw=black, fill=red]  ({3*cos(120*\i+30*\j+15*\k)},{3*sin(120*\i+30*\j+15*\k)}) circle (6pt);     
    			}

    	\foreach \i in {2}
    		\foreach \j in {1}
    			\foreach \k in {-1,1}
    			{  
    				\filldraw[draw=black, fill=red] ({3*cos(120*\i+30*\j+15*\k)},{3*sin(120*\i+30*\j+15*\k)}) circle (6pt);    
    			}

    	\foreach \i in {0}
    		\foreach \j in {-1}
    			\foreach \k in {1}
    				\foreach \l in {-1,1}
    				{ 
    					\filldraw[draw=black, fill=red] ({4*cos(120*\i+30*\j+15*\k+8*\l)},{4*sin(120*\i+30*\j+15*\k +8*\l)}) circle (5pt);

    					%\draw[red] ({4*cos(120*\i+30*\j+15*\k+8*\l)},{4*sin(120*\i+30*\j+15*\k +8*\l)}) circle (9pt);   
    				} 
	
     \end{tikzpicture}

		\caption{\( T_V \) in black and red, and~\( \mathcal{B}^-(V) \) in red}
		\end{subfigure} 
		\hfil
		\begin{subfigure}[t]{.3\textwidth}\centering
			\begin{tikzpicture}[scale=.3] 
    	
    	% level 6
    	\foreach \i in {0,...,3}
    		\foreach \j in {-1,1}
    			\foreach \k in {-1,1}
    				\foreach \l in {-1,1}
    					\foreach \m in {-1,1}
    						\foreach \n in {-1,1}
    						{
    							\draw[help lines] ({5*cos(120*\i+30*\j+15*\k+8*\l+4*\m)},{5*sin(120*\i+30*\j+15*\k +8*\l+4*\m)})
    							--
    							({5.8*cos(120*\i+30*\j+15*\k+8*\l+4*\m+2*\n)},{5.8*sin(120*\i+30*\j+15*\k +8*\l+4*\m+2*\n)});  
    						} 
    	
    	% level 5
    	\foreach \i in {0,...,3}
    		\foreach \j in {-1,1}
    			\foreach \k in {-1,1}
    				\foreach \l in {-1,1}
    					\foreach \m in {-1,1}
    					{
    						\draw[help lines] ({4*cos(120*\i+30*\j+15*\k+8*\l)},{4*sin(120*\i+30*\j+15*\k +8*\l)})
    						--
    						({5*cos(120*\i+30*\j+15*\k+8*\l+4*\m)},{5*sin(120*\i+30*\j+15*\k +8*\l+4*\m)});  
    						
    						\filldraw[help lines, fill=white] ({5*cos(120*\i+30*\j+15*\k+8*\l+4*\m)},{5*sin(120*\i+30*\j+15*\k +8*\l+4*\m)}) circle (4pt); 
    					}  
    	
    	% level 4
    	\foreach \i in {0,...,3}
    		\foreach \j in {-1,1}
    			\foreach \k in {-1,1}
    				\foreach \l in {-1,1}
    				{
    					\draw[help lines] ({3*cos(120*\i+30*\j+15*\k)},{3*sin(120*\i+30*\j+15*\k )}) -- ({4*cos(120*\i+30*\j+15*\k+8*\l)},{4*sin(120*\i+30*\j+15*\k +8*\l)});  
    					\filldraw[help lines, fill=white] ({4*cos(120*\i+30*\j+15*\k+8*\l)},{4*sin(120*\i+30*\j+15*\k +8*\l)}) circle (5pt);   
    				}

    	% level 3
    	\foreach \i in {0,...,3}
    		\foreach \j in {-1,1}
    			\foreach \k in {-1,1}
    			{
    				\draw[help lines] ({2*cos(120*\i+30*\j)},{2*sin(120*\i+30*\j)}) -- ({3*cos(120*\i+30*\j+15*\k)},{3*sin(120*\i+30*\j+15*\k )}); 
    				\filldraw[help lines, fill=white] ({3*cos(120*\i+30*\j+15*\k)},{3*sin(120*\i+30*\j+15*\k)}) circle (6pt);   
    			}   
    	
    	% level 2
    	\foreach \i in {0,...,3} 
    		\foreach \j in {-1,1}
    		{
    			\draw[help lines] ({cos(120*\i)},{sin(120*\i)}) -- ({2*cos(120*\i+30*\j)},{2*sin(120*\i+30*\j)}); 
    			\filldraw[help lines, fill=white] ({2*cos(120*\i+30*\j)},{2*sin(120*\i+30*\j)}) circle (7pt);   
    		}   
		
    	% level 1
    	\foreach \i in {0,...,3}
    	{	  
    		\draw[help lines] (0,0) -- ({cos(120*\i)},{sin(120*\i)});
    		\filldraw[help lines, fill=white] ({cos(120*\i)},{sin(120*\i)}) circle (8pt);  
    	}
    	
    	% level 0 
    	\filldraw[help lines, fill=white] (0,0) circle (10pt);

    	% ---- 

    	\foreach \i in {0,...,3}
    	{	   
    		\draw (0,0) -- ({cos(120*\i)},{sin(120*\i)}); 
    	}

    	\foreach \i in {0,...,1} 
    		\foreach \j in {-1,1}
    		{
    			\draw ({cos(120*\i)},{sin(120*\i)}) -- ({2*cos(120*\i+30*\j)},{2*sin(120*\i+30*\j)});  
    		}

    	\foreach \i in {2} 
    		\foreach \j in {1}
    		{
    			\draw ({cos(120*\i)},{sin(120*\i)}) -- ({2*cos(120*\i+30*\j)},{2*sin(120*\i+30*\j)}); 
    		}
    	
    	\foreach \i in {0}
    		\foreach \j in {-1,1}
    			\foreach \k in {-1,1}
    			{ 
    				\draw ({2*cos(120*\i+30*\j)},{2*sin(120*\i+30*\j)}) -- ({3*cos(120*\i+30*\j+15*\k)},{3*sin(120*\i+30*\j+15*\k )});  
    			}
    	
    	\foreach \i in {0}
    		\foreach \j in {1}
    			\foreach \k in {-1,1}
    			{ 
    				\draw ({2*cos(120*\i+30*\j)},{2*sin(120*\i+30*\j)}) -- ({3*cos(120*\i+30*\j+15*\k)},{3*sin(120*\i+30*\j+15*\k )}); 
    			}
    	
    	\foreach \i in {0}
    		\foreach \j in {-1}
    			\foreach \k in {-1,1}
    			{ 
    				\draw ({2*cos(120*\i+30*\j)},{2*sin(120*\i+30*\j)}) -- ({3*cos(120*\i+30*\j+15*\k)},{3*sin(120*\i+30*\j+15*\k )}); 
    			}

    	\foreach \i in {1}
    		\foreach \j in {1}
    			\foreach \k in {-1,1}
    			{ 
    				\draw ({2*cos(120*\i+30*\j)},{2*sin(120*\i+30*\j)}) -- ({3*cos(120*\i+30*\j+15*\k)},{3*sin(120*\i+30*\j+15*\k )});    
    			}

    	\foreach \i in {2}
    		\foreach \j in {1}
    			\foreach \k in {-1,1}
    			{ 
    				\draw ({2*cos(120*\i+30*\j)},{2*sin(120*\i+30*\j)}) -- ({3*cos(120*\i+30*\j+15*\k)},{3*sin(120*\i+30*\j+15*\k )});   
    			}

    	\foreach \i in {1}
    		\foreach \j in {-1}
    			\foreach \k in {1}
    			{ 
    				\draw ({2*cos(120*\i+30*\j)},{2*sin(120*\i+30*\j)}) -- ({3*cos(120*\i+30*\j+15*\k)},{3*sin(120*\i+30*\j+15*\k )});   
    			}

    	\foreach \i in {0}
    		\foreach \j in {-1}
    			\foreach \k in {1}
    				\foreach \l in {-1,1}
    				{
    					\draw ({3*cos(120*\i+30*\j+15*\k)},{3*sin(120*\i+30*\j+15*\k )}) -- ({4*cos(120*\i+30*\j+15*\k+8*\l)},{4*sin(120*\i+30*\j+15*\k +8*\l)});    
    				}

    				%-----

    \filldraw[fill=black] (0,0) circle (10pt);  
    	
    	\foreach \i in {0,...,3}
    	{	   
    		\filldraw[fill=black] ({cos(120*\i)},{sin(120*\i)}) circle (8pt);  
    	}

    	\foreach \i in {0,...,1} 
    		\foreach \j in {-1,1}
    		{
    			\filldraw[fill=black] ({2*cos(120*\i+30*\j)},{2*sin(120*\i+30*\j)}) circle (7pt);  
    		}

    	\foreach \i in {2} 
    		\foreach \j in {1}
    		{
    			\filldraw[fill=black] ({2*cos(120*\i+30*\j)},{2*sin(120*\i+30*\j)}) circle (7pt);    
    		}
    	
    	\foreach \i in {2} 
    		\foreach \j in {-1}
    		{
    			\filldraw[draw=black,fill=red] ({2*cos(120*\i+30*\j)},{2*sin(120*\i+30*\j)}) circle (7pt);    
    		}
    	
    	\foreach \i in {0}
    		\foreach \j in {-1,1}
    			\foreach \k in {-1,1}
    			{  
    				\filldraw[fill=black] ({3*cos(120*\i+30*\j+15*\k)},{3*sin(120*\i+30*\j+15*\k)}) circle (6pt);   
    			}
    	  
    	\foreach \i in {1}
    		\foreach \j in {-1}
    			\foreach \k in {1}
    			{  
    				\filldraw[draw=black, fill=black] ({3*cos(120*\i+30*\j+15*\k)},{3*sin(120*\i+30*\j+15*\k)}) circle (6pt);      
    			}  
    	  
    	\foreach \i in {1}
    		\foreach \j in {-1}
    			\foreach \k in {-1}
    			{  
    				\filldraw[draw=black, fill=red] ({3*cos(120*\i+30*\j+15*\k)},{3*sin(120*\i+30*\j+15*\k)}) circle (6pt);      
    			}  
    	
    	\foreach \i in {1}
    		\foreach \j in {1}
    			\foreach \k in {-1,1}
    			{  
    				\filldraw[draw=black, fill=black]  ({3*cos(120*\i+30*\j+15*\k)},{3*sin(120*\i+30*\j+15*\k)}) circle (6pt);     
    			}

    	\foreach \i in {2}
    		\foreach \j in {1}
    			\foreach \k in {-1,1}
    			{  
    				\filldraw[draw=black, fill=black] ({3*cos(120*\i+30*\j+15*\k)},{3*sin(120*\i+30*\j+15*\k)}) circle (6pt);    
    			}

    	\foreach \i in {0}
    		\foreach \j in {-1}
    			\foreach \k in {1}
    				\foreach \l in {-1,1}
    				{ 
    					\filldraw[draw=black, fill=black] ({4*cos(120*\i+30*\j+15*\k+8*\l)},{4*sin(120*\i+30*\j+15*\k +8*\l)}) circle (5pt); 
    					  
    				}
    				
    	\foreach \i in {0}
    		\foreach \j in {-1}
    			\foreach \k in {-1}
    				\foreach \l in {-1,1}
    				{ 
    					\filldraw[draw=black, fill=red] ({4*cos(120*\i+30*\j+15*\k+8*\l)},{4*sin(120*\i+30*\j+15*\k +8*\l)}) circle (5pt); 
    					  
    				} 
    	
    	\foreach \i in {0}
    		\foreach \j in {1}
    			\foreach \k in {-1,1}
    				\foreach \l in {-1,1}
    				{ 
    					\filldraw[draw=black, fill=red] ({4*cos(120*\i+30*\j+15*\k+8*\l)},{4*sin(120*\i+30*\j+15*\k +8*\l)}) circle (5pt); 
    					  
    				}

    	\foreach \i in {1}
    		\foreach \j in {1}
    			\foreach \k in {-1,1}
    				\foreach \l in {-1,1}
    				{ 
    					\filldraw[draw=black, fill=red] ({4*cos(120*\i+30*\j+15*\k+8*\l)},{4*sin(120*\i+30*\j+15*\k +8*\l)}) circle (5pt); 
    					  
    				}

    	\foreach \i in {1}
    		\foreach \j in {-1}
    			\foreach \k in { 1}
    				\foreach \l in {-1,1}
    				{ 
    					\filldraw[draw=black, fill=red] ({4*cos(120*\i+30*\j+15*\k+8*\l)},{4*sin(120*\i+30*\j+15*\k +8*\l)}) circle (5pt); 
    					  
    				} 
    	
    	\foreach \i in {2}
    		\foreach \j in {1}
    			\foreach \k in { -1,1}
    				\foreach \l in {-1,1}
    				{ 
    					\filldraw[draw=black, fill=red] ({4*cos(120*\i+30*\j+15*\k+8*\l)},{4*sin(120*\i+30*\j+15*\k +8*\l)}) circle (5pt); 
    					  
    				}

    	\foreach \i in {0}
    		\foreach \j in {-1}
    			\foreach \k in {1}
    				\foreach \l in {-1,1}
    					\foreach \m in {-1,1}
    					{ 
    						\filldraw[draw=black, fill=red] ({5*cos(120*\i+30*\j+15*\k+8*\l+4*\m)},{5*sin(120*\i+30*\j+15*\k +8*\l +4*\m)}) circle (4pt); 
    					  
    					}

     \end{tikzpicture}

		\caption{ \( T_V \) in black, and \( \mathcal{B}^+(V) \) in red.}
		\end{subfigure}
		\hfil 
		\caption{In the figures above, we illustrate \( T_V \), \( \mathcal{B}^-(V) \) and \( \mathcal{B}^+(V) \) for a given set \( T \) and a set \( V \subseteq V(T). \)}
\end{figure}

\subsection{More general tree-indexed Markov chains}\label{section: detailed definition}

Later in the paper, it will be useful to consider tree-indexed Markov chains where the parameters \( r \) and \( p \) are allowed to be different for different vertices and edges. The following definition generalizes  Definition~\ref {def: TIMC} to this setting.

\begin{definition}\label{def: TIMC general}
	Let \( T \) be a tree with finite vertex degree, vertex set \( V( T) \), and edge set \( E(T)\).
	For each \( v \in V(T), \) let \( r_v \in [0,1], \) and for each \( e \in E(T), \) let \( p_e \in [0,1]. \)	
	Fix a vertex \( o \in V(T). \) For each vertex \( v \in V(T) \) independently, let 
	\begin{equation*}
		R_v = \begin{cases}
			0 &\text{w.p. } r_v \cr 
			1 &\text{w.p. } 1-r_v.
		\end{cases} 
	\end{equation*} 
	The \emph{tree-indexed Markov chain  \( X \) indexed by \( T \) with parameters  \( \bigl((r_v)_{v \in V(T)},(p_e)_{e \in E(T)}  \bigr) \)} is defined as follows.
	
	\begin{enumerate}[label=(\Alph*')]
		\item Let \( X(o) = R_o \), let \( V_0 \coloneqq \{ o \},\) and for   \( j = 1,2,\dots, \) let \( V_j \coloneqq V_{j-1} \cup \mathcal{B}^+(V_{j-1}) \). Sequentially, for each \( j = 1,2, \dots \), for each \( v' \in V_j \) independently, let \( v \in V_{j-1} \) be such that \( v \sim v' \) and set 
	\begin{equation*}
		X(v') = \begin{cases}
			X(v) &\text{w.p. } 1-p_{(v,v')} \cr 
			R_{v'} & \text{w.p. }p_{(v,v')}.
		\end{cases}
	\end{equation*}\label{item: sampling A'}
	\end{enumerate}

	Equivalently, given \( o \) and \( (R_v)_{v \in V(T)}, \) we can construct a random variable with the same distribution as \( X \) as follows.
	
	\begin{enumerate}[label=(\Alph*'),resume]
		\item Let \( Y \) be a random subset of \( E(T) \) obtained by removing each \( e \in E(T) \), independently, with probability \( p_e \) (i.e., let \( Y \) be edge percolation). 
			Independently for each connected component \( C \subseteq V(T) \)  of the resulting random graph, let \( v' \) be the unique vertex in \( C \) that minimizes the distance to the origin, and set \( X(v) = R_{v'} \) for all \( v \in C. \) \label{item: sampling B'}
	\end{enumerate}
\end{definition}

Note that both (A') and~(B') above result in the same process, and their corresponding constructions reduce to (A) and (B) of Definition~\ref{def: TIMC} in the special case \( (r_v) \equiv r \) and \( (p_e) \equiv p \) for some \( r,p \in [0,1]. \) Moreover, note that if \( r_v = r \) for some \( r \in [0,1] \) and all \( v \in V(T), \) then the process is independent of the choice of \( o . \)

\subsection{Properties of the polylogarithm function}\label{section: polylog}

Throughout this paper, the polylogarithm function will naturally appear. The polylogarithm function with index \( s \in \mathbb{C} \) is defined for \( z \in \mathbb{C} \) with \( |z|< 1 \) by
\[ \Li_s(z) = \sum_{k=1}^\infty \frac{z^k}{k^s} \]
and can be extended to the whole complex plane by analytic continuation. 
Some special cases include \( \Li_1(z) = -\log (1-z), \) \( \Li_0(z) = \frac{z}{1-z} \) and \( \Li_{-1}(z) = \frac{z}{(1-z)^2}. \) When the index \( s \) is a negative integer, the polylogarithm function can also be defined using other combinatorial quantities such as Stirling numbers of the second kind and Eulerian numbers. For this reason, they naturally show up in combinatorial arguments.

For \( n \geq 2, \) the polylogarithm function \( \mathrm{Li}_{1-n}(z) \) is well-known to satisfy the following properties (see, e.g., the paragraph before~\cite[Theorem~5.6]{fgs}). 
	\begin{itemize}
		\item \( \mathrm{Li}_{1-n}(0)=0. \)
		\item \( \mathrm{Li}_{1-n} \) has \( n-1 \) real negative simple roots \( r_{n-1}^{(n)}<\dots < r_1^{(n)}<0\) which are interlacing, in the sense that \( r_{j+1}^{(n)}<   r_{j}^{(n-1)} < r_j^{(n)} \) for all \( j \in \{ 1,2, \dots, n-2\}. \)  
		\item  \( \mathrm{Li}_{1-n}(z)<0 \) for \( z \in (r_1^{(n)},0). \)   
	\end{itemize}	
	
Using the notation of the introduction, we thus have \( r^{(n)} = r^{(n)}_1 .\) Recalling that we defined \( r_1(n) \coloneqq 1/(1-r^{(n)}) .\) In particular, this implies that
\begin{equation*}
	r<r_1(n) \Leftrightarrow -(1-r)/r < r^{(n)}.
\end{equation*}

A simple consequence of the above-listed properties of the polylogarithm function, which we will use multiple times, is the following lemma.
\begin{lemma}\label{lemma: polylog}
	Let \( n \geq 3 \) and \( r \in (0,1). \) Then the following holds.
	\begin{enumerate}[label=(\alph*)]
	\item If \( r >r_1(n) , \) then \( \Li_{1-j}(-(1-r)/r) < 0 \) for all \( j \in \{  1, 2,3, \dots, n\} .\)\label{lemma: polylog a}
	\item If \( r < r_1(n) , \) then there is \( j \in \{ 3, \dots, n \}  \) such that \( \Li_{1-j}(-(1-r)/r) > 0. \) \label{lemma: polylog b} 
\end{enumerate}
\end{lemma}

\begin{proof}
	We first show that~\ref{lemma: polylog a} holds. To this end, assume that \( r > r_1(n) \) or equivalently, that \( -(1-r)/r > r_1^{(n)}. \) Since \( r_1^{(n)} > r_1^{(n-1)} > \dots > r_1^{(2)}, \) it follows that \( -(1-r)/r \in ( r_1^{(n)},0) \) for all \( j \in \{ 2,3, \dots, n \}.  \) Since \( \Li_{1-j}(z)<0 \) for all \( z \in  (r_1^{(n)},0) \), this concludes the proof of~\ref{lemma: polylog a}.
	
	Now  assume that \( r < r_1(n) \) or equivalently, that \( -(1-r)/r < r_1^{(n)}. \)  Then either \( -(1-r)/r < r_1^{(2)} \) or there is \( j \in \{ 3,4, \dots, n\} \)  such that \( -(1-r)/r \in  ( r_2^{(j)},r_1^{(j)}) .\) 
	Since for all \( j \in \{ 2,3, \dots ,n \},\)  all roots of \( \Li_{1-j} \) are simple roots and \( \Li_{1-j}(z) < 0 \) for all \( z \in (r_1^{(n)},0), \) it follows that  \( \Li_{1-j}(z)>0 \) if either \( j=2 \) and \( z \in (-\infty, r_1^{(2)} ), \) or if  \( j \geq 3 \) and \( z \in ( r_2^{(j)},r_1^{(j)}) .\)  
	Combining the above observations, we obtain~\ref{lemma: polylog b}.
\end{proof}

\section{Properties of the signed representation}\label{section: signed connected}

In this section, we state and prove two results that will be very useful in subsequent sections. Both of these results are more general versions of results that first appeared in~\cite{fgs}.

The first result of this section, Lemma~\ref{lemma: finite to infinite new} below, extends Lemma~\ref{lemma: finite to infinite} to signed measures in a finite setting.

\begin{lemma}\label{lemma: finite to infinite new}
	Let \( S \) be a finite set, let $X= (X_s)_{s\in S}$  be a \( \{ 0,1 \}\)-valued process, and let \( \nu \) be the unique signed measure which satisfies~\eqref{eq: step 1}. Let \( B \subseteq S. \) Then the following holds.
\begin{enumerate}[label=(\alph*)]

	\item For every measurable subset \( \mathcal{A} \subseteq \mathcal{P}(B)\smallsetminus \{ \emptyset \}, \) let \( \nu_B(\mathcal{A}) \coloneqq \nu\bigl( \{ A' \in \mathcal{P}(S) \colon  A' \cap B\in \mathcal{A} \} \bigr).\)  Then \( \nu_B\)  is the unique signed measure corresponding to \( X|_B. \)\label{item: finite to infinite i new}

    %\item If \( X=X^\nu \in \mathcal{R}\)  and \( B \subseteq S, \) then there is \( \nu_B \) such that \( X|_{B} = X^{\nu_B}.\) Moreover, for any non-empty measurable subset \( \mathcal{A} \subseteq \mathcal{P}(B), \) we have \( \nu_B(\mathcal{A}) = \nu\bigl( \{ A' \in \mathcal{P}(S) \colon  A' \cap B\in \mathcal{A} \} \bigr).\)\label{item: finite to infinite i}
	
	\item For every measurable subset \( \mathcal{A} \subseteq \mathcal{P}(B)\smallsetminus \{ \emptyset \}, \) let \( \nu_{B,0} (\mathcal{A}) \coloneqq \nu|_{\mathcal{P}(B)}(\mathcal{A}). \) Then  \( \nu_{B,0}\) is the unique signed measure corresponding to \(  X| \bigl\{ X(B^c)\equiv 0 \bigr\}. \)  \label{item: finite to infinite iii new}
	
    %\item If \( X=X^\nu \in \mathcal{R}\) and \( B \subseteq S,\) then there is a measure \( \nu_{B,0} \) on \( \mathcal{P}(B)\smallsetminus \{ \emptyset \} \) such that \( X| \bigl\{ X(B^c)\equiv 0 \bigr\} = X^{\nu_{B,0}}.\) Moreover, \( \nu_{B,0} = \nu|_{\mathcal{P}(B)}. \) \label{item: finite to infinite iii}
    
    %\item\label{item: finite to infinite ii} If there exist $S_1\subseteq S_2\subseteq,\ldots$ such that $S=\bigcup_i S_i$ and  $X_{S_n} = X^{\nu_n} $ for some $\nu_n \geq 0$, then $X=X^\nu$ for some $\nu \geq 0$. (The projection of $\nu$ on to each $S_n$ will simply be  $\nu_n$.)
\end{enumerate} 
\end{lemma}

\begin{proof}

	For \( I \subseteq B, \) let \( \mathcal{S}_I^{\cup ,B} \coloneqq \mathcal{S}_I^{\cup} \cap \mathcal{P}(B).\)

	We now show that~\ref{item: finite to infinite i new} holds. To this end, let \( \nu_0 \) be the unique signed measure which corresponds to \( X|_B, \) and let \( I \subseteq B. \) Then, by Lemma~\ref{lemma: unique signed}, we have  
	\begin{equation*} 
        \begin{split}
        &\nu_0(\mathcal{S}_I^\cup) 
        = -\log P\bigl(X|_B(I)\equiv 0\bigr)
        = -\log P\bigl(X(I)\equiv 0 \bigr) = \nu(\mathcal{S}_I^\cup).
        \end{split}
        \end{equation*}
   At the same time, by definition, we have
        \begin{align*}
        	\nu_B(\mathcal{S}_I^{\cup,B}) = \nu\bigl( \{ A' \in \mathcal{P}(S) \colon  A' \cap B\in \mathcal{S}_I^{\cup,B} \} \bigr)
        	=
        	\nu\bigl( \{ A' \in \mathcal{P}(S) \colon  A' \in \mathcal{S}_I^\cup \} \bigr) = \nu(  \mathcal{S}_I^{\cup} ).
        \end{align*}
        Again using Lemma~\ref{lemma: unique signed}, it follows that \( \nu_0 = \nu_B. \) This concludes the proof of~\ref{item: finite to infinite i new}.

        We now show that~\ref{item: finite to infinite iii new} holds. To this end, let \( \nu_0 \) be the unique signed measure which corresponds to \(  X| \bigl\{ X(B^c)\equiv 0 \bigr\} . \) Let \( I \subseteq B. \) Then
        \begin{align*} 
        \nu_0\bigl(\mathcal{S}_{I}^{\cup,B}\bigr) = -\log P\bigl( X(I) \equiv 0 |  X(B^c)\equiv 0  \bigr) 
        \end{align*}
        and
        \begin{align*}
        	\nu_{B,0}\bigl(\mathcal{S}_{I}^{\cup,B}\bigr) 
        	&= \nu|_{\mathcal{P}(B)}\bigl(\mathcal{S}_{I}^{\cup,B}\bigr) 
        	= \nu \bigl(\mathcal{S}_{I}^{\cup,B}\bigr) 
        	= \nu \bigl(\mathcal{S}_{I}^{\cup} \smallsetminus \mathcal{S}_{B^c}^{\cup} \bigr) 
        	= \nu \bigl(\mathcal{S}_{I \cup B^c}^{\cup} \bigr) -  \nu \bigl(  \mathcal{S}_{B^c}^{\cup} \bigr) 
        	\\&=-\log P\bigl( X(I \cup B^c) \equiv 0  \bigr) - \Bigl( -\log P\bigl( X( B^c) \equiv 0  \bigr)\Bigr)
        	 \\&= -\log P\bigl( X(I) \equiv 0 \mid  X(B^c)\equiv 0  \bigr) .
        \end{align*}
        Using Lemma~\ref{lemma: unique signed}, it follows that \( \nu_0 = \nu_{B,0}. \) This concludes the proof of~\ref{item: finite to infinite iii new}.
\end{proof}

In this paper, we mostly consider tree-indexed Markov chains, which naturally satisfy a Markov property. More generally, one could consider a Markov field on a graph. Assuming the graph is finite, let \( \nu \) be  the unique signed measure which satisfies~\eqref{eq: step 1}.
By~\cite[Proposition 3.7]{fgs}, if \( \nu \geq 0, \) then \( \nu(S) = 0 \) whenever \( S \subseteq V(G)\) is a disconnected set. The following result shows that this holds even if we drop the assumption that \( \nu \geq 0. \)

\begin{proposition} \label{prop: connected}
    Let \( X \) be a \(\{0,1\}\)-valued process on a finite connected graph \( G \) that satisfies the Markov property, and assume that for all finite \( A \subseteq V(G)\) we have
    \begin{equation*}
        P\bigl(X(A)\equiv 0 \mid  X(\mathcal{B}^+( A)) \equiv  0 \bigr) >0.
    \end{equation*}
    Let \( \nu\) be the unique signed measure that satisfies~\eqref{eq: step 1}, and let \( D\) be any disconnected subset of \( V(G).\)  Then \( \nu(D)=0.\)  
\end{proposition}

\begin{proof}%[Proof of Proposition~\ref{prop: connected}]
    We will first show that for any non-empty sets \( A_1   \subseteq V(G)\) and \( A_2 \subseteq V(G) \smallsetminus (A_1 \cup \mathcal{B}^+( A_1)),\) we have 
    \begin{equation}\label{eq: connected part 1}
        \nu\bigl( (\mathcal{S}_{A_1}^\cup \cap \mathcal{S}_{A_2}^\cup) \smallsetminus \mathcal{S}_{V(G)\smallsetminus (A_1 \cup A_2)}^\cup \bigr) =0.
    \end{equation}
    To this end, let \( A_1 \subseteq V(G)\) and \( A_2 \subseteq V(G)\smallsetminus (A_1 \cup \mathcal{B}^+( A_1))\) be two non-empty sets,  and let \( B \coloneqq (A_1 \cup A_2)^c = V(G)\smallsetminus (A_1 \cup A_2).\) 
    Since \( X\) satisfies the Markov property, we have 
    \begin{align*}
        &P\bigl(X(A_1)\equiv 0 \mid  X(\mathcal{B}^+( A_1)) \equiv  0 \bigr) =P\bigl(X(A_1)\equiv 0 \mid  X(A_1^c) \equiv  0 \bigr) =   P\bigl(X(A_1)\equiv 0 \mid  X(B) \equiv  0 \bigr).
    \end{align*}
    Since
    \begin{align*}
        & P\bigl(X(A_1)\equiv 0 \mid  X(B) \equiv  0 \bigr) =   e^{-\nu( \mathcal{S}_{A_1}^\cup \smallsetminus 
  \mathcal{S}_{B}^\cup )}  
    \end{align*} 
    and
    \begin{align*}
        & P\bigl(X(A_1)\equiv 0 \mid  X(\mathcal{B}^+( A_1)) \equiv  0 \bigr) = e^{-\nu( \mathcal{S}_{A_1}^\cup \smallsetminus \mathcal{S}_{A_1^c}^\cup)},
    \end{align*}
    and, by assumption, we have
    \begin{align*}
    	0<P\bigl(X(A_1)\equiv 0 \mid  X(\mathcal{B}^+( A_1)) \equiv  0 \bigr) \leq 1,
    \end{align*}
    it follows that 
    \begin{align*}
         0< e^{-\nu( \mathcal{S}_{A_1}^\cup \smallsetminus 
  \mathcal{S}_{B}^\cup )} = 
        e^{-\nu( \mathcal{S}_{A_1}^\cup \smallsetminus \mathcal{S}_{A_1^c}^\cup)} \leq 1,
    \end{align*}
    and hence
    \[
        0 \leq  \nu( \mathcal{S}_{A_1}^\cup \smallsetminus 
         \mathcal{S}_{B}^\cup)
        =
        \nu(\mathcal{S}_{A_1}^\cup \smallsetminus \mathcal{S}_{{A_1}^c}^\cup) <\infty.
    \]
    Since \( \partial A_1 \subseteq B, \) this yields
    \begin{equation*}
        \nu\bigl( (\mathcal{S}_{A_1}^\cup \cap \mathcal{S}_{A_2}^\cup) \smallsetminus \mathcal{S}_{B}^\cup \bigr)  = \nu\bigl( (\mathcal{S}_{A_1}^\cup \cap \mathcal{S}_{A_1^c}^\cup) \smallsetminus \mathcal{S}_{B}^\cup \bigr) 
        =
        0
    \end{equation*}
    and thus completes the proof of~\eqref{eq: connected part 1}.

    We now complete the proof of the proposition by showing that for any non-empty sets \( A_1 \subseteq V(G)\) and  \( A_2 \subseteq V(G)\smallsetminus (A_1 \cup \mathcal{B}^+( A_1)),\) we have \( \nu(A_1 \cup A_2) = 0.\)
    %
    %To this end, note first that since \( A_1 \subseteq V(G)\) and  \( A_2 \subseteq V(G)\smallsetminus (A_1 \cup \mathcal{B}^+( A_1))\) are both non-empty, then \( |A_1|,|A_2| \geq 1.\) 
    To this end, note first that if \( |A_1| = |A_2|=1,\) then, by~\eqref{eq: connected part 1}, we have
    \begin{equation*}
        \nu(A_1 \cup A_2 ) = \nu\bigl( (\mathcal{S}_{A_1}^\cup \cap \mathcal{S}_{A_2}^\cup) \smallsetminus \mathcal{S}_{B}^\cup \bigr) = 0.
    \end{equation*}
    For the more general statement, we note that
    \begin{align*}
        &0=\nu\bigl( (\mathcal{S}_{A_1}^\cup \cap \mathcal{S}_{A_2}^\cup) \smallsetminus \mathcal{S}_{B}^\cup \bigr)  = \sum_{\substack{S_1 \subseteq A_1,\, S_2 \subseteq A_2 \colon  S_1,S_2 \neq \emptyset}} \nu(S_1 \cup S_2)
         \\&\qquad= \nu(A_1 \cup A_2) + \sum_{\substack{S_1 \subseteq A_1,\, S_2 \subseteq A_2 \colon S_1,S_2 \neq \emptyset ,\\ |S_1|+|S_2| < |A_1|+|A_2|}} \nu(S_1 \cup S_2).
    \end{align*}
    Using induction, the desired conclusion immediately follows.
\end{proof}

\section{Formulas for tree-indexed Markov chains}\label{section: markov formulas}

By~\cite[Theorem~3.1]{fgs}, the formula in~\eqref{eq: Mobius inversion result} can be substantially simplified when the process \( X \) is a Markov chain. The proof relied on using the corresponding probability measure and several algebraic manipulations. 
The following lemma extends this result from Markov chains to Markov processes on graphs. When the graph is a line, the corresponding process is a regular Markov chain, and in this case we recover~\cite[Theorem~3.1]{fgs}. However, even in this special case, the proof presented below is substantially easier than the corresponding proof in~\cite{fgs}.
\begin{lemma}\label{lemma: graphs}
	    Let \( X \) be a \(\{0,1\}\)-valued process on a finite  graph \( G \) that satisfies the Markov property, and assume that for all sets \( A \subseteq V(G)\) we have
    \begin{equation*}
        P\bigl(X(A)\equiv 0 \mid  X(\mathcal{B}^+( A)) \equiv  0 \bigr) >0.
    \end{equation*}
    Let \( \nu\) be the unique signed measure corresponding to \( X ,\) and let \( S \subseteq V(G) \) be non-empty. Then
    \begin{equation}\label{eq: nu eq new}
    	\nu(S) = \begin{cases}
    		\sum_{J \subseteq \mathcal{B}^-(S)} (-1)^{|J|} 
         \log P\pigl( X\bigl(J \cup \mathcal{B}^+(S) \bigr) \equiv 0 \pigr) &\text{if } S \text{ is connected}\cr
         0 &\text{otherwise.}
    	\end{cases}
    \end{equation}
\end{lemma}

\begin{remarks}
	The formula in~\eqref{eq: nu eq new} is valid also for finite subsets of the vertex set of infinite graphs by restriction (using Lemma~\ref{lemma: finite to infinite}). However, we stress that in this case, we cannot conclude that the formula in~\eqref{eq: nu eq new} holds for infinite sets.
\end{remarks}

The main motivation for obtaining a lemma such as Lemma~\ref{lemma: graphs} is that for a graph, the number of vertices in the boundary of a connected subset is sometimes substantially smaller than the number of vertices in the set, and in this case the number of terms in~\eqref{eq: nu eq} is much lower than the number of terms in~\eqref{eq: Mobius inversion result}. This is, in particular, true when the graph is a line graph, in which case the boundary of any connected set has a cardinality of at most two.

\begin{proof}[Proof of Lemma~\ref{lemma: graphs}]
	If \( S \) is not connected, then, by Proposition~\ref{prop: connected}, we have \( \nu(S) = 0.\)     
    Now assume that \( S  \) is connected. 
    Then, since \( \nu \) has support only on connected sets by Proposition~\ref{prop: connected}, we have
    \begin{equation*}
        \nu(S) 
        = 
        \nu(\mathcal{S}_{\mathcal{B}^-(S)}^\cap \smallsetminus \mathcal{S}_{\mathcal{B}^+(S)}^\cup) .
    \end{equation*}
    Let \( \nu_{{\mathcal{B}^+(S) \cup \mathcal{B}^-(S)}} \) be the unique signed measure which corresponds to  \( X|_{\mathcal{B}^+(S) \cup \mathcal{B}^-(S)}.\) Then, by Lemma~\ref{lemma: finite to infinite new}\ref{item: finite to infinite i new}, we have 
    \begin{equation*}
        \nu(\mathcal{S}_{\mathcal{B}^-(S)}^\cap \smallsetminus \mathcal{S}_{\mathcal{B}^+(S)}^\cup)  
        =
        \nu_{\mathcal{B}^+(S) \cup \mathcal{B}^-(S)}(\mathcal{S}_{\mathcal{B}^-(S)}^\cap \smallsetminus \mathcal{S}_{\mathcal{B}^+(S)}^\cup) 
        =
        \nu_{\mathcal{B}^+(S) \cup \mathcal{B}^-(S)}(\mathcal{B}^-(S)).
    \end{equation*}
    Combining the previous equations and applying~Lemma~\ref{lemma: unique signed}, we obtain
    \begin{align*}
        \nu(S) &=\nu_{\mathcal{B}^+(S) \cup \mathcal{B}^-(S)} \bigl(\mathcal{B}^-(S) \bigr)
        =
        \sum_{I \subseteq \mathcal{B}^-(S)} (-1)^{|\mathcal{B}^-(S)|-|I|} 
         \log P\bigl( X((\mathcal{B}^-(S) \cup \mathcal{B}^+(S))\smallsetminus I) \equiv 0 \bigr)
         \\&=  
        \sum_{J \subseteq \mathcal{B}^-(S)} (-1)^{|J|} 
         \log P\bigl( X(J \cup \mathcal{B}^+(S) ) \equiv 0 \bigr),
    \end{align*}
    which is the desired conclusion.
\end{proof}

In the case of tree-indexed Markov chains, the formula in~\eqref{eq: nu eq new} can be simplified further, and this is the content of the next lemma.
    \begin{lemma}\label{lemma: trees}
    Let \( T\) be a finite tree, and let \( X\) be a tree-indexed Markov chain on \( T\) with parameters \( ((r_v),(p_e)). \)
    %. For \( v \in V(T) \) and \( e \in E(T), \) let \( r_v \in (0,1] \) and \( p_e \in [0,1], \) and let \( X\) be a tree-indexed Markov chain on \( T\) with parameters \( ((r_v),(p_e)). \) 
    %
    Let \( \nu \) be the unique signed measure corresponding to \( X \) in the sense of Lemma~\ref{lemma: unique signed}. Then, for any non-empty \( S \subseteq V(T),\) we have
    \begin{equation}\label{eq: nu eq} 
        \nu(S) = 
        \begin{cases} 
         \sum_{J \subseteq \mathcal{B}^-(S)} (-1)^{|J|} 
         \log P\pigl( X \bigl(J \cup (\mathcal{B}^+(S)\smallsetminus \partial_S J) \bigr) \equiv 0 \pigr)
          &\text{if ${S}$ is connected and } |S| \geq 2\cr 
          -\log P(X(S \cup \mathcal{B}^+(S))\equiv 0 \mid X \bigl( \mathcal{B}^+(S) \bigr) \equiv 0 ) 
            %\log \frac{\bigl( r_v \prod_{v' \sim v }(1-p_{(v,v')}+p_{(v,v')}r_{v'}) + (1-r_v) \prod_{v' \sim v}(p_{(v,v')} r_{v'})\bigr) }{r_v\prod_{v' \sim v}(1-p_{(v,v')}+p_{(v,v')}r_{v'})} 
         &\text{if } S = \{ v \} \cr
            0 &\text{otherwise.}
        \end{cases} 
    \end{equation} 
\end{lemma}

\begin{proof} 
    Let \( S\) be a non-empty subset of \( V(T). \) 
    The desired conclusion will follow from~Lemma~\ref{lemma: graphs} if we can show that for tree-indexed Markov chains, the formula in~\eqref{eq: nu eq new} simplifies to the formula in~\eqref{eq: nu eq}. To this end, note first that if \( S \) is not connected, then the two formulas trivially agree. 
    Next, assume that \( S = \{v \} \) for some \( v \in V(T). \) Then \( \mathcal{B}^-(S) = S = \{ v \}, \) and hence, by~Lemma~\ref{lemma: graphs}, we have 
    \begin{equation*} 
    	\nu(S) = 
    		\log P\pigl( X\bigl( \mathcal{B}^+(S) \bigr) \equiv 0 \pigr)  
    		-
    		\log P\pigl( X\bigl(S  \cup \mathcal{B}^+(S) \bigr) \equiv 0 \pigr)  
    \end{equation*}
    which immediately implies~\eqref{eq: nu eq} in this case.
    Now assume that \( S \) is connected and \( |S| \geq 2.\) Then, for any \( J \subseteq \mathcal{B}^-(S), \) we have 
    \begin{align*}
    	P\pigl( X\bigl(J \cup \mathcal{B}^+(S) \bigr) \equiv 0 \pigr)
    	=
    	P\bigl( X(J \cup (\mathcal{B}^+(S)\smallsetminus \partial_S J ) ) \equiv 0 \bigr)  \prod_{ \substack{v \in J ,\,  v' \in \partial_S v}} %v' \in N(\{ v \}) \cap \mathcal{B}^+(S) }} 
    	(1-p_{(v,v')}+p_{(v,v')}r_{v'}),      
	\end{align*}
	and hence
    \begin{align*}
    	&\sum_{J \subseteq \mathcal{B}^-(S)} (-1)^{|J|} 
         \log P\pigl( X\bigl(J \cup \mathcal{B}^+(S) \bigr) \equiv 0 \pigr) 
        \\&\qquad=
    	\sum_{J \subseteq \mathcal{B}^-(S)} (-1)^{|J|} 
         \biggl( \log P\pigl( X \bigl(J \cup (\mathcal{B}^+(S)\smallsetminus \partial_S J) \bigr) \equiv 0 \pigr) 
         + 
    	\sum_{ \substack{v \in J ,\, v' \in  \partial_S v}} %v' \in N(\{ v \}) \cap \mathcal{B}^+(S) }} 
    	\log (1-p_{(v,v')}+p_{(v,v')}r_{v'})\biggr) .   
    \end{align*}
    Now note that
    \begin{align*}
    	&\sum_{J \subseteq \mathcal{B}^-(S)} (-1)^{|J|} \sum_{ \substack{v \in J ,\,  v' \in \partial_S v}}% \in N(\{ v \}) \cap \mathcal{B}^+(S) }} 
    	\log (1-p_{(v,v')}+p_{(v,v')}r_{v'}) 
    	\\&\qquad
    	=\sum_{ \substack{v \in \mathcal{B}^-(S) ,\, v' \in \partial_S v}} %v' \in N(\{ v \}) \cap \mathcal{B}^+(S) }}  
    	\log (1-p_{(v,v')}+p_{(v,v')}r_{v'}) 
    	\sum_{J \subseteq \mathcal{B}^-(S) \colon v \in J} (-1)^{|J|} .
    \end{align*}
    Since \( |S|\geq 2 \) implies that  \( |\mathcal{B}^-(S)|\geq 2,\) we have \( \sum_{J \subseteq \mathcal{B}^-(S) \colon v \in J} (-1)^{|J|} = 0,\) and hence 
    \begin{equation*}
    	\sum_{J \subseteq \mathcal{B}^-(S)} (-1)^{|J|} \sum_{ \substack{v \in J ,\,  v' \in \partial_S v }} \log (1-p_{(v,v')}+p_{(v,v')}r_{v'})  =0,
    \end{equation*}
    implying in particular that~\eqref{eq: nu eq} holds in this case. 
    This concludes the proof.
\end{proof}

\begin{remark}
	\cite[Theorem~3.1]{fgs}, which in particular gives a simple formula for \( \nu \) when \( X \) is a Markov chain,  holds not only for Markov chains but also for renewal processes. Since it is not obvious how to define renewal processes on trees, Lemma~\ref{lemma: trees} does not mention renewal processes, even though in the special case that the tree is a line, the proof of Lemma~\ref{lemma: trees} gives the more general result of~\cite[Theorem~3.1]{fgs}.
\end{remark}

\section{Scaling properties of trees}\label{section: scaling}

In this section, we first define what we mean by a tree being self-similar, and then show how this property relates to a Markov process indexed by such a tree being Poisson representable.

\begin{definition}\label{def: self-similar}
	We say that an infinite tree \( T \) is \emph{self-similar} if there are arbitrarily large \( k \geq 2 \) such that if we let \(  T^k \) be the tree obtained by replacing each edge of \( T \) with a path of length \( k, \) then \(  T^k \) is a subgraph of \( T. \) 
\end{definition}

We note that with this definition, both line graphs, octopus trees, and regular trees are self-similar.

The main result in this section is the following proposition.

\begin{proposition}\label{proposition: scaling}
    Let \( T \) be a  self-similar infinite tree, and for \( r,p \in (0,1), \) let \( X_{r,p} \) be a tree-indexed Markov chain on \( T \) with parameters \( (r,p). \) Then the following holds.
    \begin{enumerate}[label=(\alph*)]
        \item If there is \( \varepsilon, r \in (0,1) \) such that \( X_{r,p} \notin \mathcal{R}\) for all  \( p \in (1-\varepsilon,1), \)  then \( X_{r,p} \notin \mathcal{R} \) for all \( p \in (0,1).\)\label{lemma: scaling a}
        
        \item If there is \( \varepsilon, r \in (0,1) \) such that  \( X_{r,p} \in \mathcal{R}\) for all  \( p \in (0,\varepsilon), \)  then \( X_{r,p} \in \mathcal{R} \) for all \( p \in (0,1).\)\label{lemma: scaling b}
    \end{enumerate}
     
\end{proposition}

\begin{proof}
	Let \( r,p \in (0,1) \) and let \( X = X_{r,p}. \) 
	
	For \( k \geq 1, \) let \( T_k \) be as in Definition~\ref{def: self-similar}, and let  \( V_0(T^k) \) denote the vertices in \( T^k \) that are the endpoints of the added paths in Definition~\ref{def: self-similar}.
	Note that if \( k \geq 1 \) is such that \( T^k \) is a subgraph of \( T,\) then \( X(V_0(T^k)) \) is equal in distribution to a tree-indexed Markov chain \( X^{(k)} \) on \( T \) with parameters \( (r, 1-(1-p)^k) .\) Using Lemma~\ref{lemma: finite to infinite new}\ref{item: finite to infinite i new}, it follows that if \( X \in \mathcal{R},\) then \( X^{(k)} \in \mathcal{R},\) and if \( X^{(k)} \notin \mathcal{R},\) then \( X \notin \mathcal{R}.\) 
	
	We now show that~\ref{lemma: scaling a} holds.
	To this end, assume that \( r \) and \( \varepsilon \in (0,1)\) are such that  such that \( X_{r,p'} \notin \mathcal{R}\) for all \( p' \in (1-\varepsilon,1). \) Further, let \( k \geq 1 \) be such that \( T^k \) is a subgraph of \( T \) and \( 1-(1-p)^k>1-\varepsilon.\) Then, by assumption,  we have \( X^{(k)} \notin \mathcal{R}.\) Since \( X^{(k)}\) is equal in distribution to \( X(V_0(T^k))\) which is a restriction of \( X,\) using Lemma~\ref{lemma: finite to infinite new}\ref{item: finite to infinite i new}, it follows that \( X \notin \mathcal{R}.\)
	
	We now show that~\ref{lemma: scaling b} holds. To this end, instead assume that \( r \) and \( \varepsilon \in (0,1) \) are such that  \( X_{r,p'} \in \mathcal{R}\) for all  \( p' \in (0,\varepsilon). \) Further, let \( k \geq 1\) be such that \( T^k \) is a subgraph of \( T \) and \( 1-(1-p)^{1/k}<\varepsilon.\) 
	Let \( \hat X \) be a tree-indexed Markov chain on \( T \) with parameters \( (r,1-(1-p)^{1/k}).\) By the choice of \( k,\) we have \( \hat X \in \mathcal{R}.\) Since \( X \) is equal in distribution to \( \hat X^{(k)} \), using Lemma~\ref{lemma: finite to infinite new}\ref{item: finite to infinite i new}, it follows that \( X \in \mathcal{R}.\)
\end{proof}

\section{Finite trees and \( p \) close to zero}\label{section: finite and p zero}

In this section, we provide proofs of parts~\ref{theorem: finite trees ii} and~\ref{theorem: finite trees iii} of Theorem~\ref{theorem: finite trees}. The main tool in these proofs is the following Proposition~\ref{proposition: new negative lemma for p close to zero} below.
 In this proposition, and throughout the rest of this paper, we will use the following notation. For a multiset \( E \subseteq E(T), \) where \( e \in E \) has multiplicity \( \mu(e), \)  we write \( \frac{d}{dp_E} \) to denote the partial derivative \( \prod_{e \in E}\frac{d^{\mu(e)}}{dp_e^{\mu(e)}}.  \) Analogously, for a multiset \( V \subseteq V(T), \) where \( v \in V \) has multiplicity \( \mu(v), \)  we write \( \frac{d}{dr_V} \) to denote the partial derivative \( \prod_{v \in V}\frac{d^{\mu(v)}}{dp_v^{\mu(v)}}.  \)

\begin{proposition}\label{proposition: new negative lemma for p close to zero}
	Let \( T \) be a finite tree, and let \( S \subseteq V(T)\) be connected and have a cardinality of at least two. 
	Let \( r \in (0,1), \) and for \( e \in E(T), \) let \( p_e \in (0,1).\) Further, let \( X \) be the \( (r,(p_e)) \)-tree indexed Markov chain on \( T, \) and let \( \nu \) be the signed measure corresponding to \( X \) in the sense of Lemma~\ref{lemma: unique signed}. 
	Then, the following holds. 
	\begin{enumerate}[label=(\alph*)]
		\item Let \( E \subseteq E(T) \) be a multiset with \( |E| < |\mathcal{B}^+( S)|. \) Then  
		\( \bigl[ \frac{d}{dp_E} \nu(S) \bigr]_{(p_e) \equiv 0} =0. \) \label{item: new tree lemma 2 2}
		\item Let \( E \subseteq E(T) \) be a multiset with \( |E| = |\mathcal{B}^+( S)|. \) Let \( E_S \) be the set of edges between \( S \) and \( S^c. \) Then
		\begin{equation*}
			\Bigl[ \frac{d}{dp_E} \nu(S) \Bigr]_{(p_e) \equiv 0} =
			\begin{cases}
				(1-r) r^{|\mathcal{B}^+(S)|-1}-(-1)^{|\mathcal{B}^+(S)|}  \tilde B_{|\mathcal{B}^+(S)|} (1-r)^{|\mathcal{B}^+(S)|}  &\text{if } E=E_S \cr 
				0 &\text{otherwise,}
			\end{cases}
		\end{equation*}  
	where we recall that for \( n \geq 0 \), \( \tilde B_n \) is the \( n \)th complementary Bell number.
		\label{item: new tree lemma 3 2} 
	\end{enumerate}  
\end{proposition}

Before we provide a proof of Proposition~\ref{proposition: new negative lemma for p close to zero}, we show how it implies~\ref{theorem: finite trees ii} and~\ref{theorem: finite trees iii} of Theorem~\ref{theorem: finite trees}.

\begin{proof}[Proof of Theorem~\ref{theorem: finite trees}\ref{theorem: finite trees ii} and \ref{theorem: finite trees iii}]
	For \( j \geq 2 \) let
	\begin{align*}
		f_j(r) \coloneqq (1-r) r^{j-1}-(-1)^{j}  \tilde B_{j} (1-r)^{j}.
	\end{align*}
	Then \( f_j(1) = 0, \) \( f_j(r) > 0 \) for all \( r \in (0,1) \) sufficiently close to one, and \( f_j \) has at most two real roots. Moreover, if \( f_j \) has exactly two roots, then both are simple, with one root at \( r = 1 \) and a second root at 
	\begin{align*}
	  	\tilde r^{(j)} = \frac{\bigl( (-1)^{j} \tilde B_{k}\bigr)^{1/(j-1)}}{1+ \bigl( (-1)^{j} \tilde B_{j}\bigr)^{1/(j-1)}}.
	\end{align*}  
	
	Let \( \nu \) be the unique signed measure corresponding to \( X. \) For \( p = 0, \) one easily verifies that \( \nu(S) = -\log r \) if \( S= V(T) \) and \( \nu_T(S) = 0 \) else. Moreover, if \( S \) is not connected, then \( \nu(S) = 0 \) by~Proposition~\ref{prop: connected}.
	Next, let \( S \subseteq V(T) \) be connected with \( S \neq V(T). \) Then, using  a Taylor expansion, it follows from Proposition~\ref{proposition: new negative lemma for p close to zero} that if \( p \) is sufficiently small, then
	\begin{equation}\label{eq: sign eq} 
		\sgn \nu(S) = \sgn \bigl( (1-r) r^{|\mathcal{B}^+(S)|-1}-(-1)^{|\mathcal{B}^+(S)|}  \tilde B_{|\mathcal{B}^+(S)|} (1-r)^{|\mathcal{B}^+(S)|} \bigr) = \sgn f_{|\mathcal{B}^+(S)|}(r), 
	\end{equation} 
	i.e. then \(  \nu(S)  \) and \( f_{|\mathcal{B}^+(S)|}(r) \) have the same sign. 
	Since \( k \) is the size of the boundary of \( T, \) we have \( |\mathcal{B}^+(S)| \leq k \), and hence, if  \( r>r_0(k), \) it follows that \( \nu(S) > 0 \) for sufficienty small \( p \). Since this holds for all connected sets \( S \subseteq V(T), \) this concludes the proof of part~\ref{theorem: finite trees ii}.
	
 	To see that part~\ref{theorem: finite trees ii} holds, assume that \( r<r_0(k). \) Then there is \( j \leq k \) such that \( f_j(r) < 0 \) and a connected set \( S \subseteq V(T) \) such that \( \mathcal{B}^+(S) = j. \) Using~\eqref{eq: sign eq}, it follows that \( \nu(S)<0 \) for sufficiently small \( p , \) and hence for such \( p , \) we must have \( X \notin \mathcal{R}, \) thus concluding the proof of~\ref{theorem: finite trees iii}.
\end{proof}

We now proceed to the proof of Proposition~\ref{proposition: new negative lemma for p close to zero}.

\begin{proof}[Proof of Proposition~\ref{proposition: new negative lemma for p close to zero}]
	
	Let \( E_S \) be the set of edges between \( S \) and \( S^c. \)
	
	For \( J \subseteq \mathcal{B}^-(S), \) consider the event 
	\begin{equation*}
		\mathcal{E}_J \coloneqq \pigl\{ X \bigl(J \cup (\mathcal{B}^+(S)\smallsetminus \partial_S J) \bigr) \equiv 0 \pigr\}.
	\end{equation*}
	Since \( S \) is a connected set and \( |S|\geq 2, \) by Lemma~\ref{lemma: trees}, we have 
	\begin{equation} 
        \nu(S) =  
         \sum_{J \subseteq \mathcal{B}^-(S)} (-1)^{|J|} 
         \log P( \mathcal{E}_J) .
    \end{equation}  
    For any fixed  \( e' \in E_S\), let \( v \) denote its endpoint in \( \mathcal{B}^-(S) . \) Then, if
\( E \subseteq E(T) \smallsetminus \{ e' \}\) is a multiset, for any  \( J \subseteq \mathcal{B}^-(S) \) we have
    \begin{align*}
    	&\Bigl[ \frac{d}{dp_E} \log P( \mathcal{E}_J) \Bigr]_{p_{e'} = 0}
    	= 
    	\Bigl[\frac{d}{dp_E} \log P( \mathcal{E}_{J \Delta \{ v \}}) \Bigr]_{p_{e'} = 0}.
    \end{align*}
    Hence, for any multi-set \(  E \subseteq E(T) \) with \( E_S \not \subseteq E, \) we have 
    \begin{align*}
    	\Bigl[ \frac{d}{dp_E} \nu(S) \Bigr]_{(p_e) \equiv 0} = 0.
    \end{align*}
    This completes the proof of~\ref{item: new tree lemma 2 2}, and also the proof of~\ref{item: new tree lemma 3 2} in the case \( E \neq E_S. \)
    
    We now prove part~\ref{item: new tree lemma 3 2} in the case \( E = E_S \). 
     To this end, note first that for any non-empty \( J \subseteq \mathcal{B}^-(S) \), the event \( \mathcal{E}_J \) does not depend on all edges in \( E_S, \) and hence 
    \begin{align*}
    	\frac{d}{dp_{E_S}} \log P ( \mathcal{E}_J  ) =0.
    \end{align*}
    Hence, 
    \begin{equation*} 
        \frac{d}{dp_{E_S}}\nu(S) =   \frac{d}{dp_{E_S}} \log P(\mathcal{E}_\emptyset ) =  
         \frac{d}{dp_{E_S}} \log P\bigl( X \bigl( \mathcal{B}^+(S) ) \equiv 0 \bigr) .
    \end{equation*} 
    Let  \( \mathcal{P}(E_S)  \) denote the set of all partitions of \( E_S \) into non-empty sets. Then
	\begin{equation*}
		\begin{split}
		&\Bigl[ \frac{d}{dp_{E_S}}  
         \log P\bigl( X (\mathcal{B}^+(S) ) \equiv 0 \bigr) \Big]_{(p_e)\equiv  0} 
         \\&\qquad=
         \sum_{(E_1, \dots , E_k) \in \mathcal{P}(E_S)}
         \frac{ \bigl[ (-1)^{k-1}\prod_{j=1}^k \frac{d}{dp_{E_j}}P\bigl( X (\mathcal{B}^+(S) ) \equiv 0 \bigr)\bigr]_{(p_e) \equiv  0}}{ \bigl[ P\bigl( X (\mathcal{B}^+(S) ) \equiv 0 \bigr)^k \bigr]_{(p_e) \equiv 0}}  .
		\end{split}
	\end{equation*}
	We now make a few additional observations about the quantities on the right-hand side of the previous equation. First, note that
	\[
		\Bigl| P\bigl( X (\mathcal{B}^+(S) ) \equiv 0 \bigr) \Bigr]_{(p_e) \equiv 0} = r.
	\]    
	Next, note that for any multi-set \( E' \subset E(T) \) that contains two copies of some edge, we have
	\begin{align*}
		\frac{d}{dp_{E'}}P\bigl( X (\mathcal{B}^+(S) ) \equiv 0 \bigr) = 0.
	\end{align*}
	We now consider the case that \( E' \subseteq E_S\) is not a multiset. For this case, note first that since \( E(T_S) \) and \(  E_S  \) are disjoint, we have 
	\begin{align*}
		\Bigl[ \frac{d}{dp_{E'}} P\bigl( X (\mathcal{B}^+(S) ) \equiv 0 \bigr)\Bigr]_{(p_e)\equiv  0} =  
		\Bigl[ \frac{d}{dp_{E'}} \pigl[ P\bigl( X (\mathcal{B}^+(S) ) \equiv 0 \bigr) \pigr]_{(p_e)_{e \in E(T_S)}\equiv 0} \Bigr]_{(p_e)\equiv  0}.
	\end{align*}
	Moreover, we have 
	\begin{align*}
		\pigl[ P\bigl( X (\mathcal{B}^+(S) ) \equiv 0 \bigr) \pigr]_{(p_e)_{e \in E(T_S)}\equiv 0}
		=
		r \prod_{e \in E_S} (1-p_e+p_er) + (1-r) \prod_{e \in E_S} p_e r 
	\end{align*}
	Using this observation, it follows that if  \( E' \subset E_S \) is not a multiset, then
	\begin{align*}
		\Bigl[ \frac{d}{dp_{E'}} P\bigl( X (\mathcal{B}^+(S) ) \equiv 0 \bigr)\Bigr]_{(p_e)\equiv  0} %= \Bigl[r \frac{d}{dp_{E'}} \prod_{e \in E'}  (1-p_e+p_er) \Bigr]_{(p_e)\equiv 0}
		= 
		r(-1+r)^{|E'|}.
	\end{align*} 
	and similarly, if \( E' = E_S, \) then
	\begin{align*}
		&\Bigl[ \frac{d}{dp_{E'}}P\bigl( X (\mathcal{B}^+(S) ) \equiv 0 \bigr)\Bigr]_{(p_e) \equiv  0} 
		%= 
		%\Bigl[ \frac{d}{dp_{E_S}} \Bigl( r \prod_{e \in E_S} (1-p_e+p_er) + (1-r) \prod_{e \in E_S} rp_e \Bigr)\Bigr]_{(p_e)\equiv 0}
		 %\\&\qquad
		 =
		 r (-1+r)^{|E_S|} + (1-r) r^{|E_S|} .
	\end{align*} 
	Combining the above equations and simplifying, it follows that 
    \begin{align*} 
        & \pigl[\frac{d}{dp_{E_S}}\nu(S) \pigr]_{(p_e)\equiv 0} =  
         -(-1+r)^{|\mathcal{B}^+(S)|} \!\!\!\!\!\!\!\! \sum_{(E_1, \dots , E_k) \in \mathcal{P}(E_S)} \!\!\!\!\!\!(-1)^{k}  
         + (1-r) r^{|\mathcal{B}^+(S)|-1} .
    \end{align*}  
    This completes the proof of~\ref{item: new tree lemma 3 2}.
\end{proof}

\section{Finite trees and \( p \) close to one}\label{section: finite and p one}

The main purpose of this section is to provide proofs of Theorem~\ref{theorem: finite trees}\ref{theorem: finite trees iv} and Theorem~\ref{theorem: finite trees}\ref{theorem: finite trees i}.
The main tool of the proofs will be the following proposition.

\begin{proposition}\label{proposition: new negative lemma for p close to one}
	Let \( T \) be a finite tree, and let \( S \) be a connected subset of \( V(T)\) with a cardinality of at least two.  
	Further, let \( r \in (0,1), \) for \( e \in E(T), \) let  \( p_e \in (0,1), \) and let \( X \) be the tree indexed Markov chain on \( T\) with parameters \( (r,(p_e)). \) Further,  let \( \nu \) be the signed measure corresponding to \( X \) in the sense of Lemma~\ref{lemma: unique signed}. 
	Then the following holds.
	\begin{enumerate}[label=(\alph*)] 
		\item Let \( E \subseteq E(T) \) be a non-empty multiset with \( E( T_S ) \nsubseteq E. \) Then  
		\( \bigl[ \frac{d}{dp_E} \nu(S) \bigr]_{(p_e) \equiv 1} =0. \)  \label{item: new tree lemma 2}
		\item Let \( E \subseteq E(T) \) be a non-empty multiset with \( E =  E( T_S ). \) Then \label{item: new tree lemma 3}
	\begin{equation}\label{eq: nu derivative}
		\Bigl[ \frac{d}{dp_{E}} \nu(S) \Bigr]_{p_{E(T)} \equiv 1}  =
			(-1)^{|\mathcal{B}^-(S) |+1} 
          \pigl( \frac{r}{1-r} \pigr)^{|V(T_S) \smallsetminus \mathcal{B}^-(S)|-1}\prod_{j=2}^\infty    \Li_{1-j}\bigl(-r^{-1}(1-r)\bigr)^{k_j}  , 
	\end{equation} 
	where \( k_j \) is the number of vertices of degree \( j \) in \( T_S. \)
	\end{enumerate} 
\end{proposition}

Before we provide a proof of Proposition~\ref{proposition: new negative lemma for p close to one}, we show how it implies~\ref{theorem: finite trees iv} and~\ref{theorem: finite trees i} of Theorem~\ref{theorem: finite trees}.

\begin{proof}[Proof of Theorem~\ref{theorem: finite trees}\ref{theorem: finite trees iv} and Theorem~\ref{theorem: finite trees}\ref{theorem: finite trees i}]
	Let \( S \subseteq V(T),\) and let \( \nu \) be the signed measure which corresponds to \( X \) as in Lemma~\ref{lemma: unique signed}.
	
	When \( (p_e) \equiv 1, \) \( X(v) \) and \( X(v') \) are independent for distinct \( v,v' \in V(T), \) implying in particular that \( X \in \mathcal{R} \) with 
	\[ \nu(S) = \begin{cases}
		-\log r  &\text{if } |S|=1 \cr 0 &\text{otherwise.}
	\end{cases} \]
	For general \( (p_e), \) we  know from Lemma~\ref{lemma: trees} that we have \( \nu(S) \geq 0 \) if \( |S| = 1 \) or if \( S \) is not connected.

	Now assume that \( |S|>1 \) is a connected set.
	Then, by Proposition~\ref{proposition: new negative lemma for p close to one} and a Taylor expansion, it follows that for all \( p \) sufficiently close to one, \( \nu(S) \) has the same sign as 	
	\begin{equation}\label{eq: sign equation} 
		 (-1)^{|E(T_S)|} (-1)^{|\mathcal{B}^-(S) |+1}  \prod_{j=2}^\infty  
         \Bigl( \mathrm{Li}_{1-j}(-(1-r)/r)  \Bigr)^{k_j}
         =
         \prod_{j=2}^\infty  
         \pigl( -\mathrm{Li}_{1-j}\bigl(-(1-r)/r \bigr) \pigr)^{k_j},
	\end{equation} 
	where for \( j \geq 2, \) \( k_j \) is the number of vertices in \( T_S \) with degree \( j. \)  
	Since by Lemma~\ref{lemma: polylog}\ref{lemma: polylog a}, \( -\Li_{1-j}(-(1-r)/r) > 0 \) for all \( j \in \{ 2,3 \dots, m\} \) if  \( r > r_1(m),\) this concludes the proof of~\ref{theorem: finite trees iv}.

	To see that part~\ref{theorem: finite trees i} holds, we recall from Lemma~\ref{lemma: polylog}\ref{lemma: polylog b} that if \( r < r_1(m), \) then there is \( m' \leq m \) such that \( -\Li_{1-m'}(-(1-r)/r)<0 . \) 
	Since \( m \) is the maximal degree of any vertex in \( T \), there is a set \( S \subseteq V(T) \) such that \( T_S \) is connected, \( T_S \) has exactly one vertex with degree \( m',\) and all other vertices have degree 1. 	For this set, since  \( -(1-r)/r \in (r_2^{(m')},r_1^{(m')}), \) it follows that~\eqref{eq: sign equation} is strictly negative, and hence the restriction of \( X \) to \( S \) is not in \( \mathcal{R} \) for \( (p_e) \) sufficiently close to one. Using Lemma~\ref{lemma: finite to infinite}\ref{item: finite to infinite i}, it follows that \( X \notin \mathcal{R} \) for \( ( p_e) \) sufficiently close to one, thus completing the proof of~part~\ref{theorem: finite trees i}.
\end{proof}

We now move to the proof of Proposition~\ref{proposition: new negative lemma for p close to one}, which will use a number of lemmas that we first state and prove. The first of these lemmas, Lemma~\ref{lemma: observation 1} below, shows that in the setting of Proposition~\ref{proposition: new negative lemma for p close to one}, if the set \( E \) is not connected, then the derivative \( \frac{d}{dp_E} \) of the terms appearing on the righ-hand side of~\eqref{eq: nu eq} of Lemma~\ref{lemma: trees} is zero.

\begin{figure}
\hfil
	\begin{subfigure}[t]{.3\textwidth}\centering
			\begin{tikzpicture}[scale=.4]

    	\foreach \i in {0,1}
    		\foreach \j in {-1,1}
    			\foreach \k in {-1,1}
    			{ 
    				\draw[] ({2*cos(120*\i+30*\j)},{2*sin(120*\i+30*\j)}) -- ({3*cos(120*\i+30*\j+15*\k)},{3*sin(120*\i+30*\j+15*\k )});
    				\filldraw[fill=white] ({3*cos(120*\i+30*\j+15*\k)},{3*sin(120*\i+30*\j+15*\k)}) circle (5pt);   
    			}
    	  
    	  \foreach \i in {1}
    		\foreach \j in {-1}
    			\foreach \k in {-1,1}
    			{ 
    				\draw[red, very thick] ({2*cos(120*\i+30*\j)},{2*sin(120*\i+30*\j)}) -- ({3*cos(120*\i+30*\j+15*\k)},{3*sin(120*\i+30*\j+15*\k )});
    				\filldraw[fill=white] ({3*cos(120*\i+30*\j+15*\k)},{3*sin(120*\i+30*\j+15*\k)}) circle (5pt);   
    			}

    	\foreach \i in {0,...,1} 
    		\foreach \j in {-1,1}
    		{
    			\draw ({cos(120*\i)},{sin(120*\i)}) -- ({2*cos(120*\i+30*\j)},{2*sin(120*\i+30*\j)});
    			\filldraw[fill=black] ({2*cos(120*\i+30*\j)},{2*sin(120*\i+30*\j)}) circle (6pt);   
    		}

    	\foreach \i in {2} 
    		\foreach \j in {-1,1}
    		{
    			\draw[red, very thick] ({cos(120*\i)},{sin(120*\i)}) -- ({2*cos(120*\i+30*\j)},{2*sin(120*\i+30*\j)});
    			\filldraw[fill=white] ({2*cos(120*\i+30*\j)},{2*sin(120*\i+30*\j)}) circle (6pt);   
    			%\draw[red] ({2*cos(120*\i+30*\j)},{2*sin(120*\i+30*\j)}) circle (11pt);   
    		}

    	\foreach \i in {0,...,3}
    	{	   
    		\draw (0,0) -- ({cos(120*\i)},{sin(120*\i)});
    		\filldraw[fill=black] ({cos(120*\i)},{sin(120*\i)}) circle (7pt);  
    	}

    	\filldraw[fill=black] (0,0) circle (8pt);  
		
		\draw (.2,1.2) node {\footnotesize \( e \)};
     \end{tikzpicture}

		\caption{The setting of Lemma~\ref{lemma: observation 1}}\label{figure: obs 1}
		\end{subfigure}
		\hfil
		\begin{subfigure}[t]{.3\textwidth}\centering
			\begin{tikzpicture}[scale=.4]

    	\foreach \i in {0}
    		\foreach \j in { -1,1}
    			\foreach \k in {-1,1}
    			{ 
    				\draw[] ({2*cos(120*\i+30*\j)},{2*sin(120*\i+30*\j)}) -- ({3*cos(120*\i+30*\j+15*\k)},{3*sin(120*\i+30*\j+15*\k )});
    				\filldraw[fill=white] ({3*cos(120*\i+30*\j+15*\k)},{3*sin(120*\i+30*\j+15*\k)}) circle (5pt);   
    			}

    	\foreach \i in {1}
    		\foreach \j in {-1,1}
    			\foreach \k in {-1,1}
    			{ 
    				\draw[] ({2*cos(120*\i+30*\j)},{2*sin(120*\i+30*\j)}) -- ({3*cos(120*\i+30*\j+15*\k)},{3*sin(120*\i+30*\j+15*\k )});
    				\filldraw[fill=white] ({3*cos(120*\i+30*\j+15*\k)},{3*sin(120*\i+30*\j+15*\k)}) circle (5pt);   
    			}

    	  \foreach \i in {1}
    		\foreach \j in {-1,1}
    			\foreach \k in {-1,1}
    			{ 
    				\draw  ({2*cos(120*\i+30*\j)},{2*sin(120*\i+30*\j)}) -- ({3*cos(120*\i+30*\j+15*\k)},{3*sin(120*\i+30*\j+15*\k )});
    				\filldraw[fill=white] ({3*cos(120*\i+30*\j+15*\k)},{3*sin(120*\i+30*\j+15*\k)}) circle (5pt);   
    			}

    	\foreach \i in {0 } 
    		\foreach \j in {-1,1}
    		{
    			\draw[] ({cos(120*\i)},{sin(120*\i)}) -- ({2*cos(120*\i+30*\j)},{2*sin(120*\i+30*\j)});
    			\filldraw[fill=black] ({2*cos(120*\i+30*\j)},{2*sin(120*\i+30*\j)}) circle (6pt);   
    		}

    		\foreach \i in {1} 
    			\foreach \j in {-1,1}
    		{
    			\draw ({cos(120*\i)},{sin(120*\i)}) -- ({2*cos(120*\i+30*\j)},{2*sin(120*\i+30*\j)});
    			\filldraw[fill=black] ({2*cos(120*\i+30*\j)},{2*sin(120*\i+30*\j)}) circle (6pt);   
    		}

    	\foreach \i in {2} 
    		\foreach \j in {-1,1}
    		{
    			\draw[red, very thick] ({cos(120*\i)},{sin(120*\i)}) -- ({2*cos(120*\i+30*\j)},{2*sin(120*\i+30*\j)});
    			\filldraw[fill=white] ({2*cos(120*\i+30*\j)},{2*sin(120*\i+30*\j)}) circle (6pt);   
    			%\draw[red] ({2*cos(120*\i+30*\j)},{2*sin(120*\i+30*\j)}) circle (11pt);   
    		}

    	\foreach \i in {1,2}
    	{	   
    		\draw[red, very thick] (0,0) -- ({cos(120*\i)},{sin(120*\i)});
    		\filldraw[fill=black] ({cos(120*\i)},{sin(120*\i)}) circle (7pt);  
    	}
    	
    	\foreach \i in {0}
    	{	   
    		\draw (0,0) -- ({cos(120*\i)},{sin(120*\i)});
    		\filldraw[fill=black] ({cos(120*\i)},{sin(120*\i)}) circle (7pt);  
    	}

    	\filldraw[fill=black] (0,0) circle (8pt);  
		
		\draw (1.3,1.5) node {\footnotesize \( j \)}; 
		
     \end{tikzpicture}

		\caption{The setting of Lemma~\ref{lemma: observation 3}}\label{figure: obs 3}
		\end{subfigure}
		\hfil
	\caption{The figures above illustrate the settings of Lemma~\ref{lemma: observation 1} %, Lemma~\ref{lemma: observation 2} 
	and Lemma~\ref{lemma: observation 3}. In both figures, we have drawn a finite tree \( T \), a set \( S \) in black, the corresponding set \( \mathcal{B}^+(S) \) in white, and a set \( E \) in red.}%, and a vertex  \(v \in \mathcal{B}^-(S) \) such that \( |\{ v' \in S \colon v \sim v' \}| \geq 2 \) in red.}
\end{figure}
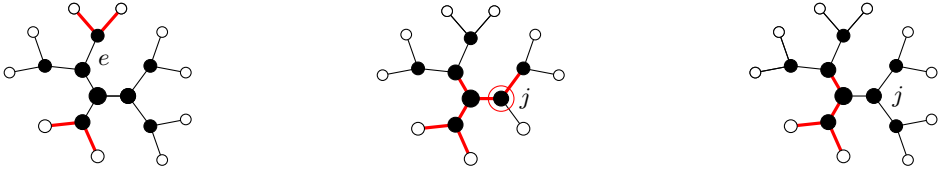

\begin{lemma}\label{lemma: observation 1}
	In the setting of Proposition~\ref{proposition: new negative lemma for p close to one}, let \( E \subseteq E(T)\) be a non-empty multi-set that is not connected (see Figure~\ref{figure: obs 1}). Then, for any \( J \subseteq  \mathcal{B}^-(S), \) we have
	\[
	\Bigl[ \frac{d}{dp_E} \log P\pigl( X \bigl(J \cup (\mathcal{B}^+(S)\smallsetminus \partial_S J) \bigr) \equiv 0 \pigr)  \Bigr]_{(p_e)  \equiv 1} = 0.
	\] 
\end{lemma}

\begin{proof}
	Let \( J \subseteq  \mathcal{B}^-(S), \) and let \( Y \) be random edge set  obtained from the percolation process corresponding to \( X. \)
 	
	Since \( E \) is not connected, there is an edge in \( e \in E(T_S)\smallsetminus E \) that separates \( J \cup (\mathcal{B}^+(S)\smallsetminus \partial_S J)  \) into two non-empty sets \( A \) and \( B.\) On the event \( e \notin Y \)  the events \( X(A) \equiv 0 \) and \( X(B) \equiv 0 \) are independent. The implies that 
	\begin{align*}
		&\Bigl[\frac{d}{dp_E} \log P\pigl( X \bigl(J \cup (\mathcal{B}^+(S)\smallsetminus \partial_S J) \bigr) \equiv 0 \pigr)  \Bigr]_{p_e= 1}
	\\&\qquad=
	\frac{d}{dp_E} \log P\bigl( X (A \sqcup B ) \equiv 0 \mid e \notin Y \bigr) 
	\\&\qquad=
	\frac{d}{dp_E} \log \pigl( P\bigl( X (A  ) \equiv 0  \mid e \notin Y \bigr) P\bigl( X ( B ) \equiv 0  \mid e \notin Y \bigr)   \pigr)
	\\&\qquad=
	\frac{d}{dp_E} \pigl( \log P\bigl( X (A  ) \equiv 0  \mid e \notin Y  \bigr)  + \log P\bigl( X ( B ) \equiv 0  \mid e \notin Y \bigr)   \pigr).
	\end{align*}
	(Here and below,  \( \sqcup\) is used for a disjoint union.)  
	Since \( E \) intersects both \( E(T_A) \) and \( E(T_B)\) and these sets are disjoint, we have
	\begin{align*}
		\frac{d}{dp_{E_A}} P\bigl( X ( B ) \equiv 0  \mid e \notin Y  \bigr) = \frac{d}{dp_{E_B}} P\bigl( X ( A ) \equiv 0  \mid e \notin Y  \bigr) = 0,
	\end{align*}
	and hence
	\[
	\Bigl[\frac{d}{dp_E} \log P\pigl( X \bigl(J \cup (\mathcal{B}^+(S)\smallsetminus \partial_S J) \bigr) \equiv 0 \pigr)  \Bigr]_{(p_e) \equiv 1} = 0. 
	\] 
	This completes the proof.
\end{proof}

The first step of the proof of Proposition~\ref{proposition: new negative lemma for p close to one} is the following lemma, which reduces the calculation of the derivative in~\eqref{eq: nu derivative} to the calculation of the derivative of the logarithm of a probability of a particular event.

\begin{lemma}\label{lemma: first step of proof}
	In the setting of Proposition~\ref{proposition: new negative lemma for p close to one}, for any set \( E \subseteq E(T) , \) we have
	\begin{align*}
		&\Bigl[ \frac{d}{dp_E} \nu(S) \Bigr]_{(p_e) \equiv 1}  =
		\begin{cases}
			0 &\text{if } E(T_S) \not \subseteq E \cr 
			(-1)^{|\mathcal{B}^-(S) |} \Bigl[\frac{d}{dp_E}    
         \log P\pigl( X \bigl(\mathcal{B}^-(S) \bigr)  \equiv 0 \pigr) \Bigr]_{(p_e) \equiv 1} &\text{if } E=E(T_S).
		\end{cases} 
	\end{align*}  
\end{lemma}

\begin{proof}
	To simplify notation, %let \( R \coloneqq \bigl\{ v \in \mathcal{B}^-(S) \colon |\{ v' \in S \colon v \sim v' \}| \geq 2 \} \) and \( R_+ \coloneqq \{ v \in \mathcal{B}^+(S) \colon v \sim R \}, \) and 
	let \( E_S \) be the set of edges between \( S \) and \( S^c. \) 

	Since \( |S|\geq 2 \) and \( S \) is connected, by Lemma~\ref{lemma: trees}, we have 
	\begin{equation} \label{eq: step 1 for sets}
        \nu(S) =  
         \sum_{J \subseteq \mathcal{B}^-(S)} (-1)^{|J|} 
         \log P\pigl( X \bigl(J \cup (\mathcal{B}^+(S)\smallsetminus \partial_S J) \bigr) \equiv 0 \pigr) .
    \end{equation}  
 	 By Lemma~\ref{lemma: observation 1}, for any multiset \( E \subseteq E(T) \) with \( E(T_S) \nsubseteq E, \) we have
 	 \begin{equation*}
 	 	\Bigl[ \frac{d}{dp_{E}} \nu(S)\Bigr]_{(p_e)\equiv 1} = 0.
 	 \end{equation*} 
 	 Next, assume \( E = E(T_S). \)
 	 By~\eqref{eq: step 1 for sets}, we have 
 	 \begin{align*}
		&\Bigl[ \frac{d}{dp_E} \nu(S) \Bigr]_{(p_e) \equiv 1} 
		=
		\Bigl[ \frac{d}{dp_E} \sum_{J \subseteq \mathcal{B}^-(S)} (-1)^{|J|} 
         \log P\pigl( X \bigl(J \cup (\mathcal{B}^+(S)\smallsetminus \partial_S J) \bigr) \equiv 0 \pigr) \Bigr]_{(p_e) \equiv 1}.
		%\\&\qquad=
		%\Bigl[ \frac{d}{dp_E} \sum_{J \subseteq \mathcal{B}^-(S)\smallsetminus R} (-1)^{|J|} 
         %\log P\pigl( X \bigl(J \cup (\mathcal{B}^+(S)\smallsetminus \partial_S J) \bigr) \equiv 0 \pigr) \Bigr]_{(p_e) \equiv 1} .
	\end{align*} 
	
	Let \( Y \) be the random edge set corresponding to the percolation process corresponding to \( X \) (see (B') of Definition~\ref{def: TIMC general}).
	Then, since \( E = E(T_S), \) we have
	\begin{align*}
		&\Bigl[ \frac{d}{dp_E} \sum_{J \subseteq \mathcal{B}^-(S) } (-1)^{|J|} 
         \log P\pigl( X \bigl(J \cup (\mathcal{B}^+(S)\smallsetminus \partial_S J) \bigr) \equiv 0 \pigr) \Bigr]_{(p_e) \equiv 1} 
         \\&\qquad=
         \Bigl[ \frac{d}{dp_E} \sum_{J \subseteq \mathcal{B}^-(S) } (-1)^{|J|} 
         \log P\pigl( X \bigl(J \cup (\mathcal{B}^+(S)\smallsetminus \partial_S J) \bigr) \equiv 0 \mid Y \cap E_S = \emptyset\pigr) \Bigr]_{(p_e) \equiv 1}  
         \\&\qquad=
         \Bigl[ \frac{d}{dp_E} \sum_{J \subseteq \mathcal{B}^-(S) } (-1)^{|J|} 
         \log \pigl[ P\bigl( X (J  ) \equiv 0    \bigr)  \prod_{j \in \mathcal{B}^-(S) \smallsetminus J } r^{|\partial_S \{ j \}|} \pigr] \Bigr]_{(p_e) \equiv 1},
	\end{align*}
	where the last equality follows since on the event \( Y \cap E_S = \emptyset, \) the random vectors \( X(S) \) and \( X(\mathcal{B}^+(S)) \) are independent. 
	Combining the previous equations, we obtain
	\begin{align*}
		&\Bigl[ \frac{d}{dp_E} \nu(S) \Bigr]_{(p_e) \equiv 1}  
		= 
		\Bigl[ \frac{d}{dp_E} \sum_{J \subseteq \mathcal{B}^-(S) } (-1)^{|J|} 
         \log \pigl[ P\bigl( X (J  ) \equiv 0   \bigr)\prod_{j \in \mathcal{B}^-(S) \smallsetminus J } r^{|\partial_S \{ j \}|}  \pigr] \Bigr]_{(p_e) \equiv 1} 
         \\&\qquad=
		\Bigl[ \frac{d}{dp_E} \sum_{J \subseteq \mathcal{B}^-(S) } (-1)^{|J|} 
         \log P\bigl( X (J  ) \equiv 0   \bigr) \Bigr]_{(p_e) \equiv 1} 
         \\&\qquad\qquad+
		\Bigl[ \frac{d}{dp_E} \sum_{J \subseteq \mathcal{B}^-(S) } (-1)^{|J|} 
         \log  \prod_{j \in \mathcal{B}^-(S) \smallsetminus J }   r^{|\partial_S \{ j \}|}\Bigr]_{(p_e) \equiv 1} 
         \\&\qquad=
		\Bigl[ \frac{d}{dp_E} \sum_{J \subseteq \mathcal{B}^-(S) } (-1)^{|J|} 
         \log P\bigl( X (J  ) \equiv 0   \bigr) \Bigr]_{(p_e) \equiv 1} .
	\end{align*}  	
	If \( J \neq \mathcal{B}^-(S) , \) then 
	\begin{equation*}
	 	\Bigl[ \frac{d}{dp_E}  
         \log P\bigl( X (J ) \equiv 0 \bigr) \Bigr]_{(p_e) \equiv 1} =0 
	\end{equation*}
	since the event \( X(J )\equiv 0 \) does not depend on all edges in \( E = E(T_S). \) From this the desired conclusion immediately follows.
\end{proof}

The following lemma, Lemma~\ref{lemma: observation 3} below, shows that in the setting of Proposition~\ref{proposition: new negative lemma for p close to one}, if the set \( E \) is not adjacent to some vertex \( v \in \mathcal{B}^-(S), \) then the derivative \( \frac{d}{dp_E} \) of the terms appearing on the righ-hand side of~\eqref{eq: nu eq} of Lemma~\ref{lemma: trees} is zero.

\begin{lemma}\label{lemma: observation 3}
	In the setting of Proposition~\ref{proposition: new negative lemma for p close to one}, let \( E \subseteq E(T) \) be a non-empty multi-set. Assume that   \( j \in \mathcal{B}^-(S) \) is such that \( E_j \cap E = \emptyset\) (see Figure~\ref{figure: obs 3}). In other words, assume there is no edge in \( E \) with one endpoint in \( j. \) Then 
	\[
		\Bigl[\frac{d}{dp_E}\nu(S) \Bigr]_{(p_e)\equiv 1}
		= 0.
	\]
\end{lemma}

\begin{proof}
	For \( I \subseteq  \mathcal{B}^-(S), \) consider the event	
	\begin{equation*}
		\mathcal{E}_I \coloneqq \bigl\{ X \bigl(I \cup (\mathcal{B}^+(S)\smallsetminus \partial_S I \bigr)\equiv 0 \bigr\}.
	\end{equation*}
	Then, by Lemma~\ref{lemma: trees}, we have
	
	\begin{equation*}
		\nu(S) = \sum_{I \subseteq \mathcal{B}^-(S)} (-1)^{|J|} 
         \log P( \mathcal{E}_I ) .
	\end{equation*}

	 Let \( J \subseteq  \mathcal{B}^-(S)\smallsetminus \{ j \}.\) Then, since \( E_j \cap E = \emptyset, \)  
	\begin{align*}
	&
         \Bigl[ \frac{d}{dp_E} 
         \log P( \mathcal{E}_{J})  \Bigr]_{(p_e)\equiv 1}
	=
	\Bigl[\frac{d}{dp_E} \log \bigl(
         P(\mathcal{E}_{J\sqcup\{ j\}} ) r^{|\partial_S \{ j \}|-1} \bigr)\Bigr]_{(p_e)\equiv 1} 
	\\&\qquad=
	\Bigl[\frac{d}{dp_E} \log  
         P(\mathcal{E}_{J\sqcup\{ j\}} ) \Bigr]_{(p_e)\equiv 1} 
	+
	\frac{d}{dp_E} \log \bigl(  r^{|\partial_S \{ j \}|-1} \bigr) 
	=
	\Bigl[\frac{d}{dp_E} \log  
         P(\mathcal{E}_{J\sqcup\{ j\}} ) \Bigr]_{(p_e)\equiv 1} .
	\end{align*}
	Since \( (-1)^{|J\sqcup \{ j\}|} = -(-1)^{|J|}, \) it follows that  
	\begin{equation*}
		\Bigl[ \frac{d}{dp_E} \sum_{J \subseteq \mathcal{B}^-(S)} (-1)^{|J|} 
         \log P\pigl( X \bigl(J \cup (\mathcal{B}^+(S)\smallsetminus \partial_S J) \bigr) \equiv 0 \pigr) \Bigr]_{(p_e)\equiv 1} = 0.
	\end{equation*}  
	From this, the desired conclusion immediately follows.
\end{proof}

The following lemma collects a useful observation obtained by combining the previous result, Lemma~\ref{lemma: first step of proof}, with~\cite[Theorem 5.6]{fgs}.

\begin{lemma}\label{lemma: polylog key}
	Let \( n \geq 2, \) and let \( T \) be a tree be a \( n \)-star graph, consisting of one vertex, called \( o \), of degree \( n \), and \( n \) leaves (see Figure~\ref{figure: key examples}), and let \( S \subseteq \mathcal{B}^-(S). \) Further, let \( X \) be a tree indexed Markov chain on \( T \) with parameters \( (r,(p_e)), \) and let \( \nu \) be the corresponding signed measure. Then 
	\begin{align*}
		&%\Bigl[ \frac{d}{dp_{E(T_S)}} \nu(S) \Bigr]_{(p_e) \equiv 1}  = 
          %(-1)^{|\mathcal{B}^-(S) |} 
          \Bigl[ \frac{d}{dp_{E(T_S)}}   
         \log P\pigl( X \bigl(\mathcal{B}^-(S) \bigr)  \equiv 0 \pigr) \Bigr]_{(p_e) \equiv 1}  
         =
         %(-1)^{|\mathcal{B}^-(S) |} 
         -\Li_{1-|\mathcal{B}^-(S) |}(-r^{-1}(1-r)) .
	\end{align*} 
\end{lemma}

\begin{figure}[h]
	\begin{subfigure}[t]{.3\textwidth}\centering
			\begin{tikzpicture}[scale=.4]

    	\foreach \i in {0,1,2} 
    			{ 
    				\draw[] (0,0)  -- ({2*cos(120*\i)},{2*sin(120*\i)})  ;
    				\filldraw[fill=white] ({2*cos(120*\i)},{2*sin(120*\i)}) circle (5pt);   
    			} 
    	\filldraw[fill=white] (0,0) circle (5pt);  
    	
    	\end{tikzpicture}	
    	\caption{\(n=3\)}
    \end{subfigure}	
    \hfil
	\begin{subfigure}[t]{.3\textwidth}\centering
			\begin{tikzpicture}[scale=.4]

    	\foreach \i in {0,1,2,3,4} 
    			{ 
    				\draw[] (0,0)  -- ({2*cos(72*\i)},{2*sin(72*\i)})  ;
    				\filldraw[fill=white] ({2*cos(72*\i)},{2*sin(72*\i)}) circle (5pt);   
    			} 
    	\filldraw[fill=white] (0,0) circle (5pt);  
    	
    	\end{tikzpicture}	
    	\caption{\(n=5\)}
    \end{subfigure}		
    \hfil
	\begin{subfigure}[t]{.3\textwidth}\centering
			\begin{tikzpicture}[scale=.4]

    	\foreach \i in {0,1,2,3,4,5,6,7} 
    			{ 
    				\draw[] (0,0)  -- ({2*cos(45*\i)},{2*sin(45*\i)})  ;
    				\filldraw[fill=white] ({2*cos(45*\i)},{2*sin(45*\i)}) circle (5pt);   
    			} 
    	\filldraw[fill=white] (0,0) circle (5pt);  
    	
    	\end{tikzpicture}	
    	\caption{\(n=8\)}
    \end{subfigure}		
    	\caption{The three figures above show examples of the graphs considered in Lemma~\ref{lemma: polylog key}.}	\label{figure: key examples}
\end{figure}

\begin{proof}
	Let \( X' \) be the restriction of \( X \) to \( \mathcal{B}^-(S). \) Let \( \nu' \) be the signed measure corresponding to \( X'. \) Then, by Proposition~\ref{prop: connected}, we have \(  \nu(S) = \nu'(\mathcal{B}^-(S)) . \) Moreover, we note 
	\begin{align*}
		X' = r \Pi_{p(1-r)} + (1-r) \Pi_{1-p+p(1-r)},
	\end{align*}
	where for \( a \in [0,1], \) \( \Pi_a \) denotes a product distribution on \( \{ 0,1 \}^n \) with \( \Pi_a(X(1)=1) = a. \)
	
	Now let \( \alpha_1 \coloneqq r, \) \( \alpha_2 \coloneqq 1-r, \) \( x_1 \coloneqq 1, \) \( x_2 \coloneqq pr/(1-p(1-r)), \) and \( q \coloneqq 1-p(1-r).\) With this notation, we have
	\begin{equation*}
		X' \sim \alpha_1 \Pi_{1-qx_1} + \alpha_2 \Pi_{1-qx_2}.
	\end{equation*}
	Then, by the proof of~\cite[Theorem 5.6]{fgs} we have
	\begin{align*}
		\Bigl[ \frac{d^\ell}{dx_2^\ell} \nu'(\mathcal{B}^-(S))\Bigr]_{x_2=1} = \begin{cases}
			0 & \text{if } \ell < |\mathcal{B}^-(S)|\cr 
			\ell! (-1)^{\ell+1} \Li_{1-\ell}(-\alpha_1^{-1}(1-\alpha_1)) &\text{if } \ell=|\mathcal{B}^-(S)|.
		\end{cases}
	\end{align*}
	Since \( x_2 = 1 \Leftrightarrow p=1\), it follows that 
	\begin{align*} 
		\Bigl[\frac{d^\ell}{dp^\ell}\nu(S) \Bigr]_{p=1} = \begin{cases}
			0 & \text{if } \ell < |\mathcal{B}^-(S)|\cr 
			\ell! (-1)^{\ell+1} \Li_{1-\ell}(-r^{-1}(1-r)) &\text{if } \ell=|\mathcal{B}^-(S)|.
		\end{cases}
	\end{align*} 
	At the same time, by Lemma~\ref{lemma: first step of proof}, for any \( E \subseteq E(T), \) we have
	\begin{align*}
		&\Bigl[ \frac{d}{dp_E} \nu(S) \Bigr]_{(p_e) \equiv 1}  = \begin{cases} 0 &\text{if } E(T_S) \nsubseteq E \cr 
          (-1)^{|\mathcal{B}^-(S) |} \Bigl[\frac{d}{dp_E}    
         \log P\pigl( X \bigl(\mathcal{B}^-(S) \bigr)  \equiv 0 \pigr) \Bigr]_{(p_e) \equiv 1} &\text{if } E=E(T_S) .
         \end{cases}
	\end{align*} 
	Combining these two equations, it follows that if \( E = E(T_S), \) then 
	\begin{align*}
		&\Bigl[ \frac{d}{dp_E} \nu(S) \Bigr]_{(p_e) \equiv 1}  = 
          (-1)^{|\mathcal{B}^-(S) |} \Bigl[\frac{d}{dp_E}    
         \log P\pigl( X \bigl(\mathcal{B}^-(S) \bigr)  \equiv 0 \pigr) \Bigr]_{(p_e) \equiv 1}  
         \\&\qquad =
         (-1)^{|E(T_S)|+1} \Li_{1-|E(T_S)|}(-r^{-1}(1-r)) 
	\end{align*} 
	Noting that for the specific tree considered here, we have \( |E(T_S)| = |\mathcal{B}^-(S)|, \) this concludes the proof.
\end{proof}

The next main result of this section, which will be used in the proof of Proposition~\ref{proposition: new negative lemma for p close to one} is the following proposition, which gives a concrete expression for derivatives of the terms appearing on the right-hand side of~\eqref{eq: nu eq} in Lemma~\ref{lemma: trees} in terms of polylogarithm functions.

\begin{proposition}\label{proposition: 2/3 trees}
	Let \( T \) be a tree, let \( S \subseteq V(T) \) be connected and have cardinality at least two.  
	%Assume further that for all \(  v \in \mathcal{B}^-(S) \)  we have \(  |N(v) \cap S|=1 .\) 
	For \( v \in V(T) \) and \( e \in E(T), \) let \( r_v \in (0,1] \) and \( p_e \in [0,1],\) and let \( X \) be the tree-indexed Markov chain on \( T \) with parameters \( ((r_v),(p_e)). \) Then 
	\begin{align*}
		&\Bigl[ \frac{d}{dp_{E(T_S)}}  
         \log P\pigl( X \bigl(\mathcal{B}^-(S) \bigr) \equiv 0 \pigr) \Bigr]_{(p_e) \equiv  1}
         =
         -\frac{1-r}{r} \prod_{v \in S \smallsetminus L} f_{ \deg v}(r) 
	\end{align*} 
	where for \( j \geq 2, \)  \( k_j \) is the number of vertices of degree \( j\) in \( T_S,\) \( [j] = \{ 1,2,\dots, j\}, \) and
	\begin{align}\label{eq: fdef}
		f_j(r) \coloneqq \frac{r}{1-r} \sum_{\substack{(F_1,\dots, F_m) \in \mathcal{P}([j]) \colon\\ |F_1|,\dots, |F_m|\geq 2}}   (r-1)^{m} 
			\prod_{i =1}^m  
			\Bigl( 1   -   ( (r-1)/r)^{|F_i|-1}  \Bigr).
	\end{align}   
\end{proposition}

Before we prove Proposition~\ref{proposition: 2/3 trees}, we will state and prove a number of shorter lemmas that will be useful in its proof.

\begin{lemma}\label{lemma: new lemma 1}
	In the setting of Proposition~\ref{proposition: 2/3 trees}, let \( E \subseteq E(T_S) \) and assume that there is \( v \in V(T_S) \smallsetminus \mathcal{B}^-(S) \) that is the end-point of exactly one edge in \( e \in E.\)  %\label{item: new lemma item 1}
		%\item there is some edge \( e \in E \) that is not connected to any vertex in \( \mathcal{B}^-(S). \) \label{item: new lemma item 2}
	%\end{enumerate}
	Then
	\begin{align*} 
		& \Bigl[ \frac{d}{dp_{E}}P\bigl( X (\mathcal{B}^-(S) ) \equiv 0 \bigr)\Bigr]_{(p_e)\equiv  1} =  0. 
	\end{align*}
\end{lemma}

\begin{proof}
	Let \( Y \) be the random edge set corresponding to the percolation process corresponding to \( X \) (see (B') of Definition~\ref{def: TIMC general}).
	Note that
	\begin{align*}
		&\Bigl[ \frac{d}{dp_{E}}P\bigl( X (\mathcal{B}^-(S) ) \equiv 0 \bigr)\Bigr]_{(p_e)\equiv  1}
		=
		\Bigl[\frac{d}{dp_{E}}P\pigl( X \bigl(\mathcal{B}^-(S) \bigr) \equiv 0 \mid (E(T_S) \smallsetminus E) \cap Y = \emptyset\pigr)\Bigr]_{(p_e)\equiv  1}.
	\end{align*}
	Let \( v \in V(T_S) \smallsetminus \mathcal{B}^-(S) \) be the end-point of exactly one edge  \( e \in E.\) Then, since \( (r_v) \equiv r \) and we obtain \( X \) using \( (B')\) of Definition~\ref{def: TIMC general}, conditioned on the event \( (E(T_S) \smallsetminus E) \cap Y = \emptyset \) the two events \( X (\mathcal{B}^-(S) ) \equiv 0  \) and \( e \in Y \) are independent, and hence
	\begin{align*}
		&\Bigl[ \frac{d}{dp_{E}}P\pigl( X \bigl(\mathcal{B}^-(S) \bigr) \equiv 0 \big| (E(T_S) \smallsetminus E) \cap Y = \emptyset\bigr) \pigr)\Bigr]_{(p_e) \equiv  1} = 0.
	\end{align*}
	Combining the previous equations, we obtain the desired conclusion.
\end{proof}

\begin{lemma}\label{lemma: proof of case 1}
		In the setting of Proposition~\ref{proposition: 2/3 trees}, let \( v \in S\smallsetminus  \mathcal{B}^-(S) \) be such that \( L \coloneqq  N(v) \cap  \mathcal{B}^-(S) \) is nonempty. Let \( E \subseteq E_L \) satisfy \( |E| \geq 2, \) and let \( L_E \) be the set of all vertices in \( L \) that are the endpoint of some edge in \( E.\) Then
		\begin{align*}
			&\Bigl[ \frac{d}{dp_{E}} P\bigl( X (\mathcal{B}^-(S) ) \equiv 0 \bigr) \Bigr]_{(p_e) \equiv  1}
			=
			r^{|\mathcal{B}^-(S)\smallsetminus L_E|} \bigl( r  (r-1)^{|E|} + (1-r)  r^{|E|} \bigr).
			%r^{|\mathcal{B}^-(S)\smallsetminus L_E|}\Bigl[ 
			%\frac{d}{dp_{E}} P\bigl( X (L_E)\equiv 0 \bigr) \Bigr]_{(p_e) \equiv  1}.
		\end{align*}
		%Moreover,
		%\begin{align*}
		%	&\Bigl[
		%	\frac{d}{dp_{E}} P\bigl( X (L_E)\equiv 0 \bigr) \Bigr]_{(p_e) \equiv  1}
		%	=
		%	 r  (r-1)^{|E|} + (1-r)  r^{|E|} .
		%\end{align*}
	\end{lemma}
	
	\begin{proof}
		%The first statement of the lemma follows immediately from the definition of \( X. \) 
		We have
		\begin{align*}
			&\Bigl[ 
			\frac{d}{dp_{E}} P\bigl( X (L_E)\equiv 0 \bigr) \Bigr]_{(p_e) \equiv  1}
			=
			\Bigl[ 
			\frac{d}{dp_{E}} \Bigl( r \prod_{e \in E} (1-p_e+p_e r) + (1-r) \prod_{e \in E } (p_e r) \Bigr) \Bigr]_{(p_e) \equiv  1}
			\\&\qquad=
			r  \prod_{e \in E} (-1+ r) + (1-r) \prod_{e \in E} r 
			=
			r  (-1+ r)^{|E|} + (1-r)  r^{|E|}  .
		\end{align*}
		This concludes the proof.
	\end{proof}

	\begin{lemma}\label{lemma: proof of case 2}
		In the setting of Proposition~\ref{proposition: 2/3 trees}, let \( v \in S\smallsetminus  \mathcal{B}^-(S) \) be such that \( L \coloneqq  N(v) \cap  \mathcal{B}^-(S) \) is nonempty, and let \( E \subseteq E(T_S)\smallsetminus E_L .\) Then
		\begin{align*}
			&\Bigl[ \frac{d}{dp_{E}} P\bigl( X (\mathcal{B}^-(S) ) \equiv 0 \bigr) \Bigr]_{(p_e)\equiv  1}
			=
			r^{|L|-1}
			\Bigl[  \frac{d}{dp_{E}} P\pigl( X \bigl( \mathcal{B}^-(S\smallsetminus L) \bigr)   \equiv 0 \pigr) \Bigr]_{(p_e)\equiv  1} .
		\end{align*}
	\end{lemma}
	
	\begin{proof}
	Recall that \( E_v \subseteq E(T_S) \) is the set of all edges adjacent to~\( v\), and note that on the event \( E_v \cap Y = \emptyset, \) the random variables \( X(\mathcal{B}^-(S)\smallsetminus L), \) \( X(v), \) and \( X(v') \) for \( v' \in L \) are independent. Using this observation, and noting that \( v \in \mathcal{B}^-(S\smallsetminus L), \) it follows that
		\begin{align*}
			&\Bigl[ \frac{d}{dp_{E}} P\bigl( X (\mathcal{B}^-(S) ) \equiv 0 \bigr) \Bigr]_{(p_e)\equiv  1}
			=
			\Bigl[ \frac{d}{dp_{E}} P\bigl( X (\mathcal{B}^-(S) ) \equiv 0 \mid E_v \cap Y = \emptyset  \bigr) \Bigr]_{(p_e)\equiv  1}
			\\&\qquad=
			\Bigl[ r^{|L|} \frac{d}{dp_{E}} P\bigl( X (\mathcal{B}^-(S)\smallsetminus L ) \equiv 0 \bigr) \Bigr]_{(p_e)\equiv  1} 
			=
			\Bigl[ r^{|L|-1} \frac{d}{dp_{E}} P\pigl( X \bigl(\mathcal{B}^-(S\smallsetminus L)\bigr)   \equiv 0 \pigr) \Bigr]_{(p_e)\equiv  1} .
		\end{align*}
		This concludes the proof.
	\end{proof}

	\begin{lemma}\label{lemma: proof of case 3}
		In the setting of Proposition~\ref{proposition: 2/3 trees}, let \( v \in S\smallsetminus  \mathcal{B}^-(S) \) be such that \( L \coloneqq  N(v) \cap  \mathcal{B}^-(S) \) is nonempty, and let \( e_v \) be the unique edge connecting \( v \) to \( E(T_S) \smallsetminus E(T_{L \cup \{ v \}}) \). Let \( E \subseteq E(T_S)\)  and assume that \( e_v \in E\) and that \( E \cap E_L \neq \emptyset .\) Let \( L_E \) be the set of all vertices in \( L \) that are the endpoint of some edge in \( E.\)  
		Then
		\begin{align*}
			& \Bigl[ \frac{d}{dp_{E}} P\bigl( X (\mathcal{B}^-(S) ) \equiv 0 \bigr)\Bigr]_{(p_e)\equiv  1}
			\\&\qquad = 
			r^{|L\smallsetminus L_E|}\bigl(  ( r-1)^{|E \cap E_L|} -   r^{|E \cap E_L|} 
			\bigr)  	\Bigl[ \frac{d}{dp_{E \smallsetminus E_L}} P\pigl( X \bigl(\mathcal{B}^-(S)\smallsetminus L  \bigr)  \equiv 0 \pigr) \Bigr)  \Bigr]_{(p_e) \equiv 1}.
		\end{align*} 
	\end{lemma}
	
	\begin{proof}
		Recall that \( Y \) is the set of retained edges in the edge percolation process used to construct X. 
		If we sum over all possiblities for \( Y \cap (E \cap E_L ),\) we get  
		\begin{align*}
			& P\bigl( X (\mathcal{B}^-(S) ) \equiv 0 \bigr)
			= 
			 \sum_{ E_0 \subseteq E \cap E_L} P\bigl( X (\mathcal{B}^-(S) ) \equiv 0 \text{ and } Y\cap (E \cap E_L)  = E_0 \bigr).
		\end{align*}
		If \( E_0 \subseteq E \cap E_L \) is empty, then the corresponding probability in the above sum does not depend on \( p_{e_v}, \) and hence
		\begin{align*}
			\frac{d}{dp_{E}}  P\pigl( X \bigl(\mathcal{B}^-(S) \bigr) \equiv 0 \text{ and } Y\cap (E \cap E_L)  = E_0 \pigr)= 0.
		\end{align*} 
		Now, instead assume that \( { E_0 \subseteq E \cap E_L} \) is non-empty. In this case, noting that \( v \in \mathcal{B}^-(S)\smallsetminus L, \) we have
		\begin{align*}
			&
			 P\bigl( X (\mathcal{B}^-(S) ) \equiv 0 \text{ and } Y \cap (E \cap E_L) = E_0 \bigr) 
			\\&\qquad=  
			r^{-1} P\bigl (X(L) \equiv 1 \text{ and } Y \cap (E \cap E_L)  = E_0  \bigr)			P\pigl( X \bigl(\mathcal{B}^-(S)\smallsetminus L\bigr)  \equiv 0 \pigr).
		\end{align*}
		Combining the two cases, we get
		\begin{align*}
			&
			 \frac{d}{dp_E} \sum_{E_0 \subseteq E \cap E_L} P\bigl( X (\mathcal{B}^-(S) ) \equiv 0 \text{ and } Y \cap (E \cap E_L) = E_0 \bigr)  
			\\&\qquad=  
			\frac{d}{dp_E} \Bigl( r^{-1}   P\bigl(X(L) \equiv 0    \text{ and } Y \cap (E \cap E_L)  \neq  \emptyset  \bigr)		P\pigl( X \bigl(\mathcal{B}^-(S)\smallsetminus L  \bigr)  \equiv 0 \pigr) \Bigr)
			\\&\qquad=  
			 r^{-1}  \frac{d}{dp_{E \cap E_L}} P\bigl(X(L) \equiv 0    \text{ and } Y \cap (E \cap E_L)  \neq  \emptyset  \bigr)		\frac{d}{dp_{E \smallsetminus E_L}} P\pigl( X \bigl(\mathcal{B}^-(S)\smallsetminus L \bigr)  \equiv 0 \pigr) \Bigr).
		\end{align*}
		Now note that for the first of these terms, we have
		\begin{align*}
			&\Bigl[ \frac{d}{dp_{E \cap E_L}}  P\bigl(X(L) \equiv 0    \text{ and } Y \cap (E \cap E_L)  \neq  \emptyset  \bigr)	\Bigr]_{{p_e}\equiv 1}
			\\&\qquad
			=
			r^{|L\smallsetminus L_E|}\Bigl[ \frac{d}{dp_{E \cap E_L}}  P\bigl(X(L_E) \equiv 0    \text{ and } Y \cap (E \cap E_L)  \neq  \emptyset  \bigr)	\Bigr]_{{p_e}\equiv 1}
			\\&\qquad
			=
			r^{|L\smallsetminus L_E|}\Bigl[ \frac{d}{dp_{E \cap E_L}}  P\bigl(X(L_E) \equiv 0   \bigr) -  P\bigl(X(L_E) \equiv 0    \text{ and } Y \cap (E \cap E_L)  =  \emptyset  \bigr)	
			\Bigr]_{{p_e}\equiv 1}.
		\end{align*}
		Since
		\begin{align*}
			&\frac{d}{dp_{E \cap E_L}} P\bigl(X(L_E) \equiv 0    \text{ and } Y \cap (E \cap E_L)  =  \emptyset  \bigr)	
			=
			\frac{d}{dp_{E \cap E_L}} \prod_{e \in E \cap E_L} rp_e 
			= r^{|E \cap E_L|}
		\end{align*}
		and 
		\begin{align*}
			&\frac{d}{dp_{E \cap E_L}} P\bigl(X(L_E) \equiv 0   \bigr)	
			=
			\frac{d}{dp_{E \cap E_L}} \Bigl( r\prod_{e \in E \cap E_L} \bigl( (1-p_e) + rp_e \bigr) + (1-r) \prod_{e \in E \cap E_L} \bigl( rp_e \bigr) \Bigr) 
			\\&\qquad=
			r (  -1+r)^{|E \cap E_L|} + (1-r) r^{|E \cap E_L|} ,
		\end{align*}
		the desired conclusion follows.
	\end{proof}

We now proceed to the proof of Proposition~\ref{proposition: 2/3 trees}
\begin{proof}[Proof of Proposition~\ref{proposition: 2/3 trees}]
To simplify the notation in the proof, for a connected set \( S_0 \subseteq S, \) we let 
	\begin{align*}
		&A(S_0) \coloneqq \Bigl[ \frac{d}{dp_{E(T_{S_0})}}  
         \log P\pigl( X \bigl(\mathcal{B}^-(S_0) \bigr) \equiv 0 \pigr) \Bigr]_{(p_e) \equiv  1}
	\end{align*} 
	and note that with this notation, we need to prove that
		\begin{align*}
		&A(S) =
         -\frac{1-r}{r} \prod_{j=2}^\infty  \biggl( \frac{\mathrm{Li}_{1-j}(-(1-r)/r)}{r^{j-1}(1-r)} \biggr)^{k_j}.
	\end{align*}

	Fix any \( o \in V(T_S) \smallsetminus \mathcal{B}^-(S). \) Let \( Y \) be the random set corresponding to the edge percolation process corresponding to \( X. \)

	Let \( \mathcal{P}(E(T_S))  \) denote the set of all partitions of \( E(T_S) \) into non-empty sets. Then 
	\begin{equation}
		\begin{split}
		&A(S) = \Bigl[ \frac{d}{dp_{E(T_S)}}  
         \log P\bigl( X (\mathcal{B}^-(S) ) \equiv 0 \bigr) \Bigr]_{(p_e)\equiv  1} 
         \\&\qquad=
          \sum_{(E_1, \dots , E_k) \in \mathcal{P}(E(T_S))}
         \frac{\bigl[(-1)^{k-1}\prod_{j=1}^k \frac{d}{dp_{E_j}}P\bigl( X (\mathcal{B}_T^-(S) ) \equiv 0 \bigr)\bigr]_{(p_e) \equiv  1}}{\bigl[ P\bigl( X (\mathcal{B}^-(S) ) \equiv 0 \bigr)^k \bigr]_{(p_e) \equiv  1}},
         \\&\qquad=
         \sum_{(E_1, \dots , E_k) \in \mathcal{P}(E(T_S))}
         (-1)^{k-1} r^{-k|\mathcal{B}^-(S)|} \prod_{j=1}^k \Bigl[\frac{d}{dp_{E_j}}P\bigl( X (\mathcal{B}^-(S) ) \equiv 0 \bigr)\Bigr]_{(p_e) \equiv  1} ,\label{eq: second line with partition}
		\end{split}
	\end{equation} 
	If  \( (E_1, \dots , E_k) \in \mathcal{P}(E(T_S)) \) ontains either a partition element that is not connected, or a partition element \( E_j \) such that a vertex \( v \in V(T_S)\smallsetminus \mathcal{B}^-(S) \) is the endpoint of exactly one edge in \( E_j \), then, by Lemma~\ref{lemma: observation 1} applied with \( J = \mathcal{B}^-(S) \), and Lemma~\ref{lemma: new lemma 1}, we have 
	\begin{equation*}
		\prod_{j=1}^k \Bigl[  \frac{d}{dp_{E_j}}P\bigl( X (\mathcal{B}^-(S) ) \equiv 0 \bigr)\Bigr]_{(p_e) \equiv  1} = 0.
	\end{equation*}
	 Hence, in~\eqref{eq: second line with partition} we only need to sum over partitions \( (E_1,\dots, E_k) \) of \( E(T_S) \) such that for each \( j \in [k],\) the following holds.
	 \begin{enumerate}[label=(\roman*)]
		\item \( E_j \) is connected, and
		\item no vertex in \( V(T_S)\smallsetminus \mathcal{B}^-(S)\) is the endpoint of exactly one edge in \( E_j. \) \label{item: hat condition ii}
	\end{enumerate} 
	Let \( \hat {\mathcal{P}}(E(T_S)) \) be the set of all such partitions.  
	%We note that if each vertex of \( T_S \) has degree at most three, then \( \hat {\mathcal{P}}(E(T_S)) = \{ (E(T_S))\}.\) Moreover, w
	We note that if \( ((E_1,\dots, E_k) \in \hat {\mathcal{P}}(E(T_S)), \) then for all \( j \in [k],\) 
	\begin{enumerate}[label=(\roman*), resume]
		\item  \( E_j \) is connected to some vertex in \( \mathcal{B}^-(S), \) and \label{item: hat condition i} 
	\end{enumerate} 

	Now let \( v \in S\smallsetminus  \mathcal{B}^-(S) \) be such that \( L \coloneqq  N(v) \cap  \mathcal{B}^-(S) \) is nonempty. 	
	Let \( e_v \) be the unique edge connecting \( v \) to \( E(T_S) \smallsetminus E(T_{L \cup \{ v \}}) \).
	Assume that \( E \subseteq E(T_S) \) satisfies~\ref{item: hat condition i} and~\ref{item: hat condition ii},  let \( L_E \) be the set of all vertices in \( L \) that are the endpoint of some edge in \( E,\) and recall that  \( E_L \) is the set of all edges in \( E(T_S) \) which have one endpoint in \( L. \) Note that since we can assume \( E \) to be connected, if it doesn't contain the edge \( e_v \), it either contains only edges from \(E_L\) or no edge from \(E_L\).
We thus have three cases;~\ref{item: case 1} \( E \subseteq E_L \), \ref{item: case 2}~\( E \cap E_L = \emptyset \), and~\ref{item: case 3}~\( e_v \in E\)  (see Figure~\ref{figure: cases}).
	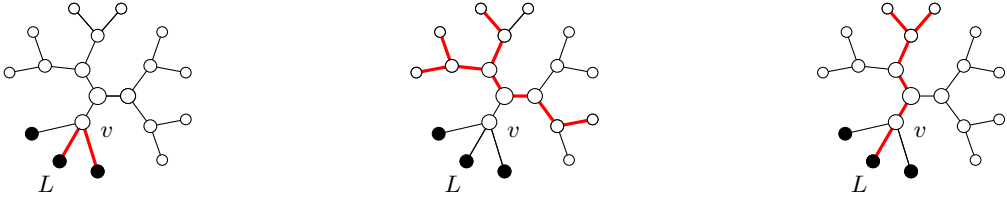
\begin{figure}[]
		\begin{subfigure}[t]{.3\textwidth}\centering
			\begin{tikzpicture}[scale=.4]

    	\foreach \i in {0,1}
    		\foreach \j in {-1,1}
    			\foreach \k in {-1,1}
    			{ 
    				\draw[] ({2*cos(120*\i+30*\j)},{2*sin(120*\i+30*\j)}) -- ({3*cos(120*\i+30*\j+15*\k)},{3*sin(120*\i+30*\j+15*\k )});
    				\filldraw[fill=white] ({3*cos(120*\i+30*\j+15*\k)},{3*sin(120*\i+30*\j+15*\k)}) circle (5pt);   
    			}
    	  
    	  \foreach \i in {1}
    		\foreach \j in {-1}
    			\foreach \k in {-1,1}
    			{ 
    				\draw[] ({2*cos(120*\i+30*\j)},{2*sin(120*\i+30*\j)}) -- ({3*cos(120*\i+30*\j+15*\k)},{3*sin(120*\i+30*\j+15*\k )});
    				\filldraw[fill=white] ({3*cos(120*\i+30*\j+15*\k)},{3*sin(120*\i+30*\j+15*\k)}) circle (5pt);   
    			}

    	\foreach \i in {0,...,1} 
    		\foreach \j in {-1,1}
    		{
    			\draw ({cos(120*\i)},{sin(120*\i)}) -- ({2*cos(120*\i+30*\j)},{2*sin(120*\i+30*\j)});
    			\filldraw[fill=white] ({2*cos(120*\i+30*\j)},{2*sin(120*\i+30*\j)}) circle (6pt);   
    		}

    	\foreach \i in {2} 
    		\foreach \j in {-1,0,1}
    		{
    			\draw[] ({1.4*cos(120*\i)},{1.4*sin(120*\i)}) -- ({2.9*cos(120*\i+30*\j)},{2.9*sin(120*\i+30*\j)});
    			
    			\filldraw[fill=black] ({2.9*cos(120*\i+30*\j)},{2.9*sin(120*\i+30*\j)}) circle (6pt);   
    			%\draw[red] ({2*cos(120*\i+30*\j)},{2*sin(120*\i+30*\j)}) circle (11pt);   
    		}
    	
    	\foreach \i in {2} 
    		\foreach \j in {0,1}
    		{
    			\draw[red, very thick] ({1.4*cos(120*\i)},{1.4*sin(120*\i)}) -- ({2.9*cos(120*\i+30*\j)},{2.9*sin(120*\i+30*\j)});
    			
    			\filldraw[fill=black] ({2.9*cos(120*\i+30*\j)},{2.9*sin(120*\i+30*\j)}) circle (6pt);   
    			%\draw[red] ({2*cos(120*\i+30*\j)},{2*sin(120*\i+30*\j)}) circle (11pt);   
    		}

    	\foreach \i in {0,...,1}
    	{	   
    		\draw (0,0) -- ({cos(120*\i)},{sin(120*\i)});
    		\filldraw[fill=white] ({cos(120*\i)},{sin(120*\i)}) circle (7pt);  
    	}

    	\foreach \i in {2}
    	{	   
    		\draw (0,0) -- ({1.4*cos(120*\i)},{1.4*sin(120*\i)});
    		\filldraw[fill=white] ({1.4*cos(120*\i)},{1.4*sin(120*\i)}) circle (7pt);  
    	}

    	\filldraw[fill=white] (0,0) circle (8pt);  
		
		\draw (-2.1,-2.9) node {\footnotesize \( L \)};
		\draw (-.1,-1.3) node {\footnotesize \( v \)};
		\draw (-.8,-.4) node {\footnotesize \( e_v \)};
     \end{tikzpicture}

		\caption{The setting of~\ref{item: case 1} in the proof of Proposition~\ref{proposition: 2/3 trees}}\label{figure: case 1}
		\end{subfigure}
		\hfil
		\begin{subfigure}[t]{.3\textwidth}\centering
			\begin{tikzpicture}[scale=.4]

    	\foreach \i in {0,1}
    		\foreach \j in {-1,1}
    			\foreach \k in {-1,1}
    			{ 
    				\draw[] ({2*cos(120*\i+30*\j)},{2*sin(120*\i+30*\j)}) -- ({3*cos(120*\i+30*\j+15*\k)},{3*sin(120*\i+30*\j+15*\k )});
    				\filldraw[fill=white] ({3*cos(120*\i+30*\j+15*\k)},{3*sin(120*\i+30*\j+15*\k)}) circle (5pt);   
    			}
    	  
    	  \foreach \i in {1}
    		\foreach \j in {-1}
    			\foreach \k in {-1,1}
    			{ 
    				\draw[] ({2*cos(120*\i+30*\j)},{2*sin(120*\i+30*\j)}) -- ({3*cos(120*\i+30*\j+15*\k)},{3*sin(120*\i+30*\j+15*\k )});
    				\filldraw[fill=white] ({3*cos(120*\i+30*\j+15*\k)},{3*sin(120*\i+30*\j+15*\k)}) circle (5pt);   
    			}
    	
    	\foreach \i in {0,1}
    		\foreach \j in {-1}
    			\foreach \k in {1}
    			{ 
    				\draw[red, very thick] ({2*cos(120*\i+30*\j)},{2*sin(120*\i+30*\j)}) -- ({3*cos(120*\i+30*\j+15*\k)},{3*sin(120*\i+30*\j+15*\k )});
    				\filldraw[fill=white] ({3*cos(120*\i+30*\j+15*\k)},{3*sin(120*\i+30*\j+15*\k)}) circle (5pt);   
    			}
    	
    	\foreach \i in {1}
    		\foreach \j in {1}
    			\foreach \k in {-1,1}
    			{ 
    				\draw[red, very thick] ({2*cos(120*\i+30*\j)},{2*sin(120*\i+30*\j)}) -- ({3*cos(120*\i+30*\j+15*\k)},{3*sin(120*\i+30*\j+15*\k )});
    				\filldraw[fill=white] ({3*cos(120*\i+30*\j+15*\k)},{3*sin(120*\i+30*\j+15*\k)}) circle (5pt);   
    			}

    	\foreach \i in {0,1} 
    		\foreach \j in {-1,1}
    		{
    			\draw ({cos(120*\i)},{sin(120*\i)}) -- ({2*cos(120*\i+30*\j)},{2*sin(120*\i+30*\j)});
    			\filldraw[fill=white] ({2*cos(120*\i+30*\j)},{2*sin(120*\i+30*\j)}) circle (6pt);   
    		}
    		
    		\foreach \i in {0} 
    		\foreach \j in {-1}
    		{
    			\draw[red, very thick] ({cos(120*\i)},{sin(120*\i)}) -- ({2*cos(120*\i+30*\j)},{2*sin(120*\i+30*\j)});
    			\filldraw[fill=white] ({2*cos(120*\i+30*\j)},{2*sin(120*\i+30*\j)}) circle (6pt);   
    		}
    		
    		\foreach \i in {1} 
    		\foreach \j in {-1,1}
    		{
    			\draw[red, very thick] ({cos(120*\i)},{sin(120*\i)}) -- ({2*cos(120*\i+30*\j)},{2*sin(120*\i+30*\j)});
    			\filldraw[fill=white] ({2*cos(120*\i+30*\j)},{2*sin(120*\i+30*\j)}) circle (6pt);   
    		}

    	\foreach \i in {2} 
    		\foreach \j in {-1,0,1}
    		{
    			\draw[] ({1.4*cos(120*\i)},{1.4*sin(120*\i)}) -- ({2.9*cos(120*\i+30*\j)},{2.9*sin(120*\i+30*\j)});
    			
    			\filldraw[fill=black] ({2.9*cos(120*\i+30*\j)},{2.9*sin(120*\i+30*\j)}) circle (6pt);   
    			%\draw[red] ({2*cos(120*\i+30*\j)},{2*sin(120*\i+30*\j)}) circle (11pt);   
    		}
    	
    	\foreach \i in {2} 
    		\foreach \j in {0,1}
    		{
    			\draw[] ({1.4*cos(120*\i)},{1.4*sin(120*\i)}) -- ({2.9*cos(120*\i+30*\j)},{2.9*sin(120*\i+30*\j)});
    			
    			\filldraw[fill=black] ({2.9*cos(120*\i+30*\j)},{2.9*sin(120*\i+30*\j)}) circle (6pt);   
    			%\draw[red] ({2*cos(120*\i+30*\j)},{2*sin(120*\i+30*\j)}) circle (11pt);   
    		}

    	\foreach \i in {0,1}
    	{	   
    		\draw[red, very thick] (0,0) -- ({cos(120*\i)},{sin(120*\i)});
    		\filldraw[fill=white] ({cos(120*\i)},{sin(120*\i)}) circle (7pt);  
    	}
    	 
    	\foreach \i in {2}
    	{	   
    		\draw[] (0,0) -- ({1.4*cos(120*\i)},{1.4*sin(120*\i)});
    		\filldraw[fill=white] ({1.4*cos(120*\i)},{1.4*sin(120*\i)}) circle (7pt);  
    	}

    	\filldraw[fill=white] (0,0) circle (8pt);  
		
		\draw (-2.1,-2.9) node {\footnotesize \( L \)};
		\draw (-.1,-1.3) node {\footnotesize \( v \)};
		\draw (-.8,-.4) node {\footnotesize \( e_v \)};
     \end{tikzpicture}

		\caption{The setting of~\ref{item: case 2} in the proof of Proposition~\ref{proposition: 2/3 trees}}\label{figure: case 1}
		\end{subfigure}
		\hfil
		\begin{subfigure}[t]{.3\textwidth}\centering
			\begin{tikzpicture}[scale=.4]

    	\foreach \i in {0,1}
    		\foreach \j in {-1,1}
    			\foreach \k in {-1,1}
    			{ 
    				\draw[] ({2*cos(120*\i+30*\j)},{2*sin(120*\i+30*\j)}) -- ({3*cos(120*\i+30*\j+15*\k)},{3*sin(120*\i+30*\j+15*\k )});
    				\filldraw[fill=white] ({3*cos(120*\i+30*\j+15*\k)},{3*sin(120*\i+30*\j+15*\k)}) circle (5pt);   
    			}
    	  
    	  \foreach \i in {1}
    		\foreach \j in {-1}
    			\foreach \k in {-1,1}
    			{ 
    				\draw[red, very thick] ({2*cos(120*\i+30*\j)},{2*sin(120*\i+30*\j)}) -- ({3*cos(120*\i+30*\j+15*\k)},{3*sin(120*\i+30*\j+15*\k )});
    				\filldraw[fill=white] ({3*cos(120*\i+30*\j+15*\k)},{3*sin(120*\i+30*\j+15*\k)}) circle (5pt);   
    			}

    	\foreach \i in {0,1} 
    		\foreach \j in {-1,1}
    		{
    			\draw ({cos(120*\i)},{sin(120*\i)}) -- ({2*cos(120*\i+30*\j)},{2*sin(120*\i+30*\j)});
    			\filldraw[fill=white] ({2*cos(120*\i+30*\j)},{2*sin(120*\i+30*\j)}) circle (6pt);   
    		}
    	 
    	\foreach \i in {1} 
    		\foreach \j in {-1}
    		{
    			\draw[red, very thick] ({cos(120*\i)},{sin(120*\i)}) -- ({2*cos(120*\i+30*\j)},{2*sin(120*\i+30*\j)});
    			\filldraw[fill=white] ({2*cos(120*\i+30*\j)},{2*sin(120*\i+30*\j)}) circle (6pt);   
    		}

    	\foreach \i in {2} 
    		\foreach \j in {-1,0,1}
    		{
    			\draw[] ({1.4*cos(120*\i)},{1.4*sin(120*\i)}) -- ({2.9*cos(120*\i+30*\j)},{2.9*sin(120*\i+30*\j)});
    			
    			\filldraw[fill=black] ({2.9*cos(120*\i+30*\j)},{2.9*sin(120*\i+30*\j)}) circle (6pt);   
    			%\draw[red] ({2*cos(120*\i+30*\j)},{2*sin(120*\i+30*\j)}) circle (11pt);   
    		}
    	
    	\foreach \i in {2} 
    		\foreach \j in {0}
    		{
    			\draw[red, very thick] ({1.4*cos(120*\i)},{1.4*sin(120*\i)}) -- ({2.9*cos(120*\i+30*\j)},{2.9*sin(120*\i+30*\j)});
    			
    			\filldraw[fill=black] ({2.9*cos(120*\i+30*\j)},{2.9*sin(120*\i+30*\j)}) circle (6pt);   
    			%\draw[red] ({2*cos(120*\i+30*\j)},{2*sin(120*\i+30*\j)}) circle (11pt);   
    		}
    		
    	\foreach \i in {2} 
    		\foreach \j in {1}
    		{
    			\draw ({1.4*cos(120*\i)},{1.4*sin(120*\i)}) -- ({2.9*cos(120*\i+30*\j)},{2.9*sin(120*\i+30*\j)});
    			
    			\filldraw[fill=black] ({2.9*cos(120*\i+30*\j)},{2.9*sin(120*\i+30*\j)}) circle (6pt);   
    			%\draw[red] ({2*cos(120*\i+30*\j)},{2*sin(120*\i+30*\j)}) circle (11pt);   
    		}

    	\foreach \i in {0}
    	{	   
    		\draw (0,0) -- ({cos(120*\i)},{sin(120*\i)});
    		\filldraw[fill=white] ({cos(120*\i)},{sin(120*\i)}) circle (7pt);  
    	}
    	
    	\foreach \i in {1}
    	{	   
    		\draw[red, very thick] (0,0) -- ({cos(120*\i)},{sin(120*\i)});
    		\filldraw[fill=white] ({cos(120*\i)},{sin(120*\i)}) circle (7pt);  
    	}
    	
    	\foreach \i in {2}
    	{	   
    		\draw[red, very thick] (0,0) -- ({1.4*cos(120*\i)},{1.4*sin(120*\i)});
    		\filldraw[fill=white] ({1.4*cos(120*\i)},{1.4*sin(120*\i)}) circle (7pt);  
    	}

    	\filldraw[fill=white] (0,0) circle (8pt);  
		
		\draw (-2.1,-2.9) node {\footnotesize \( L \)};
		\draw (-.1,-1.3) node {\footnotesize \( v \)};
		\draw (-.8,-.4) node {\footnotesize \( e_v \)};
     \end{tikzpicture}

		\caption{The setting of~\ref{item: case 3} in the proof of Proposition~\ref{proposition: 2/3 trees}}\label{figure: case 1}
		\end{subfigure}
		\caption{The settings of~\ref{item: case 1},~\ref{item: case 2}, and~\ref{item: case 3} in the proof of Proposition~\ref{proposition: 2/3 trees}. In the figures, the set \( L \) is drawn in black, and a set \( E \) in red.}\label{figure: cases}
	\end{figure}

	\begin{enumerate}[label=(\arabic*)] 
		
		\item Assume that \( E \subseteq E_L. \) Equivalently, assume that each edge in \( E \) is connected to some vertex in \(L.\) Note that since \( E \) satisfies~\ref{item: hat condition ii}, we have \( |E| \geq 2. \) Consequently, by Lemma~\ref{lemma: proof of case 1}, we have
		\begin{align*}
			&\Bigl[ \frac{d}{dp_{E}} P\bigl( X (\mathcal{B}^-(S) ) \equiv 0 \bigr) \Bigr]_{(p_e) \equiv  1}
			= 
			r^{|\mathcal{B}^-(S)\smallsetminus L_E|}
			\bigl( r  (r-1)^{|E|} + (1-r)  r^{|E|} \bigr) 
			\\&\qquad =  
			r^{|L\smallsetminus L_E|-1}
			\bigl( r  (r-1)^{|E|} + (1-r)  r^{|E|} \bigr) \Bigl[  \frac{d}{dp_{E\smallsetminus E_L}} P\pigl( X \bigl(\mathcal{B}^-(S\smallsetminus L) \bigr)   \pigr) \equiv 0 \bigr) \Bigr]_{(p_e)\equiv  1} 
		\end{align*}\label{item: case 1}
		where the last equality follows from the observation that if \( E \subseteq E_L, \) then \( E \smallsetminus E_L = \emptyset\) and \( L_E = L,\) and hence
		\begin{align*}
			&\Bigl[  \frac{d}{dp_{E\smallsetminus E_L}} P\pigl( X \bigl(\mathcal{B}^-(S\smallsetminus L) \bigr)   \pigr) \equiv 0 \bigr) \Bigr]_{(p_e)\equiv  1} 
			=
			\Bigl[ P\pigl( X \bigl(\mathcal{B}^-(S\smallsetminus L)   \bigr)   \pigr) \equiv 0 \bigr) \Bigr]_{(p_e)\equiv  1} 
			\\&\qquad = r^{|\mathcal{B}^-(S)\smallsetminus L_E|+1}.
		\end{align*}

		\item Assume that \( E\cap E_L = \emptyset. \) Equivalently, since \( E \) satisfies~\ref{item: hat condition ii} by assumption, assume that \( E \) is not connected to \( L \cup \{ v \}.\) Noting that in this case, we have \( E = E\smallsetminus E_L\) and \( L_E = \emptyset, \) it follows from Lemma~\ref{lemma: proof of case 2} that
		\begin{align*}
			&\Bigl[ \frac{d}{dp_{E}} P\bigl( X (\mathcal{B}^-(S) ) \equiv 0 \bigr) \Bigr]_{(p_e)\equiv  1}
			=
			r^{|L\smallsetminus L_E|-1}
			\Bigl[  \frac{d}{dp_{E\smallsetminus E_L}} P\pigl( X \bigl(\mathcal{B}^-(S\smallsetminus L)  \bigr)   \pigr) \equiv 0 \bigr) \Bigr]_{(p_e)\equiv  1} .
		\end{align*}
		\label{item: case 2}

		\item Assume that \( e_v \in E. \) Then, since \( E \) satisfies~\ref{item: hat condition ii} by assumption, we have \( E \cap E_L \neq \emptyset .\) Consequently, by Lemma~\ref{lemma: proof of case 3}, we have 
		 		\begin{align*}
			& \Bigl[ \frac{d}{dp_{E}} P\bigl( X (\mathcal{B}^-(S) ) \equiv 0 \bigr)\Bigr]_{(p_e)\equiv  1}
			\\&\qquad= 
			r^{|L\smallsetminus L_E|-1} 
			r\pigl( (r-1)^{|E \cap E_L|} -r^{|E \cap E_L|}  \pigr) \Bigl[
			\frac{d}{dp_{E\smallsetminus E_L}}   P\pigl( X \bigl(\mathcal{B}^-(S\smallsetminus L)    \bigr)  \equiv 0 \pigr) \Bigr]_{(p_e)\equiv 1}	  .
		\end{align*}\label{item: case 3}
	\end{enumerate}
	Combining these cases, we obtain
	\begin{equation}\label{eq: summary}
			\begin{split}
			&	
		\Bigl[\frac{d}{dp_{E}} P\bigl( X (\mathcal{B}^-(S) ) \equiv 0 \bigr)\Bigr]_{(p_e)\equiv  1}
			\\&\quad= 
			r^{|L\smallsetminus L_E|-1} 
			\Bigl[
			\frac{d}{dp_{E\smallsetminus E_L}}   P\pigl( X \bigl(\mathcal{B}^-(S\smallsetminus L)    \bigr)  \equiv 0 \pigr) \Bigr]_{(p_e)\equiv 1}
			\begin{cases}
				 1 &\text{if } E \cap E_L = \emptyset 
				 \cr 
			r  (r-1)^{|E|} + (1-r)  r^{|E|}  
			&\text{if } E \subseteq E_L
			\cr 
			r(r-1)^{|E \cap E_L|} -r^{|E \cap E_L|}r  &\text{otherwise.}%\text{if } E \cap E_L \neq \emptyset \text{ and }E \not \subseteq E_L   
			\end{cases} 
			\end{split}
	\end{equation} 
	
	We now return to the sum in~\eqref{eq: second line with partition}. 
	To this end, note that there is a natural bijection from the set of partitions \(  \hat {\mathcal{P}}(E(T_{S })) \) to the direct product of \( \hat {\mathcal{P}}(E(T_{S \smallsetminus L})) \) and \( \hat {\mathcal{P}}(E_v) \) (see Figure~\ref{figure: bijection}). For \((F_1,\dots, F_\ell) \in \hat {\mathcal{P}}(E(T_{S \smallsetminus L})) \) and \( (G_1,\dots, G_m) \in \hat P(E_v), \) we will without loss of generality assume that \( e_v \in F_1 \cap G_1. \)
	%
	%
	%Moreover, we can note that a partition in \(  P (E(\{ v \} \cup \mathcal{B}^+(\{ v \}))) \) is in \( \hat P (E(\{ v \} \cup \mathcal{B}^+(\{ v \}))) \) exactly if each partition element has degree at least two. 
	% 
	Using this bijection together with~\eqref{eq: summary}, it follows that
	\begin{align*}
		&
		A(S) = \sum_{(E_1, \dots , E_k) \in \hat P(E(T_S))}
         (-1)^{k-1} r^{-k|\mathcal{B}^-(S)|} \prod_{j=1}^k \Bigl[\frac{d}{dp_{E_j}}P\bigl( X (\mathcal{B}^-(S) ) \equiv 0 \bigr)\Bigr]_{(p_e) \equiv  1} 
         \\&\qquad=
		\sum_{\substack{
		(F_1, \dots , F_\ell) \in \hat P(E(T_{S \smallsetminus L})) \\
        (G_1,\dots, G_m) \in \hat P (E_v) }} (-1)^{(\ell+m-1)-1} r^{-(\ell+m-1)|\mathcal{B}^-(S)|} 
        \Bigl[\frac{d}{dp_{F_1 \cup G_1}}P\bigl( X (\mathcal{B}^-(S) ) \equiv 0 \bigr)\Bigr]_{(p_e) \equiv  1} 
        \\&\qquad\qquad\qquad\qquad\cdot 
        \prod_{i=2}^\ell \Bigl[\frac{d}{dp_{F_i}}P\bigl( X (\mathcal{B}^-(S) ) \equiv 0 \bigr)\Bigr]_{(p_e) \equiv  1} 
        \prod_{j=2}^m \Bigl[\frac{d}{dp_{G_i}}P\bigl( X (\mathcal{B}^-(S) ) \equiv 0 \bigr)\Bigr]_{(p_e) \equiv  1}.
	\end{align*} 
	\begin{figure}[h]
	\begin{subfigure}[t]{.3\textwidth}\centering
			\begin{tikzpicture}[scale=.4]

    	\foreach \i in {0}
    		\foreach \j in {-1}
    			\foreach \k in {-1,1}
    			{ 
    				\draw[red, very thick] ({2*cos(120*\i+35*\j)},{2*sin(120*\i+35*\j)}) -- ({3*cos(120*\i+35*\j+12*\k)},{3*sin(120*\i+35*\j+12*\k )});
    				\filldraw[fill=white] ({3*cos(120*\i+35*\j+12*\k)},{3*sin(120*\i+35*\j+12*\k)}) circle (5pt);   
    			}
    	
    	\foreach \i in {0}
    		\foreach \j in {1}
    			\foreach \k in {-1,1}
    			{ 
    				\draw[Dandelion!90!yellow, very thick] ({2*cos(120*\i+35*\j)},{2*sin(120*\i+35*\j)}) -- ({3*cos(120*\i+35*\j+12*\k)},{3*sin(120*\i+35*\j+12*\k )});
    				\filldraw[fill=white] ({3*cos(120*\i+35*\j+12*\k)},{3*sin(120*\i+35*\j+12*\k)}) circle (5pt);   
    			}

    	\foreach \i in {0}
    		\foreach \j in {0}
    			\foreach \k in {0}
    			{ 
    				\draw[Dandelion!90!yellow, very thick] ({2*cos(120*\i+35*\j)},{2*sin(120*\i+35*\j)}) -- ({3*cos(120*\i+35*\j+12*\k)},{3*sin(120*\i+35*\j+12*\k )});
    				\filldraw[fill=white] ({3*cos(120*\i+35*\j+12*\k)},{3*sin(120*\i+35*\j+12*\k)}) circle (5pt);   
    			}

    	\foreach \i in {1}
    		\foreach \j in {1}
    			\foreach \k in {-1 }
    			{ 
    				\draw[red, very thick] ({2*cos(120*\i+30*\j)},{2*sin(120*\i+30*\j)}) -- ({3*cos(120*\i+30*\j+15*\k)},{3*sin(120*\i+30*\j+15*\k )});
    				\filldraw[fill=white] ({3*cos(120*\i+30*\j+15*\k)},{3*sin(120*\i+30*\j+15*\k)}) circle (5pt);   
    			}

    	\foreach \i in {1}
    		\foreach \j in {1}
    			\foreach \k in {0,1}
    			{ 
    				\draw[green, very thick] ({2*cos(120*\i+30*\j)},{2*sin(120*\i+30*\j)}) -- ({3*cos(120*\i+30*\j+15*\k)},{3*sin(120*\i+30*\j+15*\k )});
    				\filldraw[fill=white] ({3*cos(120*\i+30*\j+15*\k)},{3*sin(120*\i+30*\j+15*\k)}) circle (5pt);   
    			}
    	  
    	  \foreach \i in {1}
    		\foreach \j in {-1}
    			\foreach \k in {-1,1}
    			{ 
    				\draw[red, very thick] ({2*cos(120*\i+30*\j)},{2*sin(120*\i+30*\j)}) -- ({3*cos(120*\i+30*\j+15*\k)},{3*sin(120*\i+30*\j+15*\k )});
    				\filldraw[fill=white] ({3*cos(120*\i+30*\j+15*\k)},{3*sin(120*\i+30*\j+15*\k)}) circle (5pt);   
    			}

    	\foreach \i in {1} 
    		\foreach \j in {-1,1}
    		{
    			\draw[red, very thick] ({cos(120*\i)},{sin(120*\i)}) -- ({2*cos(120*\i+30*\j)},{2*sin(120*\i+30*\j)});
    			\filldraw[fill=white] ({2*cos(120*\i+30*\j)},{2*sin(120*\i+30*\j)}) circle (6pt);   
    		}

    	\foreach \i in {0} 
    		\foreach \j in {-1 }
    		{
    			\draw[red, very thick] ({cos(120*\i)},{sin(120*\i)}) -- ({2*cos(120*\i+35*\j)},{2*sin(120*\i+35*\j)});
    			\filldraw[fill=white] ({2*cos(120*\i+35*\j)},{2*sin(120*\i+35*\j)}) circle (6pt);   
    		}
    	 
    	 \foreach \i in {0} 
    		\foreach \j in {0,1}
    		{
    			\draw[Dandelion!90!yellow, very thick] ({cos(120*\i)},{sin(120*\i)}) -- ({2*cos(120*\i+35*\j)},{2*sin(120*\i+35*\j)});
    			\filldraw[fill=white] ({2*cos(120*\i+35*\j)},{2*sin(120*\i+35*\j)}) circle (6pt);   
    		}

    	\foreach \i in {2} 
    		\foreach \j in {0,1}
    		{
    			\draw[cyan!60, very thick] ({cos(120*\i)},{sin(120*\i)}) -- ({2.5*cos(120*\i+30*\j)},{2.5*sin(120*\i+30*\j)});
    			
    			\filldraw[fill=black] ({2.5*cos(120*\i+30*\j)},{2.5*sin(120*\i+30*\j)}) circle (6pt);   
    			%\draw[red] ({2*cos(120*\i+30*\j)},{2*sin(120*\i+30*\j)}) circle (11pt);   
    		}
    	
    	\foreach \i in {2} 
    		\foreach \j in {-1}
    		{
    			\draw[red, very thick] ({cos(120*\i)},{sin(120*\i)}) -- ({2.5*cos(120*\i+30*\j)},{2.5*sin(120*\i+30*\j)});
    			
    			\filldraw[fill=black] ({2.5*cos(120*\i+30*\j)},{2.5*sin(120*\i+30*\j)}) circle (6pt);   
    			%\draw[red] ({2*cos(120*\i+30*\j)},{2*sin(120*\i+30*\j)}) circle (11pt);   
    		}

    	\foreach \i in {0,...,2}
    	{	   
    		\draw[very thick, red] (0,0) -- ({cos(120*\i)},{sin(120*\i)});
    		\filldraw[fill=white] ({cos(120*\i)},{sin(120*\i)}) circle (7pt);  
    	}

    	\filldraw[fill=white] (0,0) circle (8pt);  
		
		\draw (-1.7,-2.9) node {\footnotesize \( L \)};
		\draw (.3,-1.2) node {\footnotesize \( v \)};
     \end{tikzpicture}     
     
		\caption{A tree \( T_S \) with a partition \( (E_1,E_2,E_3,E_4) \in \hat P(E(T_S)) \) (red, blue, green, and orange).}\label{figure: case 1}
		\end{subfigure}
		\hfil 
		\begin{subfigure}[t]{.3\textwidth}\centering
			\begin{tikzpicture}[scale=.4]

    	\foreach \i in {0}
    		\foreach \j in {-1}
    			\foreach \k in {-1,1}
    			{ 
    				\draw[] ({2*cos(120*\i+35*\j)},{2*sin(120*\i+35*\j)}) -- ({3*cos(120*\i+35*\j+12*\k)},{3*sin(120*\i+35*\j+12*\k )});
    				\filldraw[fill=white] ({3*cos(120*\i+35*\j+12*\k)},{3*sin(120*\i+35*\j+12*\k)}) circle (5pt);   
    			}
    	
    	\foreach \i in {0}
    		\foreach \j in {1}
    			\foreach \k in {-1,1}
    			{ 
    				\draw[] ({2*cos(120*\i+35*\j)},{2*sin(120*\i+35*\j)}) -- ({3*cos(120*\i+35*\j+12*\k)},{3*sin(120*\i+35*\j+12*\k )});
    				\filldraw[fill=white] ({3*cos(120*\i+35*\j+12*\k)},{3*sin(120*\i+35*\j+12*\k)}) circle (5pt);   
    			}

    	\foreach \i in {0}
    		\foreach \j in {0}
    			\foreach \k in {0}
    			{ 
    				\draw[] ({2*cos(120*\i+35*\j)},{2*sin(120*\i+35*\j)}) -- ({3*cos(120*\i+35*\j+12*\k)},{3*sin(120*\i+35*\j+12*\k )});
    				\filldraw[fill=white] ({3*cos(120*\i+35*\j+12*\k)},{3*sin(120*\i+35*\j+12*\k)}) circle (5pt);   
    			}

    	\foreach \i in {1}
    		\foreach \j in {-1}
    			\foreach \k in {-1,1}
    			{ 
    				\draw[] ({2*cos(120*\i+30*\j)},{2*sin(120*\i+30*\j)}) -- ({3*cos(120*\i+30*\j+15*\k)},{3*sin(120*\i+30*\j+15*\k )});
    				\filldraw[fill=white] ({3*cos(120*\i+30*\j+15*\k)},{3*sin(120*\i+30*\j+15*\k)}) circle (5pt);   
    			}

    	\foreach \i in {1}
    		\foreach \j in {1}
    			\foreach \k in {-1,0,1}
    			{ 
    				\draw[] ({2*cos(120*\i+30*\j)},{2*sin(120*\i+30*\j)}) -- ({3*cos(120*\i+30*\j+15*\k)},{3*sin(120*\i+30*\j+15*\k )});
    				\filldraw[fill=white] ({3*cos(120*\i+30*\j+15*\k)},{3*sin(120*\i+30*\j+15*\k)}) circle (5pt);   
    			}
    	  
    	  \foreach \i in {1}
    		\foreach \j in {-1}
    			\foreach \k in {-1,1}
    			{ 
    				\draw[] ({2*cos(120*\i+30*\j)},{2*sin(120*\i+30*\j)}) -- ({3*cos(120*\i+30*\j+15*\k)},{3*sin(120*\i+30*\j+15*\k )});
    				\filldraw[fill=white] ({3*cos(120*\i+30*\j+15*\k)},{3*sin(120*\i+30*\j+15*\k)}) circle (5pt);   
    			}

    	\foreach \i in {1} 
    		\foreach \j in {-1,1}
    		{
    			\draw ({cos(120*\i)},{sin(120*\i)}) -- ({2*cos(120*\i+30*\j)},{2*sin(120*\i+30*\j)});
    			\filldraw[fill=white] ({2*cos(120*\i+30*\j)},{2*sin(120*\i+30*\j)}) circle (6pt);   
    		}

    	\foreach \i in {0} 
    		\foreach \j in {-1,0,1}
    		{
    			\draw ({cos(120*\i)},{sin(120*\i)}) -- ({2*cos(120*\i+35*\j)},{2*sin(120*\i+35*\j)});
    			\filldraw[fill=white] ({2*cos(120*\i+35*\j)},{2*sin(120*\i+35*\j)}) circle (6pt);   
    		}

    	\foreach \i in {2} 
    		\foreach \j in {-1,0,1}
    		{
    			\draw[] ({cos(120*\i)},{sin(120*\i)}) -- ({2.5*cos(120*\i+30*\j)},{2.5*sin(120*\i+30*\j)});
    			
    			\filldraw[fill=black] ({2.5*cos(120*\i+30*\j)},{2.5*sin(120*\i+30*\j)}) circle (6pt);   
    			%\draw[red] ({2*cos(120*\i+30*\j)},{2*sin(120*\i+30*\j)}) circle (11pt);   
    		}
    	
    	\foreach \i in {2} 
    		\foreach \j in {0,1}
    		{
    			\draw[cyan!60, very thick] ({cos(120*\i)},{sin(120*\i)}) -- ({2.5*cos(120*\i+30*\j)},{2.5*sin(120*\i+30*\j)}); 
    			
    			\filldraw[fill=black] ({2.5*cos(120*\i+30*\j)},{2.5*sin(120*\i+30*\j)}) circle (6pt);   
    			%\draw[red] ({2*cos(120*\i+30*\j)},{2*sin(120*\i+30*\j)}) circle (11pt);   
    		}
    	
    	\foreach \i in {2} 
    		\foreach \j in {-1}
    		{
    			\draw[red, very thick] ({cos(120*\i)},{sin(120*\i)}) -- ({2.5*cos(120*\i+30*\j)},{2.5*sin(120*\i+30*\j)});
    			
    			\filldraw[fill=black] ({2.5*cos(120*\i+30*\j)},{2.5*sin(120*\i+30*\j)}) circle (6pt);   
    			%\draw[red] ({2*cos(120*\i+30*\j)},{2*sin(120*\i+30*\j)}) circle (11pt);   
    		}

    	\foreach \i in {0,...,2}
    	{	   
    		\draw (0,0) -- ({cos(120*\i)},{sin(120*\i)});
    		\filldraw[fill=white] ({cos(120*\i)},{sin(120*\i)}) circle (7pt);  
    	}
    	 
    	\foreach \i in {2}
    	{	   
    		\draw[very thick, red] (0,0) -- ({cos(120*\i)},{sin(120*\i)});
    		\filldraw[fill=white] ({cos(120*\i)},{sin(120*\i)}) circle (7pt);  
    	}
    	
    	\filldraw[fill=white] (0,0) circle (8pt);  
		
		\draw (-1.7,-2.9) node {\footnotesize \( L \)};
		\draw (.3,-1.2) node {\footnotesize \( v \)};
     \end{tikzpicture}

		\caption{A tree \( T_S \) and a {partition} \( (E_1,E_2) \in \hat P(E_v) \) (red and blue).}\label{figure: case 1}
		\end{subfigure}
		\hfil
		\begin{subfigure}[t]{.3\textwidth}\centering
		
		\begin{tikzpicture}[scale=.4]

    	\foreach \i in {0}
    		\foreach \j in {-1}
    			\foreach \k in {-1,1}
    			{ 
    				\draw[red, very thick] ({2*cos(120*\i+35*\j)},{2*sin(120*\i+35*\j)}) -- ({3*cos(120*\i+35*\j+12*\k)},{3*sin(120*\i+35*\j+12*\k )});
    				\filldraw[fill=white] ({3*cos(120*\i+35*\j+12*\k)},{3*sin(120*\i+35*\j+12*\k)}) circle (5pt);   
    			}
    	
    	\foreach \i in {0}
    		\foreach \j in {1}
    			\foreach \k in {-1,1}
    			{ 
    				\draw[Dandelion!90!yellow, very thick] ({2*cos(120*\i+35*\j)},{2*sin(120*\i+35*\j)}) -- ({3*cos(120*\i+35*\j+12*\k)},{3*sin(120*\i+35*\j+12*\k )});
    				\filldraw[fill=white] ({3*cos(120*\i+35*\j+12*\k)},{3*sin(120*\i+35*\j+12*\k)}) circle (5pt);   
    			}

    	\foreach \i in {0}
    		\foreach \j in {0}
    			\foreach \k in {0}
    			{ 
    				\draw[Dandelion!90!yellow, very thick] ({2*cos(120*\i+35*\j)},{2*sin(120*\i+35*\j)}) -- ({3*cos(120*\i+35*\j+12*\k)},{3*sin(120*\i+35*\j+12*\k )});
    				\filldraw[fill=white] ({3*cos(120*\i+35*\j+12*\k)},{3*sin(120*\i+35*\j+12*\k)}) circle (5pt);   
    			}

    	\foreach \i in {1}
    		\foreach \j in {1}
    			\foreach \k in {-1 }
    			{ 
    				\draw[red, very thick] ({2*cos(120*\i+30*\j)},{2*sin(120*\i+30*\j)}) -- ({3*cos(120*\i+30*\j+15*\k)},{3*sin(120*\i+30*\j+15*\k )});
    				\filldraw[fill=white] ({3*cos(120*\i+30*\j+15*\k)},{3*sin(120*\i+30*\j+15*\k)}) circle (5pt);   
    			}

    	\foreach \i in {1}
    		\foreach \j in {1}
    			\foreach \k in {0,1}
    			{ 
    				\draw[green, very thick] ({2*cos(120*\i+30*\j)},{2*sin(120*\i+30*\j)}) -- ({3*cos(120*\i+30*\j+15*\k)},{3*sin(120*\i+30*\j+15*\k )});
    				\filldraw[fill=white] ({3*cos(120*\i+30*\j+15*\k)},{3*sin(120*\i+30*\j+15*\k)}) circle (5pt);   
    			}
    	  
    	  \foreach \i in {1}
    		\foreach \j in {-1}
    			\foreach \k in {-1,1}
    			{ 
    				\draw[red, very thick] ({2*cos(120*\i+30*\j)},{2*sin(120*\i+30*\j)}) -- ({3*cos(120*\i+30*\j+15*\k)},{3*sin(120*\i+30*\j+15*\k )});
    				\filldraw[fill=white] ({3*cos(120*\i+30*\j+15*\k)},{3*sin(120*\i+30*\j+15*\k)}) circle (5pt);   
    			}

    	\foreach \i in {1} 
    		\foreach \j in {-1,1}
    		{
    			\draw[red, very thick] ({cos(120*\i)},{sin(120*\i)}) -- ({2*cos(120*\i+30*\j)},{2*sin(120*\i+30*\j)});
    			\filldraw[fill=white] ({2*cos(120*\i+30*\j)},{2*sin(120*\i+30*\j)}) circle (6pt);   
    		}

    	\foreach \i in {0} 
    		\foreach \j in {-1 }
    		{
    			\draw[red, very thick] ({cos(120*\i)},{sin(120*\i)}) -- ({2*cos(120*\i+35*\j)},{2*sin(120*\i+35*\j)});
    			\filldraw[fill=white] ({2*cos(120*\i+35*\j)},{2*sin(120*\i+35*\j)}) circle (6pt);   
    		}
    	 
    	 \foreach \i in {0} 
    		\foreach \j in {0,1}
    		{
    			\draw[Dandelion!90!yellow, very thick] ({cos(120*\i)},{sin(120*\i)}) -- ({2*cos(120*\i+35*\j)},{2*sin(120*\i+35*\j)});
    			\filldraw[fill=white] ({2*cos(120*\i+35*\j)},{2*sin(120*\i+35*\j)}) circle (6pt);   
    		}

    	\foreach \i in {2} 
    		\foreach \j in {0,1}
    		{
    			\draw[] ({cos(120*\i)},{sin(120*\i)}) -- ({2.5*cos(120*\i+30*\j)},{2.5*sin(120*\i+30*\j)});
    			
    			\filldraw[fill=black] ({2.5*cos(120*\i+30*\j)},{2.5*sin(120*\i+30*\j)}) circle (6pt);   
    			%\draw[red] ({2*cos(120*\i+30*\j)},{2*sin(120*\i+30*\j)}) circle (11pt);   
    		}
    	
    	\foreach \i in {2} 
    		\foreach \j in {-1}
    		{
    			\draw[] ({cos(120*\i)},{sin(120*\i)}) -- ({2.5*cos(120*\i+30*\j)},{2.5*sin(120*\i+30*\j)});
    			
    			\filldraw[fill=black] ({2.5*cos(120*\i+30*\j)},{2.5*sin(120*\i+30*\j)}) circle (6pt);   
    			%\draw[red] ({2*cos(120*\i+30*\j)},{2*sin(120*\i+30*\j)}) circle (11pt);   
    		}

    	\foreach \i in {0,...,2}
    	{	   
    		\draw[very thick, red] (0,0) -- ({cos(120*\i)},{sin(120*\i)});
    		\filldraw[fill=white] ({cos(120*\i)},{sin(120*\i)}) circle (7pt);  
    	}

    	\filldraw[fill=white] (0,0) circle (8pt);  
		
		\draw (-1.7,-2.9) node {\footnotesize \( L \)};
		\draw (.3,-1.2) node {\footnotesize \( v \)};
     \end{tikzpicture}

		\caption{A tree \( T_S \) with a partition \( (E_1,E_2,E_3) \in \hat P(E(T_{S \smallsetminus L})) \) (red, green, and orange).}\label{figure: case 1}
		\end{subfigure}
		\caption{The setting of the end of the proof of~Proposition~\ref{proposition: 2/3 trees}, illustrating the bijection between \( \hat P(E(T_S)) \) and the direct product of \( \hat P(E(T_{S\smallsetminus L})) \) and \( \hat P(E_v). \)}\label{figure: bijection}
	\end{figure}
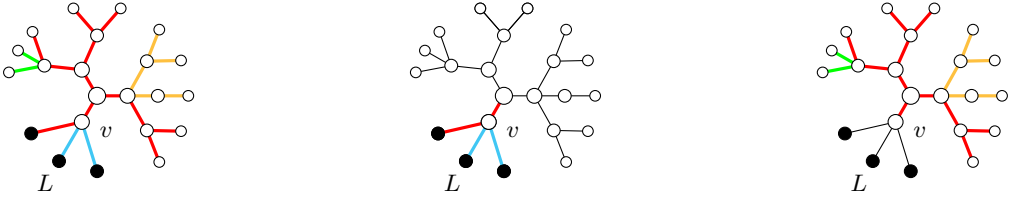
	
	For \( i = 2,3,\dots, \ell, \) using the first case of~\eqref{eq: summary}, we have
	\begin{align*}
		&\Bigl[\frac{d}{dp_{F_i}}P\bigl( X (\mathcal{B}^-(S) ) \equiv 0 \bigr)\Bigr]_{(p_e) \equiv  1} 
		= 
		r^{|L|-1} \Bigl[\frac{d}{dp_{F_i}}P\bigl( X (\mathcal{B}^-(S \smallsetminus L)) \equiv 0   \bigr)\Bigr]_{(p_e) \equiv  1} ,
	\end{align*}
	for \( j = 2,3,\dots, m, \) using the second case of~\eqref{eq: summary}, we have 
	\begin{align*}
		\Bigl[\frac{d}{dp_{G_j}}P\bigl( X (\mathcal{B}^-(S) ) \equiv 0 \bigr)\Bigr]_{(p_e) \equiv  1} = r^{|\mathcal{B}^-(S)|-|L_{G_j}|}
			\Bigl( r  (-1+ r)^{|G_j|} + (1-r)  r^{|G_j|} \Bigr)  ,
	\end{align*}
	and finally, using the third case of~\eqref{eq: summary}, since \(  L_{F_1 \cup G_1} = L_{G_1}\), \( (F_1 \cup G_1) \cap E_L = G_1 \smallsetminus \{ e_v\} \), and \( (F_1 \cup G_1)\smallsetminus E_L = F_1 \), we have
	\begin{align*}
		&\Bigl[\frac{d}{dp_{F_1 \cup G_1}}
		P\bigl( X (\mathcal{B}^-(S) ) \equiv 0 \bigr)\Bigr]_{(p_e) \equiv  1} 
		\\&\qquad= 
		r^{|L|-| L_{G_1}|} 
			\pigl( (r-1)^{|G_1|-1} -r^{|G_1|-1}  \pigr)
			\pigl[ \frac{d}{dp_{F_1}}   P\pigl( X \bigl(\mathcal{B}^-(S\smallsetminus L)\bigr)  \equiv 0 \pigr) \pigr]_{(p_e)\equiv 1}	.
	\end{align*}
	Combining these equations and noting that \( \sum_{j=1}^m|L_{G_j}| = |L|\), \( (\sum_{j=1}^m |G_j|)-1 = |L|, \) and \( |\mathcal{B}^-(S)|-|L| +1) = |\mathcal{B}^-(S\smallsetminus L)|\), we obtain
	\begin{align*}
		&
		A(S) 
		=
		\sum_{\substack{
		(F_1, \dots , F_\ell) \in \hat P(E(T_{S \smallsetminus L})) \\
        (G_1,\dots, G_m) \in \hat P (E_v) }} (-1)^{(\ell+m-1)-1} r^{-\ell|\mathcal{B}^-(S)|+(\ell-1)(|L|-1)} 
        \cdot 
        \pigl( (r-1)^{|G_1|-1} -r^{|G_1|-1}  \pigr) 
        \\&\qquad\qquad\qquad\qquad\cdot 
        \prod_{i=1}^\ell  \Bigl[\frac{d}{dp_{F_i}}P\bigl( X (\mathcal{B}^-(S \smallsetminus L)) \equiv 0   \bigr)\Bigr]_{(p_e) \equiv  1}  
        \cdot  
            \prod_{j=2}^m 
			\Bigl( r  (r-1)^{|G_j|} - (r-1)  r^{|G_j|} \Bigr)  
		\\&\qquad
		= 
		\sum_{\substack{
		(F_1, \dots , F_\ell) \in \hat P(E(T_{S \smallsetminus L})) \\
        (G_1,\dots, G_m) \in \hat P (E_v) }} (-1)^{\ell-1} r^{-\ell |\mathcal{B}^-(S\smallsetminus L)|  }  
        \cdot 
        \prod_{i=1}^\ell  \Bigl[\frac{d}{dp_{F_i}}P\bigl( X (\mathcal{B}^-(S \smallsetminus L)) \equiv 0   \bigr)\Bigr]_{(p_e) \equiv  1}  
        \\&\qquad\qquad\qquad\qquad\cdot  
            (-r)(r-1)^{m-1}  
            \prod_{j=1}^m  \Bigl(  1 - ((r-1)/r)^{|G_j|-1}  \Bigr) 
            \\&\qquad= 
            A(S\smallsetminus L) \cdot 
            \frac{r}{1-r} \sum_{\substack{
		(G_1,\dots, G_m) \in \hat P (E_v) }}  
            (r-1)^{m}  
            \prod_{j=1}^m  \Bigl(  1 - \bigl((r-1)/r \bigr)^{|G_j|-1}  \Bigr) .
	\end{align*} 
	Consequently, if we for \( j \geq 2, \) let
	\begin{align*}
		f_j(r) \coloneqq \frac{r}{1-r} \sum_{\substack{(F_1,\dots, F_m) \in \mathcal{P}([j]) \colon\\ |F_1|,\dots, |F_m|\geq 2}}   (r-1)^{m} 
			\prod_{i =1}^m  
			\Bigl( 1   -   ( (r-1)/r)^{|F_i|-1}  \Bigr),
	\end{align*} 
	then we have shown that
	\begin{equation}\label{eq: induction}
		A(S) = A(S\smallsetminus L) f_{\deg v}(r). 
	\end{equation}
	%
	%
	%Now recall that, by assumption, for all \(  v \in \mathcal{B}^-(S) \)  we have \(  |N(v) \cap S|=1 .\) By the choice of \( L , \) if \( |S\smallsetminus L| \leq 2, \) it follows that for all \(  v \in \mathcal{B}^-(S\smallsetminus L), \) we have  \(  |N(v) \cap (S \smallsetminus L)|=1 .\)  
	%
	Now note that if \( S_0 \) is a connected set with \( |S_0|=2\), then
	\[
	A(S_0)=  
         \Bigl[ \frac{d}{dp_{E(T_{S_0})}}  \log P\bigl( X (S_0 ) \equiv 0 \bigr) \Bigr]_{(p_e)\equiv  1} = \frac{r^2-r}{r^2} = - \frac{1-r}{r}.
	\]
	Combining this observation with~\eqref{eq: induction}, we obtain
	\begin{align*}
		A(S) = -\frac{1-r}{r} \prod_{v \in S \smallsetminus L} f_{ \deg v}(r).
	\end{align*}
	This concludes the proof.
	\end{proof}

We are now finally ready to give a proof of Proposition~\ref{proposition: new negative lemma for p close to one}.

\begin{proof}[Proof of Proposition~\ref{proposition: new negative lemma for p close to one}]
	To simplify notation, let \( E_S \) be the set of edges between \( S \) and \( S^c \), and note that since \( |S|\geq 2 \) and \( S \) is connected, by Lemma~\ref{lemma: trees}, we have 
	\begin{equation*} \label{eq: step 1 for sets}
        \nu(S) =  
         \sum_{J \subseteq \mathcal{B}^-(S)} (-1)^{|J|} 
         \log P\pigl( X \bigl(J \cup (\mathcal{B}^+(S)\smallsetminus \partial_S J) \bigr) \equiv 0 \pigr) .
    \end{equation*} 
    
    Now note that \(  |\partial_{S} \{ j \}|\) is equal to  the number of children in \( \mathcal{B}^+(S) \) of \( j \in \mathcal{B}^-(T_S), \)  Thus, by definition, for any \( J \subseteq \mathcal{B}^-(S) \) we have
    \[
 		\Bigl[ P\pigl( X \bigl(J \cup (\mathcal{B}^+(S)\smallsetminus \partial_S J) \bigr) \equiv 0 \pigr) \Bigr]_{(p_e) \equiv 1} = r^{|J| + \sum_{j \in J} |\partial_{S} \{ j \}|}.
 	\]
 	Using Lemma~\ref{lemma: trees}, we thus obtain
 	\begin{align*}
 		&\bigl[\nu(S) \bigr]_{(p_e) \equiv 1}  
 		= 
 		\Bigr[\sum_{J \subseteq \mathcal{B}^-(S)} (-1)^{|J|} 
         \log P\pigl( X \bigl(J \cup (\mathcal{B}^+(S)\smallsetminus \partial_S J) \bigr) \equiv 0 \pigr) \Bigr]_{(p_e) \equiv 1}
        \\&\qquad =
 		\sum_{J \subseteq \mathcal{B}^-(S)} (-1)^{|J|} 
         \log \bigl( r^{|J| + \sum_{j \in J} |\partial_{S} \{ j \}|} \bigr)
        \\&\qquad =
 		\log r \sum_{J \subseteq \mathcal{B}^-(S)} (-1)^{|J|} 
         \bigl( |J| + \sum_{j \in J} |\partial_{S} \{ j \}| \bigr)
        \\&\qquad =
 		\log r \sum_{J \subseteq \mathcal{B}^-(S)} (-1)^{|J|} |J| 
         + 
 		\log r \sum_{J \subseteq \mathcal{B}^-(S)} (-1)^{|J|} 
         \sum_{j \in J} |\partial_{S} \{ j \}| 
        \\&\qquad =
 		\log r \sum_{J \subseteq \mathcal{B}^-(S)} (-1)^{|J|} |J| 
         + 
 		\log r  \sum_{j \in \mathcal{B}^-(S) } |\partial_{S} \{ j \}|  \sum_{J \subseteq \mathcal{B}^-(S) \colon j \in J} (-1)^{|J|} = 0,
 	\end{align*}
 	since both sums in the last sum of the previous equation are equal to zero if \( |\mathcal{B}^-(S)| \geq 2. \) This completes the proof of~\ref{item: new tree lemma 2} in the case~\( {E=\emptyset.}\)
 	
 	Next, note that by Lemma~\ref{lemma: first step of proof}, for any \( E \subseteq E(T),\) we have
	\begin{align*}
		&\Bigl[ \frac{d}{dp_E} \nu(S) \Bigr]_{(p_e) \equiv 1}  =
		\begin{cases}
			0 &\text{if }  E(T_S) \nsubseteq E \cr
          (-1)^{|\mathcal{B}^-(S) |} \Bigl[\frac{d}{dp_E}    
         \log P\pigl( X \bigl(\mathcal{B}^-(S) \bigr)  \equiv 0 \pigr) \Bigr]_{(p_e) \equiv 1} &\text{if } E=E(T_S)
		\end{cases}
	\end{align*}  
	By Proposition~\ref{proposition: 2/3 trees}, we have
	\begin{align*}
		&\Bigl[ \frac{d}{dp_{E(T_S)}}  
         \log P\pigl( X \bigl(\mathcal{B}^-(S) \bigr) \equiv 0 \pigr) \Bigr]_{(p_e) \equiv  1}
         =
         -\frac{1-r}{r} \prod_{v \in S \smallsetminus L} f_{ \deg v}(r),
	\end{align*} 
	where \( f_j \) is defined in~\eqref{eq: fdef}.  
	Now note that for a tree \( T_0 \) with exactly one vertex \( v \) with \( \deg v  >1 \) and \( S_0=V(T_0), \) we have 
	\begin{align*}
		&\Bigl[ \frac{d}{dp_{E(T_S)}}  
         \log P\pigl( X \bigl(\mathcal{B}^-(S) \bigr) \equiv 0 \pigr) \Bigr]_{(p_e) \equiv  1}
         =
         -\frac{1-r}{r} f_{ \deg v}(r) .
	\end{align*}  
	Using Lemma~\ref{lemma: polylog key} and noting that \( \deg v = |\mathcal{B}^-(S_0) |, \)  it thus follows that 
		\begin{align*}
		&
		f_{ \deg v}(r) 
         = 
         \frac{r}{1-r} \Li_{1-|\mathcal{B}^-(S) |}(-r^{-1}(1-r)),
	\end{align*}  
	where for \( j \geq 1, \) \( \mathrm{Li}_{1-j} \) is the polylogarithm function with index \( 1-j .\)  
	Combining the previous equations, it follows that if \( E = E(T_S), \) then 
	\begin{align*}
		&\Bigl[ \frac{d}{dp_E} \nu(S) \Bigr]_{(p_e) \equiv 1}  
		= 
          (-1)^{|\mathcal{B}^-(S) |} \Bigl[\frac{d}{dp_E}    
         \log P\pigl( X \bigl(\mathcal{B}^-(S) \bigr)  \equiv 0 \pigr) \Bigr]_{(p_e) \equiv 1}  		
         \\&\qquad = 
          (-1)^{|\mathcal{B}^-(S) |+1} 
          \frac{1-r}{r}  \prod_{v \in S \smallsetminus L} f_{ \deg v}(r) 		
         \\&\qquad = 
          (-1)^{|\mathcal{B}^-(S) |+1} 
          \frac{1-r}{r} \prod_{v \in S \smallsetminus L} \Bigl(  \frac{r}{1-r} \Li_{1-\deg v}(-r^{-1}(1-r)) \Bigr) 
         \\&\qquad = 
          (-1)^{|\mathcal{B}^-(S) |+1} 
          \Bigl( \frac{r}{1-r} \Bigr)^{|V(T_S)| - |\mathcal{B}^-(S)|-1}\prod_{j=2}^\infty   \Li_{1-j}\bigl(-r^{-1}(1-r)\bigr)^{k_j} ,
	\end{align*}   
	where for \( j \geq 2, \)  \( k_j \) is the number of vertices of degree \( j\) in \( T_S,\) This concludes the proof. 
\end{proof}

\section{Octopus trees}\label{section: octopus trees}

In this section, the main goal will be to provide a proof of  Theorem~\ref{theorem: general octopus trees}.

\begin{figure}[h]
	
\centering

    \begin{subfigure}[t]{0.5\textwidth}\centering
\begin{tikzpicture}[scale=1.2]
	
	\draw[dotted] (0,-2) -- (0,-2.5);
	\draw[dotted] (0.866025*2,0.5*2) -- (0.866025*2.5,0.5*2.5);
	\draw[dotted] (-0.866025*2,0.5*2) -- (-0.866025*2.5,0.5*2.5);
	
	\draw[help lines] (0,0) -- (0,-2);
	\draw[help lines] (0,0) -- (0.866025*2,0.5*2);
	\draw[help lines] (0,0) -- (-0.866025*2,0.5*2);

	\draw (0,0) node[anchor=north west, xshift=1mm] {\( r_o \)};
	
	\draw (0,-1) node[anchor= west, xshift=1mm] {\( v_{2,1} \)};
	\draw (0,-2) node[anchor= west, xshift=1mm] {\( v_{2,2} \)};
	
	\draw (0.866025*1,0.5*1) node[anchor=north west] {\( v_{1,1} \)};
	\draw (0.866025*2,0.5*2) node[anchor=north west] {\( v_{1,2} \)};
	
	\draw (-0.866025*1,0.5*1) node[anchor=south west] {\( v_{0,1} \)};
	\draw (-0.866025*2,0.5*2) node[anchor=south west] {\( v_{0,2} \)};

	\foreach \x in {0,...,2}
		{\filldraw[fill=white] (0,-\x) circle (3pt);
		\filldraw[fill=white] (0.866025*\x,.5*\x) circle (3pt);
		\filldraw[fill=white] (-0.866025*\x,.5*\x) circle (3pt);}  
	 
\end{tikzpicture}
\end{subfigure}%
    \begin{subfigure}[t]{0.5\textwidth}\centering
\begin{tikzpicture}[scale=1.2]

	\draw[dotted] (0,-2) -- (0,-2.5);
	\draw[dotted] (0.866025*2,0.5*2) -- (0.866025*2.5,0.5*2.5);
	\draw[dotted] (-0.866025*2,0.5*2) -- (-0.866025*2.5,0.5*2.5);
	
	\draw (0,0) -- (0,-2);
	\draw (0,0) -- (0.866025*2,0.5*2);
	\draw (0,0) -- (-0.866025*2,0.5*2);

	\draw (0,-0.5) node[anchor = west] {\( e_{2,1} \)};
	\draw (0,-1.5) node[anchor = west] {\( e_{2,2} \)};
	
	\draw (0.866025*.5,0.5*.5) node[anchor=north west] {\( e_{1,1} \)};
	\draw (0.866025*1.5,0.5*1.5) node[anchor=north west] {\( e_{1,2} \)};
	
	\draw (-0.866025*.5,0.5*.5) node[anchor= south, xshift=2mm] {\( e_{0,1} \)};
	\draw (-0.866025*1.5,0.5*1.5) node[anchor= south, xshift=2mm] {\( e_{0,2} \)};

	\foreach \x in {0,...,2}
		{\filldraw[help lines, fill=white] (0,-\x) circle (3pt);
		\filldraw[help lines, fill=white] (0.866025*\x,.5*\x) circle (3pt);
		\filldraw[help lines, fill=white] (-0.866025*\x,.5*\x) circle (3pt);}  
	
	\draw[red, opacity=0] (-2.2,-2.2) -- (2.2,1.6);
	
\end{tikzpicture}
\end{subfigure}
\caption{The tree  \( T_3\) and the notation of Proposition~\ref{proposition: proposition at zero} and Proposition~\ref{proposition: general small 3-tree}.}\label{figure: 3 tree labelled}
\end{figure}
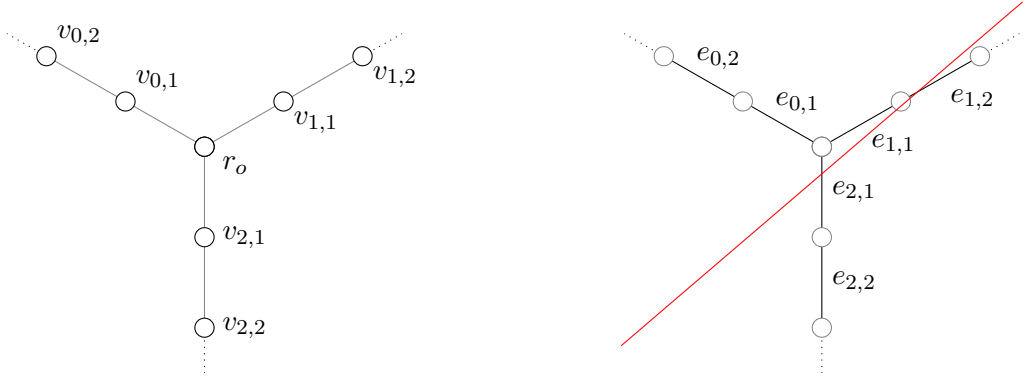

\subsection{\( r \) close to zero}

The main tool in the proof of Theorem~\ref{theorem: general octopus trees}\ref{item: general octopus trees a} is Proposition~\ref{proposition: octopus part 1}, which we now state and prove.
 
\begin{proposition}\label{proposition: octopus part 1}
	    Let \( m \geq 3\) and \( r,p \in (0,1), \) and let \( X\) be a tree-indexed Markov chain  on \( T_m\) with parameters \( (r,p).\) 
	    Then, if \( r<r_1(m)  \in (0,1), \) and \( (p_e)_{e \in E(T)} \) is sufficiently close to 1, we have  \( X \notin \mathcal{R}. \) 
       % \item \( X \in \mathcal{R}\) if \( r \geq r_2(m).\) \label{item: general octopus trees b} 
\end{proposition}

We note that Proposition~\ref{proposition: octopus part 1} is equivalent to~\cite[Theorem 6.1]{fgs}, but since the proof below is different from the corresponding proof in~\cite{fgs}, we include the argument here anyway for completeness.

\begin{proof}[Proof of Proposition~\ref{proposition: octopus part 1}]

	Let \( m' \in  \{ 3,4, \dots, m \} \) be such that \( \Li_{1-m'}(-(1-r)/r)>0. \)
	Since \( r<r_1(m), \)  Lemma~\ref{lemma: polylog}\ref{lemma: polylog b} guarantees that such \( m' \) exists	
	
	    Let \( S = \{ o \} \cup \mathcal{N}(o), \)  let \( T \) be the restriction of \( T_{m} \) to \( V(T_{m'}) \cap S \cup \mathcal{B}^+(S) , \) and let \( \nu_T \) be the unique signed measure corresponding to \( X(V(T)). \)

When \( (p_e) \equiv 1, \) \( (X(v))_{v \in V(T)} \) are independent random variables, and hence we have
\begin{align*}
	\bigl[ \nu_T(\{ v \}) \bigr]_{p= 1} = - \log P(X(v) = 0) = - \log r > 0 \quad \forall  v \in S,
\end{align*}
 and moreover, \( \bigl[ \nu_T(S) \bigr]_{p = 1}= 0 \) for all \( S \in \mathcal{P}(S)\smallsetminus \{ \emptyset \} \) with \( |S| \neq 1. \) Consequently, by Proposition~\ref{proposition: new negative lemma for p close to one} and a Taylor expansion, for \( (p_e)_{e \in E(T)} \) sufficiently close to one, \( \nu_T(S) \) has the same sign as
\begin{align*}
	(-1)^{|E(T_S)|}
			(-1)^{|\mathcal{B}^-(S) |+1}  
			 \Li_{1-m'}\bigl(-r^{-1}(1-r)\bigr)  = -\Li_{1-m'}\bigl(-r^{-1}(1-r)\bigr).
\end{align*}
Since \( -\Li_{1-m'}\bigl(-r^{-1}(1-r)<0, \) we must have \( X|_{V(T)} \notin \mathcal{R}. \) Using Lemma~\ref{lemma: finite to infinite}\ref{item: finite to infinite i}, it follows that \( X \notin \mathcal{R}. \) This concludes the proof.
\end{proof}

\begin{proof}[Proof of Theorem~\ref{theorem: general octopus trees}\ref{item: general octopus trees a}]
	Since \( T_m \) is self-similar, the desired conclusion follows by combining Proposition~\ref{proposition: octopus part 1} and  Proposition~\ref{proposition: scaling} that \( X \notin \mathcal{R} \) for any \( p \in (0,1). \) This concludes the proof.  
\end{proof}

\subsection{\( r \) close to one}

In this section, we provide a proof of Theorem~\ref{theorem: general octopus trees}\ref{item: general octopus trees b}. The main tool in this proof is the following proposition.

\begin{proposition}\label{proposition: proposition at zero}
	Let \( m \geq 3 \) and let \( T=T_m .\) For \( v \in V(T) \) and \( e \in E(T), \) let  \( r_v \in (0,1)\) and \( p_e \in (0,1). \)  Let \( X \) be the tree-indexed Markov chain on \( T \) with parameters \( ((r_v)_{v \in V},(p_e)_{e \in E}). \)
	Let \( \nu \) be the corresponding signed measure in the sense of Lemma~\ref{lemma: unique signed}, and let \( S \coloneqq \{ o \} \cup N(\{ o \}).\) Let the vertices and edges of \( T \) be  labeled as in Figure~\ref{figure: 3 tree labelled}.
	Then the following holds.
	\begin{enumerate}[label=(\alph*)]
		\item \( \bigl[ \nu(S) \bigr]_{r_o = 1}=0. \) \label{item: r der 1}
		\item For any multiset \( K \subseteq V(T) \smallsetminus \{ o \}, \) we have
	\( \bigl[ \frac{d}{dr_K} \nu(S) \bigr]_{(r_v) \equiv 1} =0. \)\label{item: r der 2}
	\item \( \bigl[\frac{d}{dr_{o}} \nu(S) \bigr]_{(r_v) \equiv 1} = -\prod_{j=1}^m (1-p_{e_{j,1}})p_{e_{j,2}} .\) \label{item: r der 3}
	\item There is a constant \( C_m >0\) that depends only on \(m \) such that for all \( v,v'\in V(T)\), we have  \[ 
		\bigl| \frac{d^2}{dr_{v} dr_{v'}} \nu(S) \bigr| \leq \frac{C_m }{\prod_{j=1}^m r_{v_{j,1}}^{2^{m+1}}} \cdot \prod_{j=1}^m (1-p_{e_{j,1}})p_{e_{j,2}} . \] \label{item: r der 4}
	\end{enumerate}
\end{proposition}

Before we provide a proof of Proposition~\ref{proposition: proposition at zero}, we show how we use it to prove Theorem~\ref{theorem: general octopus trees}\ref{item: general octopus trees b}.

\begin{proof}[Proof of Theorem~\ref{theorem: general octopus trees}\ref{item: general octopus trees b}]

	Let \( S \) be any finite subset of \( V(T_m). \) Let \( T \) be any finite subtree of \( T_m \) so that \( S \cup \mathcal{B}^+(S) \subseteq V(T)\). Let \( \nu_T \) be the unique signed measure corresponding to \( X(T) \) in the sense of Lemma~\ref{lemma: unique signed}.

	We first make the following observations.
	\begin{enumerate}[label=(\roman*)]
		\item If \(|S|=1 \) then \( \nu_T(S) \geq 0\) by Lemma~\ref{lemma: trees}.
		\item  If  \( S  \) is not connected, then \( \nu_T(S) = 0 \) by~Proposition~\ref{prop: connected}.
		
		\item If \( S  \)  is a connected set that does not contain the origin \( o, \) then \( S \) is an interval of length \( j \coloneqq |S|. \)  By Lemma~\ref{lemma: trees}, we have
		\begin{align*}
			&\nu_T(S) = \sum_{J \subseteq \mathcal{B}^-(S)} (-1)^{|J|} 
         \log P\pigl( X \bigl(J \cup (\mathcal{B}^+(S)\smallsetminus \partial_S J) \bigr) \equiv 0 \pigr)
         \\&\qquad=
         \log  \bigl( (1-p)^{j+1} r + (1-(1-p)^{j+1})r^2  \bigr) + 
         \log  \bigl( (1-p)^{j-1}r+(1-(1-p)^{j-1})r^2  \bigr)
         \\&\qquad\qquad-2\log  \bigl( (1-p)^{j}r +(1-(1-p)^{j})r^2 \bigr)
         \\&\qquad= 
         \log \Bigl( 1 +\frac{p^2(1 - p)^{j-1} r (1 - r) }{(r+(1 - p)^j  - (1 - p)^j r)^2 }\Bigr) > 0.
		\end{align*} 
	\end{enumerate}
	Now assume that \( |S|>1 \) and \( o \in S. \) 
	Let \( V_1, \dots, V_m \subseteq V(T_m) \) denote the vertex sets corresponding to the connected components of the restriction of \( T_m \) to \( V(T_m)\smallsetminus \{ o \}. \)
	By Proposition~\ref{proposition: proposition at zero} applied with \( p_{e_{j,1}} = 1-(1-p)^{|S \cap V_j|}\) and \( p_{e_{j,2}} = p \) for all \( j \in [m], \) there exists a constant \( r_2(m) \in (0,1), \) independent of \( S, \) \( T ,\) and \( p ,\) such that, by Taylor's theorem, \( \nu_T(S) > 0 \) for all \( r \geq r_2(m). \)

	To finish the proof, note that by Lemma~\ref{lemma: trees}, \( \nu_T(S) \) is independent of the choice of \( T. \)
	For \( n \geq 1, \) let \( T_{m}^{\leq n} \) be the subtree of \( T_m \) which contains all vertices at distance at most \( n \) from the origin. Applying Lemma~\ref{lemma: finite to infinite}\ref{item: finite to infinite ii} with \( S_n = V(T_m) \smallsetminus V(T_m^{\leq n})\), the desired conclusion immediately follows.
\end{proof}

\begin{proof}[Proof of Proposition~\ref{proposition: proposition at zero}]
	For \( j \in [m], \) let \( V_j \coloneqq \{ v_{j,1}, v_{j,2} \}, \) and note that \( S = \{ o \} \cup \{ v_{j,1},\, j \in [m] \}. \) 
	For \( J \subseteq \mathcal{B}^-(S) \) consider the event
	\begin{equation*}
		\mathcal{E}_J \coloneqq \bigl\{ X \bigl(J \cup (\mathcal{B}^+(S)\smallsetminus \partial_S J) \bigr) \equiv 0 \bigr\},
	\end{equation*}
	and for \( J \subseteq \mathcal{B}^-(S) \) and  \( j \in [m] \), consider the event
	\begin{equation*}
		\mathcal{E}_{J}^j \coloneqq \begin{cases}
			X(v_{j,1}) = 0 &\text{if } J \cap V_j \neq \emptyset \cr 
			X(v_{j,2}) = 0 &\text{otherwise.}
		\end{cases}
	\end{equation*}
	
	Since \( |S| \geq 2, \) by Lemma~\ref{lemma: trees}, we have
	\begin{align*} 
        \nu(S) =  
         \sum_{J \subseteq \mathcal{B}^-(S)} (-1)^{|J|} 
         \log P( \mathcal{E}_J ). 
 	\end{align*} 
 	Let \( (R_v)_{v \in T} \) be as in A' and B' of Section~\ref{section: detailed definition}. Then
	 \begin{align*} 
		& \bigl[ \nu(S) \bigr]_{r_o = 1}
		=  
		\Bigl[ \sum_{J \subseteq \mathcal{B}^-(S)} (-1)^{|J|} 
         \log P ( \mathcal{E}_J ) \Bigr]_{r_o = 1}
		=
		\sum_{J \subseteq \mathcal{B}^-(S)} (-1)^{|J|} 
		\log P\bigl(\mathcal{E}_J\mid R_o = 0\bigr) .
	\end{align*}
	Conditioned on \( R_o, \) the random vectors \( X(V_j)\), \( j \in [m] \) are independent. Using this observation, it follows that
	\begin{align*}
		&\sum_{J \subseteq \mathcal{B}^-(S)} (-1)^{|J|} 
		\log P\bigl(\mathcal{E}_J\mid R_o = 0\bigr) 
		=
		\sum_{J \subseteq \mathcal{B}^-(S)} (-1)^{|J|} 
		\log \prod_{j=1}^m P\bigl(\mathcal{E}_{J}^j \mid R_o = 0\bigr) 
		\\&\qquad=
		\sum_{J_1 \subseteq \mathcal{B}^-(S)\cap V_1}  \!\!\! \dots \!\!\! \sum_{J_m \subseteq  \mathcal{B}^-(S)\cap V_m}(-1)^{|J_1|+|J_2| + \dots + |J_m|} 
		\log \prod_{j=1}^m  P(  \mathcal{E}^j_{J_j} \mid R_o = 0 ) 
		\\&\qquad= 
		\sum_{J_1 \subseteq \mathcal{B}^-(S)\cap V_1}  \!\!\! \dots \!\!\! \sum_{J_m \subseteq \mathcal{B}^-(S)\cap V_m}(-1)^{|J_1|+|J_2| + \dots + |J_m|} 
		\sum_{j=1}^m  \log P (\mathcal{E}_{J_j}^j \mid R_o = 0 )
		\\&\qquad= 
		\sum_{j=1}^m  \biggl( \sum_{J_j \subseteq \mathcal{B}^-(S)\cap V_j}  (-1)^{|J_j|}  \log P (\mathcal{E}_{J_j}^j \mid R_o = 0 )    \prod_{k \in [m]\smallsetminus \{ j\}} \sum_{J_k \subseteq \mathcal{B}^-(S)\cap V_k}(-1)^{|J_k|} 
		 \biggr).
	\end{align*}  
	Now note that, for any \( j \in [m] \) and \( k \in  [m ]\smallsetminus \{ j\} \), the set \( \mathcal{B}^-(S)\cap V_k \) is non-empty (for any  \( j \in [m] \), such  \(k \) exist since \( m \geq 2 \)). This implies in particular that for any \( j \in [m] \) and \( k \in [m]\smallsetminus \{ j \}, \) we have \[ \sum_{J_k \subseteq \mathcal{B}^-(S)\cap V_k}(-1)^{|J_k|}  = 0. \]
	Combining the previous equations, we obtain
	\begin{equation*}
		\bigl[ \nu(S) \bigr]_{r_o = 1}=0.
	\end{equation*} 
	This completes the proof of~\ref{item: r der 1}.
	Since differentiation is a linear operator, the same conclusion holds for any combination of derivatives of \( \nu(S) \) that do not include \( \frac{d}{dr_o}. \)
	In particular, for any multiset \( K \subseteq V(T)\smallsetminus \{ o \}, \) we have
	\begin{align*}
		\Bigl[\frac{d}{dr_K} \nu(S) \Bigr]_{r_o = 1} =0.
	\end{align*}
	This completes the proof of~\ref{item: r der 2}.

	We now calculate \( \frac{d}{dr_o} \nu(S)|_{(r_v)  \equiv 1}. \) To this end, note first that for any \( J \subseteq \mathcal{B}^-(S), \) we have
	\begin{align*}
		\frac{d}{dr_o} P ( \mathcal{E}_J )=P ( \mathcal{E}_J \mid R_o = 0 )-P ( \mathcal{E}_J \mid R_o = 1 ).
	\end{align*}
	Using this observation, it follows that
	\begin{align*} 
        & \frac{d}{dr_o} \nu(S) =  
         \frac{d}{dr_o} \sum_{J \subseteq \mathcal{B}^-(S)} (-1)^{|J|} 
         \log P ( \mathcal{E}_J )  
         =
         \sum_{J \subseteq \mathcal{B}^-(S)} (-1)^{|J|} 
         \frac{\frac{d}{dr_o} P ( \mathcal{E}_J ) }{ P (\mathcal{E}_J )  }
         \\&\qquad=
         \sum_{J \subseteq \mathcal{B}^-(S)} (-1)^{|J|} 
         \frac{P ( \mathcal{E}_J \mid R_o = 0 )-P ( \mathcal{E}_J \mid R_o = 1 )   }{P(\mathcal{E}_J)   },
         \end{align*}
         and hence 
         \begin{align*}
         &\Bigl[ \frac{d}{dr_o} \nu(S) \Bigr]_{r_o = 1}
         = 
         \sum_{J \subseteq \mathcal{B}^-(S)} (-1)^{|J|} \Biggl[ 1-
         \frac{P( \mathcal{E}_J \mid R_o = 1 )   }{P ( \mathcal{E}_J\mid R_o = 0 )   }  \Biggr]
         = 
         -\sum_{J \subseteq \mathcal{B}^-(S)} (-1)^{|J|} 
         \frac{P ( \mathcal{E}_J\mid R_o = 1 )   }{P ( \mathcal{E}_J\mid R_o =0 )   }.
	\end{align*}
	Now recall that given \( r_o, \) the random vectors \( X(V_1), \) \dots, \( X(V_m) \) are independent. From this it follows that
	\begin{align*}
         &\sum_{J \subseteq \mathcal{B}^-(S)} (-1)^{|J|} 
         \frac{P ( \mathcal{E}_J \mid R_o = 1 )   }{P ( \mathcal{E}_J  \mid R_o = 0  )   } 
         =
         \sum_{J \subseteq \mathcal{B}^-(S)} (-1)^{|J|} 
         \prod_{j = 0 }^m 
         \frac{P ( \mathcal{E}_{J \cap V_j}^j   \mid R_o = 1\bigr)   }{P\bigl( \mathcal{E}_{J \cap V_j}^j\mid R_o = 0\bigr)   }
         \\&\qquad
         \sum_{J_1 \subseteq \mathcal{B}^-(S) \cap V_1 }  \dots \sum_{J_m \subseteq \mathcal{B}^-(S) \cap V_m } 
         \prod_{j =1}^m (-1)^{|J_j|}
         \frac{P (\mathcal{E}_{J_j}^j  \mid R_o = 1 )   }{P ( \mathcal{E}_{J_j}^j\mid R_o = 0 )   }.
	\end{align*}
	Changing the order of multiplication and summation, we can rewrite the previous expression as
	\begin{align*}
         & \prod_{j =1}^m \sum_{J_j \subseteq \mathcal{B}^-(S)\cap V_j }          
         (-1)^{|J_j|}
         \frac{P (\mathcal{E}_{J_j}^j \mid R_o = 1 )   }{P ( \mathcal{E}_{J_j}^j \mid R_o = 0 )   } .
	\end{align*} 
	Now note that for any \( j \in [m], \) we have
	\begin{align*}
		&  
		\sum_{J_j \subseteq \mathcal{B}^-(S) \cap V_j}          
         (-1)^{|J_j|}
         \frac{P (\mathcal{E}_{J_j}^j  \mid R_o = 1 )   }{P ( \mathcal{E}_{J_j}^j  \mid R_o = 0 )   }
         \\&\qquad=
         \frac{P\bigl( X ( v_{j,2} ) = 0 \mid R_o = 1\bigr)   }{P\bigl( X (  v_{j,2}) = 0 \mid R_o = 0\bigr) }
         -
         \frac{P\bigl( X ( v_{j,1} ) = 0 \mid R_o = 1\bigr)   }{P\bigl( X (  v_{j,1}) = 0 \mid R_o = 0\bigr) }
         \\&\qquad=
         \frac{p_{e_{j,2}}r_{v_{j,2}}+p_{e_{j,1}}(1-p_{e_{j,2}}) r_{v_{j,2}}}{p_{e_{j,2}} r_{v_{j,2}} + p_{e_{j,1}}(1-p_{e_{j,2}}) r_{v_{j,1}} + (1-p_{e_{j,1}})(1-p_{e_{j,2}})} 
         - 
         \frac{p_jr_{v_{j,1}}}{p_{e_{j,1}} r_{v_{j,1}} + (1-p_{e_{j,1}})}
         \\&\qquad= 
         \frac{
         (1-p_{e_{j,1}}) p_{e_{j,2}} r_{v_{j,2}} 
         }{
         (p_{e_{j,2}}r_{v_{j,2}} + p_{e_{j,1}}(1-p_{e_{j,2}}) r_{v_{j,1}} + (1-p_{e_{j,1}})(1-p_{e_{j,2}}))(p_{e_{j,1}} r_{v_{j,1}} + (1-p_{e_{j,1}}))
         }.
	\end{align*}   
	Combining the previous equations, it follows that
	\begin{align*}
		\frac{d}{dr_{o}} \nu(S) \big|_{(r_v)  \equiv 1} = -\prod_{j=1}^m (1-p_{e_{j,1}})p_{e_{j,2}},
	\end{align*}
	thus completing the proof of~\ref{item: r der 3}.
	
	It remains to show that~\ref{item: r der 4} holds. To this end, let  \( K \subseteq \{ o \} \cup \{ v_{j,i}\colon i \in [2],\, j \in [m] \}\) be a non-empty multiset.
	Then, for any \( j \in [m] \) and \( J \subseteq \mathcal{B}^-(S), \) we have
	\begin{equation*}
		\bigl[ P ( \mathcal{E}_J ) \bigr]_{ p_{e_{j,2}} = 0 }
		= \bigl[ P( \mathcal{E}_{J \Delta \{ v_{j,1} \}} )\bigr]_{ p_{e_{j,2}} = 0 }
	\end{equation*}
	and hence
	\begin{align*}
		 &\Bigl[ \frac{d}{dr_K}\nu(S) \Bigr]_{ p_{e_{j,2}} = 0 }
		 = 
		\biggl[\sum_{J \subseteq \mathcal{B}^-(S)} (-1)^{|J|} 
         \frac{d}{dr_{K}} \log P( \mathcal{E}_J )  \biggr]_{ p_{e_{j,2}} = 0 }
         \\&\qquad=
		\biggl[\sum_{J \subseteq \mathcal{B}^-(S) \smallsetminus \{  v_{j,1} \}} (-1)^{|J|} 
         \frac{d}{dr_{K}} \log P( \mathcal{E}_{J} ) - \log P( \mathcal{E}_{J \cup \{  v_{j,1} \}} ) \biggr]_{ p_{e_{j,2}} = 0 }
         \\&\qquad=
		\biggl[\sum_{J \subseteq \mathcal{B}^-(S) \smallsetminus \{  v_{j,1} \}} (-1)^{|J|} 
         \frac{d}{dr_{K}} \log P( \mathcal{E}_{J} ) - \log P( \mathcal{E}_J ) \biggr]_{ p_{e_{j,2}} = 0 }
         = 0.
	\end{align*} 
	 
	Next, note that if \( p_{e_{j,1}} = 1 \) for some \( j \in [m], \) then \( X(V_j) \) and \( X(S\smallsetminus V_j ) \) are independent. Consequently, we have
	\begin{align*} 
        & \Bigl[ \frac{d}{dr_{K}} \nu(S) \Bigr]_{p_{e_{j,1}} = 1 }
        =  
        \biggl[\sum_{J \subseteq \mathcal{B}^-(S)} (-1)^{|J|} \frac{d}{dr_{K}} 
         \log P ( \mathcal{E}_J )  \biggr]_{p_{e_{j,1}} = 1 }
        \\&\qquad=  
        \biggl[ \sum_{J_j \subseteq  \mathcal{B}^-(S) \cap V_j} 
        \sum_{J \subseteq \mathcal{B}^-(S) \smallsetminus V_j} (-1)^{|J_j|+|J|} 
         \frac{d}{dr_{K}}  \log P(\mathcal{E}_{J_j}^j) P\pigl( X \bigl(J \cup (\mathcal{B}^+(S) \smallsetminus (V_j\cup \partial_S J)) \bigr) \equiv 0 \pigr) \biggr]_{p_{e_{j,1}} = 1 }
        \\&\qquad=  
        \biggl[ \sum_{J_j \subseteq  \mathcal{B}^-(S) \cap V_j} 
        \sum_{J \subseteq \mathcal{B}^-(S) \smallsetminus V_j} (-1)^{|J_j|+|J|} 
        \frac{d}{dr_{K}}  \log P(\mathcal{E}_{J_i}^i)    \biggr]_{p_{e_{j,1}} = 1 }      \\&\qquad\qquad+ 
        \biggl[\sum_{J_j \subseteq \mathcal{B}^-(S) \cap V_j} 
        \sum_{J \subseteq \mathcal{B}^-(S) \smallsetminus V_j } (-1)^{|J_j|+|J|} 
         \frac{d}{dr_{K}}  \log P\pigl( X \bigl(J \cup (\mathcal{B}^+(S) \smallsetminus (V_j \cup \partial_S J)) \bigr) \equiv 0 \pigr)\biggr]_{p_{e_{j,1}} = 1 }. 
	\end{align*} 
	Noting that, since \( \mathcal{B}^-(S) \smallsetminus V_j \neq \emptyset \) and \( \mathcal{B}^-(S) \cap V_j \neq \emptyset \), we have
	\begin{align*}
		& \sum_{J_j \subseteq  \mathcal{B}^-(S) \cap V_j} 
        \sum_{J \subseteq \mathcal{B}^-(S) \smallsetminus V_j} (-1)^{|J_j|+|J|} 
        \frac{d}{dr_{K}}  \log P(\mathcal{E}_{J_i}^i)    
        \\&\qquad=
        \sum_{J_j \subseteq  \mathcal{B}^-(S) \cap V_j} 
         (-1)^{|J_j|} 
        \frac{d}{dr_{K}}  \log P(\mathcal{E}_{J_i}^i) \sum_{J \subseteq \mathcal{B}^-(S) \smallsetminus V_j}(-1)^{|J|} =0
	\end{align*}
	and  
	\begin{align*}
		&\sum_{J_j \subseteq \mathcal{B}^-(S) \cap V_j} 
        \sum_{J \subseteq \mathcal{B}^-(S) \smallsetminus V_j } (-1)^{|J_j|+|J|} 
         \frac{d}{dr_{K}}  \log P\pigl( X \bigl(J \cup (\mathcal{B}^+(S) \smallsetminus (V_j \cup \partial_S J)) \bigr) \equiv 0 \pigr)
         \\&\qquad= 
        \sum_{J \subseteq \mathcal{B}^-(S) \smallsetminus V_j } (-1)^{|J|} 
         \frac{d}{dr_{K}}  \log P\pigl( X \bigl(J \cup (\mathcal{B}^+(S) \smallsetminus (V_j \cup \partial_S J)) \bigr) \equiv 0 \pigr) 
         \sum_{J_j \subseteq \mathcal{B}^-(S) \cap V_j} (-1)^{|J_j|} =0
	\end{align*}
	it follows that \( \bigl[ \frac{d}{dr_{K}}  \nu(S) \bigr]_{p_{e_{j,1}}=1} = 0. \) 
	To sum up, we have shown that
	\begin{equation*}
		\bigl[ \frac{d}{dr_K}\nu(S) \bigr]_{ p_{e_{j,2}} = 0 } = \bigl[ \frac{d}{dr_{K}}  \nu(S) \bigr]_{p_{e_{j,1}}=1} = 0  ,\quad \forall j \in [m].
	\end{equation*}

	%Since \( \frac{d}{dr_{K}}  \nu(S) \) is a rational function, it follows that in this case,  \( \frac{d}{dr_{K}} \nu(S) \) is divisible by \( \prod_{j=1}^n (1-p_{e_{j,1}}) . \)
	%
	%
	%
	Now note that for any distinct \( v,v' \in V(T) \) we have  
	\begin{align}
	  	& \frac{d^2}{dr_{v} dr_{v'}} \nu(S)  \nonumber
        =  
        \frac{d^2}{dr_{v} dr_{v'}} 
        \sum_{J \subseteq \mathcal{B}^-(S)} (-1)^{|J|} 
        \log P(\mathcal{E}_J)  
        \\&\qquad=  \nonumber
        \sum_{J \subseteq \mathcal{B}^-(S)} (-1)^{|J|} 
         \frac{
         P(\mathcal{E}_J)  \frac{d^2}{  dr_{v} dr_{v'}} P(\mathcal{E}_J) 
         -
         \frac{d}{dr_{v}  } P(\mathcal{E}_J) 
         \frac{d}{ dr_{v'}} P(\mathcal{E}_J)  }{P(\mathcal{E}_J) ^2}
        \\&\qquad=  
                 \frac{
         \sum_{J \subseteq \mathcal{B}^-(S)} (-1)^{|J|} \bigl( 
         P(\mathcal{E}_J)  \frac{d^2}{dr_{v} dr_{v'}} P(\mathcal{E}_J) 
         -
         \frac{d}{dr_{v}  } P(\mathcal{E}_J) 
         \frac{d}{ dr_{v'}} P(\mathcal{E}_J) 
         \bigr)
         \prod_{I \subseteq \mathcal{B}^-(S),\, I\neq J}  P(\mathcal{E}_I) ^2
         }{ \prod_{J \subseteq \mathcal{B}^-(S)} P(\mathcal{E}_J)^2}\label{eq: last line of expr}
	  \end{align}
	  Note that both the numerator and the denominator are polynomials in   \( (r_{v})_{v \in V(T)}, \) and \( (p_e)_{e \in E(T)}. \) Further, note that by positive association, we have \( P(\mathcal{E}_J) \geq \prod_{j=1}^m r_{v_{j,1}} \) for all \( J \subseteq \mathcal{B}^-(S). \) Since \( \bigl[ \frac{d^2}{dr_{v} dr_{v'}} \nu(S) \bigr]_{ p_{e_{j,2}} = 0 } = 0\) and \( \bigl[ \frac{d^2}{dr_{v} dr_{v'}}  \nu(S) \bigr]_{p_{e_{j,1}}=1} = 0 \) for all \( j \in [m] \), it follows that the numerator of~\eqref{eq: last line of expr} must have \( \prod_{j=1}^m p_{e_{j,2}} (1-p_{e_{j,1}}) \) as a factor. Since each \(  P(\mathcal{E}_J)  \) is a sum of \( 2^{2m+1} \) monomials of degree at most one in each variable, we obtain~\ref{item: r der 4} as desired.
	  This completes the proof.
\end{proof}

\subsection{\( m = 3 \)}

In this section, we provide a proof of Theorem~\ref{theorem: general octopus trees}\ref{item: general octopus trees c}. The main additional tool we will need in this proof is the following proposition.

\begin{proposition}\label{proposition: general small 3-tree}
	Let \( T \) be the restriction of \( T_3 \) to the set of vertices at distance at most two from the root, let the vertices and edges be labeled as in  Figure~\ref{figure: 3 tree labelled}, and let \( S \coloneqq \{ o \} \cup \mathcal{B}^+(o). \) 
	For \( j \in [3] \) and \( k \in [2], \) let \( r_{o} = r_{v_{j,k}} = r \in [\frac{1}{2},1] \) and  \( p_{e_{j,k}} \in [0,1]. \) 
	Further, let \( X \) be the tree-indexed Markov chain on \( T \) with parameters \( ((r_v),(p_e))\), and let \( \nu \) be the corresponding signed measure as in Lemma~\ref{lemma: unique signed}. Then \( \nu(S) \geq 0. \)
\end{proposition}

Before we give a proof of Proposition~\ref{proposition: general small 3-tree}, we show how it implies Theorem~\ref{theorem: general octopus trees}\ref{item: general octopus trees c}.
\begin{proof}[Proof of Theorem~\ref{theorem: general octopus trees}\ref{item: general octopus trees c}]
	Using Theorem~\ref{theorem: general octopus trees}\ref{item: general octopus trees a}, it suffices to show that \( X \in \mathcal{R} \) if \( r \geq 1/2. \)

	Let \( S \subseteq V(T_3)\) be finite and let \( T \) be a finite subtree of \( T_3 \) such that \( S \cup \mathcal{B}^+(S)\subseteq V(T) . \) Let \( \nu_T \) be the signed measure corresponding to \( X\bigl(V(T)\bigr) \) in the sense of Lemma~\ref{lemma: unique signed}. 
	We now show that \( \nu_T(S) \geq 0. \) To this, end, note first that as in the proof of Theorem~\ref{theorem: general octopus trees}\ref{item: general octopus trees b}, we know that \( \nu_T(S) \geq 0 \) if (i) \( |S|=1 \), (ii) \( S\) is not connected, or (iii) \( o \notin S. \)
	%
	%\begin{enumerate}[label=(\roman*)]
	%	\item If \(|S|=1 \) then \( \nu_T(S) \geq 0\) by Lemma~\ref{lemma: trees}.
	%	\item  If  \( S  \) is not connected, then \( \nu_T(S) = 0 \) by~Proposition~\ref{prop: connected}.
	%	\item If \( S  \)  is a connected set that does not contain the origin \( o, \) then \( S \) is an interval of length \( j \coloneqq |S|. \)  By Lemma~\ref{lemma: trees}, we have
	%	\begin{align*}
	%		&\nu_T(S) = \sum_{J \subseteq \mathcal{B}^-(S)} (-1)^{|J|} 
    %     \log P\pigl( X \bigl(J \cup (\mathcal{B}^+(S)\smallsetminus \partial_S J) \bigr) \equiv 0 \pigr)
     %    \\&\qquad=
    %     \log  \bigl( (1-p)^{j+1} r + (1-(1-p)^{j+1})r^2  \bigr) + 
    %     \log  \bigl( (1-p)^{j-1}r+(1-(1-p)^{j-1})r^2  \bigr)
    %     \\&\qquad\qquad-2\log  \bigl( (1-p)^{j}r +(1-(1-p)^{j})r^2 \bigr)
    %     \\&\qquad= 
    %     \log \Bigl( 1 +\frac{p^2(1 - p)^{j-1} r (1 - r) }{(r+(1 - p)^j  - (1 - p)^j r)^2 }\Bigr) > 0.
	%	\end{align*} 
	%\end{enumerate}
 	%
   
   Now assume that \( S \) is a connected set that contains the origin. Let \( V_1, \) \( V_2, \) and \( V_3, \) be the sets of vertices on each arm of \( T_3,\) so that \( V(T_3) = \{ o \} \sqcup V_1 \sqcup V_2 \sqcup V_3, \) and the restriction of \( T_3 \) to \( V_j \) is an infinite path for each \( j \in [3]. \) 
   Let \( T' \) be the restriction of \( T_3 \) to the set of vertices on distance at most two from the root, and let \( X' \) be the corresponding tree-indexed Markov chain with parameters \( (r,(p_e)), \) where for \( j \in [3], \) we let \( p_{e_{j,1}} = 1-(1-p)^{|S \cap V_J|} \) and \( p_{e_{j,2}}= p. \) Then, by Lemma~\ref{lemma: trees}, we have   
   \begin{equation}\label{equation: equality}
   		\nu_T(S) = \nu_{T'}(N(o))
   \end{equation}
   and hence, by Proposition~\ref{proposition: general small 3-tree}, we have \( \nu_T(S) \geq 0\). 
   
   Finally, note that by Lemma~\ref{equation: equality}, \( \nu_T(S) \) is independent of the choice of \( T \) whenever \( S \cup \mathcal{B}^+(S) \subseteq T. \) Using Lemma~\ref{lemma: finite to infinite}\ref{item: finite to infinite ii}, the desired conclusion immediately follows.  
\end{proof}

\begin{proof}[Proof of Proposition~\ref{proposition: general small 3-tree}]
	To simplify notation, for \( J \subseteq \mathcal{B}^-(S),\) consider the event
	\[
	\mathcal{E}_J \coloneqq \bigl\{ X \bigl(J \cup (\mathcal{B}^+(S)\smallsetminus \partial_S J) \bigr) \equiv 0\bigr\}.
	\]
	By Lemma~\ref{lemma: trees}, we  have
    \begin{equation*} 
        \begin{split} 
        &\nu (S) =  
        \sum_{J \subseteq \mathcal{B}^-(S)} (-1)^{|J|} 
         \log P( \mathcal{E}_J)
         = 
         \log \prod_{J \subseteq \mathcal{B}^-(S)} P( \mathcal{E}_J )^{ (-1)^{|J|} }
         \\&\qquad=\log \Biggl( 1 + \frac{\prod_{\substack{J \subseteq \mathcal{B}^-(S) {\colon} \\ |J| \text{ even}}} P(\mathcal{E}_J)-\prod_{\substack{J \subseteq \mathcal{B}^-(S) {\colon} \\ |J| \text{ odd}}} P(\mathcal{E}_J)}{\prod_{\substack{J \subseteq \mathcal{B}^-(S) {\colon} \\ |J| \text{ odd}}} P(\mathcal{E}_J)} \Biggr).
        \end{split}
    \end{equation*} 
    Hence \( \nu(S) \geq 0 \) if and only if    
    \begin{equation}\label{eq: numerator}
    	\prod_{\substack{J \subseteq \mathcal{B}^-(S) {\colon} \\ |J| \text{ even}}} P(\mathcal{E}_J)
    	- 
    	\prod_{\substack{J \subseteq \mathcal{B}^-(S) {\colon} \\ |J| \text{ odd}}} P(\mathcal{E}_J) \geq 0.
    \end{equation}
    To simplify notation, let
    \begin{equation*}
    	\begin{cases} 
    		q_{j,1}^0 \coloneqq P(X_{v_{j,1}}=0 \mid X_o = 0) = rp_{e_{j,1}} + (1-p_{e_{j,1}})\cr 
    		q_{j,1}^1 \coloneqq P(X_{v_{j,1}}=0 \mid X_o = 1) = r p_{e_{j,1}} \cr 
    		q_{j,2}^0 \coloneqq P(X(v_{j,2}) = 0\mid X_0=0) = r  p_{e_{j,2}}  + r p_{e_{j,1}}(1-p_{e_{j,2}})  + (1-p_{e_{j,1}})(1-p_{e_{j,2}}) \cr 
    		q_{j,2}^1 \coloneqq P(X(v_{j,2}) = 0\mid X_0=1) = r p_{e_{j,2}} + r  p_{e_{j,1}} (1-p_{e_{j,2}}).
    		%q_{j,2}^0 \coloneqq P(X(v_{j,2}) = 0\mid X_0=0) = r  p_{e_{j,1}}  + r (1-p_{e_{j,1}})p_{e_{j,2}}  + (1-p_{e_{j,1}})(1-p_{e_{j,2}}) \cr 
    		%q_{j,2}^1 \coloneqq P(X(v_{j,2}) = 0\mid X_0=1) = r p_{e_{j,1}} + r (1-p_{e_{j,1}}) p_{e_{j,2}}.
    	\end{cases}
    \end{equation*} 
    Then, for the left-hand side of~\eqref{eq: numerator}, we have     \begin{align*}
    	&
    	\prod_{\substack{J \subseteq \mathcal{B}^-(S) {\colon}  |J| \text{ even}}} P(\mathcal{E}_J)
    	- 
    	\prod_{\substack{J \subseteq \mathcal{B}^-(S) {\colon}  |J| \text{ odd}}} P(\mathcal{E}_J)  
    	\\&\qquad=
    	\Bigl[r  \prod_{j \in [3]} q_{j,2}^0
    	+
    	(1-r)  \prod_{j \in [3]} q_{j,2}^1 \Bigr]
    	\cdot \prod_{k \in [3]}  \Bigl[ r q_{k,2}^0 \prod_{j \in [3]\smallsetminus \{ k \}}  q_{j,1}^0 
    	+(1-r) q_{k,2}^1 \prod_{j \in [3]\smallsetminus \{ k \}} q_{j,1}^1 \Bigr] 
    	\\&\qquad\qquad -
    	(r \prod_{j \in [3]} q_{j,1}^0 +(1-r)  \prod_{j \in [3]} q_{j,1}^1 )
    	\cdot 
    	\prod_{k \in [3]} \Bigl[ r  q_{k,1}^0\prod_{j \in [3]\smallsetminus \{ k \}} q_{j,2}^0
    	+(1-r) q_{k,1}^1 
    	\prod_{j \in [3]\smallsetminus \{ k \}} q_{j,2}^1\Bigr] 
    	\\&\qquad=
    	r(1-r)\prod_{j \in [3]}( q_{j,1}^0 q_{j,2}^1-q_{j,1}^1 q_{j,2}^0) \Bigl(  r^2\prod_{j \in [3],k \in [2]} q_{j,k}^0 - (1 - r)^2\prod_{j \in [3],k \in [2]} q_{j,k}^1\Bigr).
    \end{align*} 
    Now note that for any \( j \in [3], \) we have \( q_{j,1}^0 q_{j,2}^1-q_{j,1}^1 q_{j,2}^0 = r(1-p_{e_{j,1}})p_{e_{j,2}} >0 \),  that for all \( j \in [3] \) and \( k \in [2],\) by positive association, we have \( q_{j,k}^0 \geq q_{j,k}^1 ,\) and that since \( r \geq 1/2, \) we have \( r^2 \geq (1-r)^2. \) Combining these observations with the above equations, the desired conclusion follows.
\end{proof}

\section{Proof of Theorem~\ref{theorem: p close to zero}}\label{section: general infinite trees}

In this section, we prove Theorem~\ref{theorem: p close to zero}. To this end, we first note that \ref{theorem: infinite trees i} and \ref{theorem: infinite trees iii} of Theorem~\ref{theorem: p close to zero} are essentially an immediate consequence of  Theorem~\ref{theorem: finite trees}\ref{theorem: finite trees iii} and Theorem~\ref{theorem: general octopus trees}\ref{item: general octopus trees a} respectively.

\begin{proof}[Proof of Theorem~\ref{theorem: p close to zero}\ref{theorem: infinite trees i}]
  	Using Lemma~\ref{lemma: finite to infinite}\ref{item: finite to infinite iii}, Theorem~\ref{theorem: p close to zero}\ref{theorem: infinite trees i} is an immediate consequence of Theorem~\ref{theorem: finite trees}\ref{theorem: finite trees iii}. This concludes the proof.\end{proof}

\begin{proof}[Proof of Theorem~\ref{theorem: p close to zero}\ref{theorem: infinite trees iii}]
	By assumption, \( T_{n_3} \) is a subtree of \( T . \) By Theorem~\ref{theorem: general octopus trees}\ref{item: general octopus trees a}, \( X(T_{n_3}) \notin \mathcal{R}, \) and hence, by Lemma~\ref{lemma: finite to infinite}\ref{item: finite to infinite iii}, \( X \notin \mathcal{R}. \) This concludes the proof.
\end{proof}

We now give a proof of Theorem~\ref{theorem: p close to zero}\ref{theorem: infinite trees ii}. The main idea of the proof is that for any \( m \geq 3, \) if \( p \) is very small and  \( T \) is a tree which contains \( m \) eventually disjoint infinite paths starting from \( o, \) then \( T \) contains a subtree \( X' \) such that \( X(T') \) with parameters \( (p,r) \) is very close to \( X(T_m) \) with parameters \( (p',r) \), except possible very close to the root vertex \( o \). Since these edges are only finitely many, by continuity, the result essentially follows from Theorem~\ref{theorem: general octopus trees}\ref{item: general octopus trees a}.

\begin{proof}[Proof of Theorem~\ref{theorem: p close to zero}\ref{theorem: infinite trees ii}] 

	 Let \( m \geq 1 \) be such that there are  \( m \) distinct infinite paths \( P_1,\dots, P_m\) in \( T \) starting from \( o ,\) and such that \( -(1-r)/r \in (r_2^{(m)},r_1^{(m)}). \) Since \(r<r_1(n_2)\) implies \(  -(1-r)/r < r_1^{(n_2)}, \) such \( m \) exists, and moreover, for such \(m, \) we have \( \Li_{1-m}(-(1-r)/r) >0. \)

	For \( j \geq 1, \) let \( T^{ j} \) be the subgraph of \( T \) which corresponds to all vertices at distance at most \( j \) from the origin. Let \(k = k(m) \) be the smallest \( k \geq 0 \) such that the sets   \( P_1 \smallsetminus V(T^{ k}) \), \dots , and \( P_m \smallsetminus V(T^{ k}) \) are disjoint. Note that  \( E(T^{ k}) \) is finite.
	Let \( V_k \coloneqq (P_1 \cup \dots \cup P_m ) \cap V(T^{ k}). \)
	For \( j\geq k, \) let \( V_{k,j} \) be the union of \( V_k \) and the \( (k+j)\)th and \((k+2j)\)th vertices on each of the paths \( P_1,\dots, P_m\) (see Figure~\ref{figure: last fig 3}).

	\begin{figure}[h]
		\begin{subfigure}[t]{.23\textwidth}\centering
			\begin{tikzpicture}[scale=.3] 
    	
    	% level 6
    	\foreach \i in {0,...,3}
    		\foreach \j in {-1,1}
    			\foreach \k in {-1,1}
    				\foreach \l in {-1,1}
    					\foreach \m in {-1,1}
    						\foreach \n in {-1,1}
    						{
    							\draw ({5*cos(120*\i+30*\j+15*\k+8*\l+4*\m)},{5*sin(120*\i+30*\j+15*\k +8*\l+4*\m)})
    							--
    							({5.8*cos(120*\i+30*\j+15*\k+8*\l+4*\m+2*\n)},{5.8*sin(120*\i+30*\j+15*\k +8*\l+4*\m+2*\n)});  
    						} 
    	
    	% level 5
    	\foreach \i in {0,...,3}
    		\foreach \j in {-1,1}
    			\foreach \k in {-1,1}
    				\foreach \l in {-1,1}
    					\foreach \m in {-1,1}
    					{
    						\draw ({4*cos(120*\i+30*\j+15*\k+8*\l)},{4*sin(120*\i+30*\j+15*\k +8*\l)})
    						--
    						({5*cos(120*\i+30*\j+15*\k+8*\l+4*\m)},{5*sin(120*\i+30*\j+15*\k +8*\l+4*\m)});  
    						
    						\filldraw[fill=white] ({5*cos(120*\i+30*\j+15*\k+8*\l+4*\m)},{5*sin(120*\i+30*\j+15*\k +8*\l+4*\m)}) circle (4pt); 
    					}  
    	
    	% level 4
    	\foreach \i in {0,...,3}
    		\foreach \j in {-1,1}
    			\foreach \k in {-1,1}
    				\foreach \l in {-1,1}
    				{
    					\draw ({3*cos(120*\i+30*\j+15*\k)},{3*sin(120*\i+30*\j+15*\k )}) -- ({4*cos(120*\i+30*\j+15*\k+8*\l)},{4*sin(120*\i+30*\j+15*\k +8*\l)});  
    					\filldraw[fill=white] ({4*cos(120*\i+30*\j+15*\k+8*\l)},{4*sin(120*\i+30*\j+15*\k +8*\l)}) circle (5pt);   
    				}

    	% level 3
    	\foreach \i in {0,...,3}
    		\foreach \j in {-1,1}
    			\foreach \k in {-1,1}
    			{
    				\draw ({2*cos(120*\i+30*\j)},{2*sin(120*\i+30*\j)}) -- ({3*cos(120*\i+30*\j+15*\k)},{3*sin(120*\i+30*\j+15*\k )}); 
    				\filldraw[fill=white] ({3*cos(120*\i+30*\j+15*\k)},{3*sin(120*\i+30*\j+15*\k)}) circle (6pt);   
    			}   
    	
    	% level 2
    	\foreach \i in {0,...,3} 
    		\foreach \j in {-1,1}
    		{
    			\draw ({cos(120*\i)},{sin(120*\i)}) -- ({2*cos(120*\i+30*\j)},{2*sin(120*\i+30*\j)}); 
    			\filldraw[fill=white] ({2*cos(120*\i+30*\j)},{2*sin(120*\i+30*\j)}) circle (7pt);   
    		}   
		
    	% level 1
    	\foreach \i in {0,...,3}
    	{	  
    		\draw (0,0) -- ({cos(120*\i)},{sin(120*\i)});
    		\filldraw[fill=white] ({cos(120*\i)},{sin(120*\i)}) circle (8pt);  
    	}
    	
    	% level 0 
    	\filldraw[fill=white] (0,0) circle (10pt);   
	
     \end{tikzpicture}
		\caption{\( T \)}
		\end{subfigure}
		\hfil
		\begin{subfigure}[t]{.23\textwidth}\centering
			\begin{tikzpicture}[scale=.3] 
    	
    	% level 6
    	\foreach \i in {0,...,3}
    		\foreach \j in {-1,1}
    			\foreach \k in {-1,1}
    				\foreach \l in {-1,1}
    					\foreach \m in {-1,1}
    						\foreach \n in {-1,1}
    						{
    							\draw[help lines] ({5*cos(120*\i+30*\j+15*\k+8*\l+4*\m)},{5*sin(120*\i+30*\j+15*\k +8*\l+4*\m)})
    							--
    							({5.8*cos(120*\i+30*\j+15*\k+8*\l+4*\m+2*\n)},{5.8*sin(120*\i+30*\j+15*\k +8*\l+4*\m+2*\n)});  
    						} 
    	
    	{ %P1
    		\def\i{0}\def\j{1}\def\k{-1}\def\l{1}\def\m{1}\def\n{1}
    		\draw ({5*cos(120*\i+30*\j+15*\k+8*\l+4*\m)},{5*sin(120*\i+30*\j+15*\k +8*\l+4*\m)})
    							--
    							({5.8*cos(120*\i+30*\j+15*\k+8*\l+4*\m+2*\n)},{5.8*sin(120*\i+30*\j+15*\k +8*\l+4*\m+2*\n)});   
    	} 
    	
    	{ %P2
    		\def\i{1}\def\j{1}\def\k{1}\def\l{1}\def\m{1}\def\n{1}
    		\draw ({5*cos(120*\i+30*\j+15*\k+8*\l+4*\m)},{5*sin(120*\i+30*\j+15*\k +8*\l+4*\m)})
    							--
    							({5.8*cos(120*\i+30*\j+15*\k+8*\l+4*\m+2*\n)},{5.8*sin(120*\i+30*\j+15*\k +8*\l+4*\m+2*\n)});  
    							    	} 
    	
    	{ %P3
    		\def\i{2}\def\j{1}\def\k{1}\def\l{-1}\def\m{-1}\def\n{1}
    		\draw ({5*cos(120*\i+30*\j+15*\k+8*\l+4*\m)},{5*sin(120*\i+30*\j+15*\k +8*\l+4*\m)})
    							--
    							({5.8*cos(120*\i+30*\j+15*\k+8*\l+4*\m+2*\n)},{5.8*sin(120*\i+30*\j+15*\k +8*\l+4*\m+2*\n)});
    	} 
    	
    	{ %P4
    		\def\i{2}\def\j{1}\def\k{-1}\def\l{1}\def\m{1}\def\n{1}
    		\draw ({5*cos(120*\i+30*\j+15*\k+8*\l+4*\m)},{5*sin(120*\i+30*\j+15*\k +8*\l+4*\m)})
    							--
    							({5.8*cos(120*\i+30*\j+15*\k+8*\l+4*\m+2*\n)},{5.8*sin(120*\i+30*\j+15*\k +8*\l+4*\m+2*\n)});  
    	} 
    	
    	{ %P5
    		\def\i{2}\def\j{-1}\def\k{-1}\def\l{1}\def\m{-1}\def\n{1}
    		\draw ({5*cos(120*\i+30*\j+15*\k+8*\l+4*\m)},{5*sin(120*\i+30*\j+15*\k +8*\l+4*\m)})
    							--
    							({5.8*cos(120*\i+30*\j+15*\k+8*\l+4*\m+2*\n)},{5.8*sin(120*\i+30*\j+15*\k +8*\l+4*\m+2*\n)});  
    	} 
    	
    	{ %P6
    		\def\i{0}\def\j{-1}\def\k{-1}\def\l{-1}\def\m{1}\def\n{1}
    		\draw ({5*cos(120*\i+30*\j+15*\k+8*\l+4*\m)},{5*sin(120*\i+30*\j+15*\k +8*\l+4*\m)})
    							--
    							({5.8*cos(120*\i+30*\j+15*\k+8*\l+4*\m+2*\n)},{5.8*sin(120*\i+30*\j+15*\k +8*\l+4*\m+2*\n)});  
    	} 
    	
    	{ %P7
    		\def\i{0}\def\j{-1}\def\k{1}\def\l{-1}\def\m{1}\def\n{1}
    		\draw ({5*cos(120*\i+30*\j+15*\k+8*\l+4*\m)},{5*sin(120*\i+30*\j+15*\k +8*\l+4*\m)})
    							--
    							({5.8*cos(120*\i+30*\j+15*\k+8*\l+4*\m+2*\n)},{5.8*sin(120*\i+30*\j+15*\k +8*\l+4*\m+2*\n)});  
    	} 
    				
    	% ------------------------	

    	% level 5
    	\foreach \i in {0,...,3}
    		\foreach \j in {-1,1}
    			\foreach \k in {-1,1}
    				\foreach \l in {-1,1}
    					\foreach \m in {-1,1}
    					{
    						\draw[help lines] ({4*cos(120*\i+30*\j+15*\k+8*\l)},{4*sin(120*\i+30*\j+15*\k +8*\l)})
    						--
    						({5*cos(120*\i+30*\j+15*\k+8*\l+4*\m)},{5*sin(120*\i+30*\j+15*\k +8*\l+4*\m)});  
    						
    						\filldraw[help lines,fill=white] ({5*cos(120*\i+30*\j+15*\k+8*\l+4*\m)},{5*sin(120*\i+30*\j+15*\k +8*\l+4*\m)}) circle (4pt); 
    					}    
    	
    	{ %P1
    		\def\i{0}\def\j{1}\def\k{-1}\def\l{1}\def\m{1}\def\n{1}
    		\draw ({4*cos(120*\i+30*\j+15*\k+8*\l)},{4*sin(120*\i+30*\j+15*\k +8*\l)})
    						--
    						({5*cos(120*\i+30*\j+15*\k+8*\l+4*\m)},{5*sin(120*\i+30*\j+15*\k +8*\l+4*\m)});  
    						
    		\filldraw[fill=black] ({5*cos(120*\i+30*\j+15*\k+8*\l+4*\m)},{5*sin(120*\i+30*\j+15*\k +8*\l+4*\m)}) circle (4pt); 
    	} 
    	
    	{ %P2
    		\def\i{1}\def\j{1}\def\k{1}\def\l{1}\def\m{1}\def\n{1}
    		\draw ({4*cos(120*\i+30*\j+15*\k+8*\l)},{4*sin(120*\i+30*\j+15*\k +8*\l)})
    						--
    						({5*cos(120*\i+30*\j+15*\k+8*\l+4*\m)},{5*sin(120*\i+30*\j+15*\k +8*\l+4*\m)});  
    						
    		\filldraw[fill=black] ({5*cos(120*\i+30*\j+15*\k+8*\l+4*\m)},{5*sin(120*\i+30*\j+15*\k +8*\l+4*\m)}) circle (4pt);
    	} 
    	
    	{ %P3
    		\def\i{2}\def\j{1}\def\k{1}\def\l{-1}\def\m{-1}\def\n{1}
    		\draw ({4*cos(120*\i+30*\j+15*\k+8*\l)},{4*sin(120*\i+30*\j+15*\k +8*\l)})
    						--
    						({5*cos(120*\i+30*\j+15*\k+8*\l+4*\m)},{5*sin(120*\i+30*\j+15*\k +8*\l+4*\m)});  
    						
    		\filldraw[fill=black] ({5*cos(120*\i+30*\j+15*\k+8*\l+4*\m)},{5*sin(120*\i+30*\j+15*\k +8*\l+4*\m)}) circle (4pt);
    	} 
    	
    	{ %P4
    		\def\i{2}\def\j{1}\def\k{-1}\def\l{1}\def\m{1}\def\n{1}
    		\draw ({4*cos(120*\i+30*\j+15*\k+8*\l)},{4*sin(120*\i+30*\j+15*\k +8*\l)})
    						--
    						({5*cos(120*\i+30*\j+15*\k+8*\l+4*\m)},{5*sin(120*\i+30*\j+15*\k +8*\l+4*\m)});  
    						
    		\filldraw[fill=black] ({5*cos(120*\i+30*\j+15*\k+8*\l+4*\m)},{5*sin(120*\i+30*\j+15*\k +8*\l+4*\m)}) circle (4pt);
    	} 
    	
    	{ %P5
    		\def\i{2}\def\j{-1}\def\k{-1}\def\l{1}\def\m{-1}\def\n{1}
    		\draw ({4*cos(120*\i+30*\j+15*\k+8*\l)},{4*sin(120*\i+30*\j+15*\k +8*\l)})
    						--
    						({5*cos(120*\i+30*\j+15*\k+8*\l+4*\m)},{5*sin(120*\i+30*\j+15*\k +8*\l+4*\m)});  
    						
    		\filldraw[fill=black] ({5*cos(120*\i+30*\j+15*\k+8*\l+4*\m)},{5*sin(120*\i+30*\j+15*\k +8*\l+4*\m)}) circle (4pt);
    	} 
    	
    	{ %P6
    		\def\i{0}\def\j{-1}\def\k{-1}\def\l{-1}\def\m{1}\def\n{1}
    		\draw ({4*cos(120*\i+30*\j+15*\k+8*\l)},{4*sin(120*\i+30*\j+15*\k +8*\l)})
    						--
    						({5*cos(120*\i+30*\j+15*\k+8*\l+4*\m)},{5*sin(120*\i+30*\j+15*\k +8*\l+4*\m)});  
    						
    		\filldraw[fill=black] ({5*cos(120*\i+30*\j+15*\k+8*\l+4*\m)},{5*sin(120*\i+30*\j+15*\k +8*\l+4*\m)}) circle (4pt);
    	} 
    	
    	{ %P7
    		\def\i{0}\def\j{-1}\def\k{1}\def\l{-1}\def\m{1}\def\n{1}
    		\draw ({4*cos(120*\i+30*\j+15*\k+8*\l)},{4*sin(120*\i+30*\j+15*\k +8*\l)})
    						--
    						({5*cos(120*\i+30*\j+15*\k+8*\l+4*\m)},{5*sin(120*\i+30*\j+15*\k +8*\l+4*\m)});  
    						
    		\filldraw[fill=black] ({5*cos(120*\i+30*\j+15*\k+8*\l+4*\m)},{5*sin(120*\i+30*\j+15*\k +8*\l+4*\m)}) circle (4pt);
    	} 
    				
    	% ------------------------	
    	
    	% level 4
    	\foreach \i in {0,...,3}
    		\foreach \j in {-1,1}
    			\foreach \k in {-1,1}
    				\foreach \l in {-1,1}
    				{
    					\draw[help lines] ({3*cos(120*\i+30*\j+15*\k)},{3*sin(120*\i+30*\j+15*\k )}) -- ({4*cos(120*\i+30*\j+15*\k+8*\l)},{4*sin(120*\i+30*\j+15*\k +8*\l)});  
    					\filldraw[help lines,fill=white] ({4*cos(120*\i+30*\j+15*\k+8*\l)},{4*sin(120*\i+30*\j+15*\k +8*\l)}) circle (5pt);   
    				}   
    	
    	{ %P1
    		\def\i{0}\def\j{1}\def\k{-1}\def\l{1}\def\m{1}\def\n{1}
    		\draw[] ({3*cos(120*\i+30*\j+15*\k)},{3*sin(120*\i+30*\j+15*\k )}) -- ({4*cos(120*\i+30*\j+15*\k+8*\l)},{4*sin(120*\i+30*\j+15*\k +8*\l)});  
    		\filldraw[fill=white] ({4*cos(120*\i+30*\j+15*\k+8*\l)},{4*sin(120*\i+30*\j+15*\k +8*\l)}) circle (5pt); 
    	} 
    	
    	{ %P2
    		\def\i{1}\def\j{1}\def\k{1}\def\l{1}\def\m{1}\def\n{1}
    		\draw[help lines] ({3*cos(120*\i+30*\j+15*\k)},{3*sin(120*\i+30*\j+15*\k )}) -- ({4*cos(120*\i+30*\j+15*\k+8*\l)},{4*sin(120*\i+30*\j+15*\k +8*\l)});  
    		\filldraw[fill=white] ({4*cos(120*\i+30*\j+15*\k+8*\l)},{4*sin(120*\i+30*\j+15*\k +8*\l)}) circle (5pt); 
    	} 
    	
    	{ %P3
    		\def\i{2}\def\j{1}\def\k{1}\def\l{-1}\def\m{1}\def\n{1}
    		\draw[] ({3*cos(120*\i+30*\j+15*\k)},{3*sin(120*\i+30*\j+15*\k )}) -- ({4*cos(120*\i+30*\j+15*\k+8*\l)},{4*sin(120*\i+30*\j+15*\k +8*\l)});  
    		\filldraw[fill=white] ({4*cos(120*\i+30*\j+15*\k+8*\l)},{4*sin(120*\i+30*\j+15*\k +8*\l)}) circle (5pt); 
    	} 
    	
    	{ %P4
    		\def\i{2}\def\j{1}\def\k{-1}\def\l{1}\def\m{1}\def\n{1}
    		\draw[] ({3*cos(120*\i+30*\j+15*\k)},{3*sin(120*\i+30*\j+15*\k )}) -- ({4*cos(120*\i+30*\j+15*\k+8*\l)},{4*sin(120*\i+30*\j+15*\k +8*\l)});  
    		\filldraw[fill=white] ({4*cos(120*\i+30*\j+15*\k+8*\l)},{4*sin(120*\i+30*\j+15*\k +8*\l)}) circle (5pt); 
    	} 
    	
    	{ %P5
    		\def\i{2}\def\j{-1}\def\k{-1}\def\l{1}\def\m{1}\def\n{1}
    		\draw[] ({3*cos(120*\i+30*\j+15*\k)},{3*sin(120*\i+30*\j+15*\k )}) -- ({4*cos(120*\i+30*\j+15*\k+8*\l)},{4*sin(120*\i+30*\j+15*\k +8*\l)});  
    		\filldraw[fill=white] ({4*cos(120*\i+30*\j+15*\k+8*\l)},{4*sin(120*\i+30*\j+15*\k +8*\l)}) circle (5pt); 
    	} 
    	
    	{ %P6
    		\def\i{0}\def\j{-1}\def\k{-1}\def\l{-1}\def\m{1}\def\n{1}
    		\draw[] ({3*cos(120*\i+30*\j+15*\k)},{3*sin(120*\i+30*\j+15*\k )}) -- ({4*cos(120*\i+30*\j+15*\k+8*\l)},{4*sin(120*\i+30*\j+15*\k +8*\l)});  
    		\filldraw[fill=white] ({4*cos(120*\i+30*\j+15*\k+8*\l)},{4*sin(120*\i+30*\j+15*\k +8*\l)}) circle (5pt); 
    	} 
    	
    	{ %P7
    		\def\i{0}\def\j{-1}\def\k{1}\def\l{-1}\def\m{1}\def\n{1}
    		\draw[] ({3*cos(120*\i+30*\j+15*\k)},{3*sin(120*\i+30*\j+15*\k )}) -- ({4*cos(120*\i+30*\j+15*\k+8*\l)},{4*sin(120*\i+30*\j+15*\k +8*\l)});  
    		\filldraw[fill=white] ({4*cos(120*\i+30*\j+15*\k+8*\l)},{4*sin(120*\i+30*\j+15*\k +8*\l)}) circle (5pt); 
    	} 
    				
    	% ------------------------			
		
    	% level 3
    	\foreach \i in {0,...,3}
    		\foreach \j in {-1,1}
    			\foreach \k in {-1,1}
    			{
    				\draw[help lines] ({2*cos(120*\i+30*\j)},{2*sin(120*\i+30*\j)}) -- ({3*cos(120*\i+30*\j+15*\k)},{3*sin(120*\i+30*\j+15*\k )}); 
    				\filldraw[help lines,fill=white] ({3*cos(120*\i+30*\j+15*\k)},{3*sin(120*\i+30*\j+15*\k)}) circle (6pt);   
    			}   
    	
    	{ %P1
    		\def\i{0}\def\j{1}\def\k{-1}\def\l{1}\def\m{1}\def\n{1}
    		\draw ({2*cos(120*\i+30*\j)},{2*sin(120*\i+30*\j)}) 
    			-- ({3*cos(120*\i+30*\j+15*\k)},{3*sin(120*\i+30*\j+15*\k )}); 
    		\filldraw[fill=white] ({3*cos(120*\i+30*\j+15*\k)},{3*sin(120*\i+30*\j+15*\k)}) circle (6pt); 
    	} 
    	
    	{ %P2
    		\def\i{1}\def\j{1}\def\k{1}\def\l{1}\def\m{1}\def\n{1}
    		\draw ({2*cos(120*\i+30*\j)},{2*sin(120*\i+30*\j)}) 
    			-- ({3*cos(120*\i+30*\j+15*\k)},{3*sin(120*\i+30*\j+15*\k )}); 
    		\filldraw[fill=white] ({3*cos(120*\i+30*\j+15*\k)},{3*sin(120*\i+30*\j+15*\k)}) circle (6pt); 
    	} 
    	
    	{ %P3
    		\def\i{2}\def\j{1}\def\k{1}\def\l{1}\def\m{1}\def\n{1}
    		\draw ({2*cos(120*\i+30*\j)},{2*sin(120*\i+30*\j)}) 
    			-- ({3*cos(120*\i+30*\j+15*\k)},{3*sin(120*\i+30*\j+15*\k )}); 
    		\filldraw[fill=white] ({3*cos(120*\i+30*\j+15*\k)},{3*sin(120*\i+30*\j+15*\k)}) circle (6pt); 
    	} 
    	
    	{ %P4
    		\def\i{2}\def\j{1}\def\k{-1}\def\l{1}\def\m{1}\def\n{1}
    		\draw ({2*cos(120*\i+30*\j)},{2*sin(120*\i+30*\j)}) 
    			-- ({3*cos(120*\i+30*\j+15*\k)},{3*sin(120*\i+30*\j+15*\k )}); 
    		\filldraw[fill=white] ({3*cos(120*\i+30*\j+15*\k)},{3*sin(120*\i+30*\j+15*\k)}) circle (6pt); 
    	} 
    	
    	{ %P5
    		\def\i{2}\def\j{-1}\def\k{-1}\def\l{1}\def\m{1}\def\n{1}
    		\draw ({2*cos(120*\i+30*\j)},{2*sin(120*\i+30*\j)}) 
    			-- ({3*cos(120*\i+30*\j+15*\k)},{3*sin(120*\i+30*\j+15*\k )}); 
    		\filldraw[fill=white] ({3*cos(120*\i+30*\j+15*\k)},{3*sin(120*\i+30*\j+15*\k)}) circle (6pt); 
    	} 
    	
    	{ %P6
    		\def\i{0}\def\j{-1}\def\k{-1}\def\l{1}\def\m{1}\def\n{1}
    		\draw ({2*cos(120*\i+30*\j)},{2*sin(120*\i+30*\j)}) 
    			-- ({3*cos(120*\i+30*\j+15*\k)},{3*sin(120*\i+30*\j+15*\k )}); 
    		\filldraw[fill=white] ({3*cos(120*\i+30*\j+15*\k)},{3*sin(120*\i+30*\j+15*\k)}) circle (6pt); 
    	} 
    	
    	{ %P7
    		\def\i{0}\def\j{-1}\def\k{1}\def\l{1}\def\m{1}\def\n{1}
    		\draw ({2*cos(120*\i+30*\j)},{2*sin(120*\i+30*\j)}) 
    			-- ({3*cos(120*\i+30*\j+15*\k)},{3*sin(120*\i+30*\j+15*\k )}); 
    		\filldraw[fill=white] ({3*cos(120*\i+30*\j+15*\k)},{3*sin(120*\i+30*\j+15*\k)}) circle (6pt); 
    	} 
    	
    	% level 2
    	\foreach \i in {0,...,3} 
    		\foreach \j in {-1,1}
    		{
    			\draw ({cos(120*\i)},{sin(120*\i)}) -- ({2*cos(120*\i+30*\j)},{2*sin(120*\i+30*\j)}); 
    			\filldraw[fill=red] ({2*cos(120*\i+30*\j)},{2*sin(120*\i+30*\j)}) circle (7pt);   
    		}   
		
    	% level 1
    	\foreach \i in {0,...,3}
    	{	  
    		\draw (0,0) -- ({cos(120*\i)},{sin(120*\i)});
    		\filldraw[fill=red] ({cos(120*\i)},{sin(120*\i)}) circle (8pt);  
    	}
    	
    	% level 0 
    	\filldraw[fill=red] (0,0) circle (10pt);   
	
     \end{tikzpicture}
		\caption{\( T^2 \) in red, paths \( P_1, \dots, P_7,\) and on each path, a vertex at distance 3 from \( T^2 \) in black. }
		\label{fig: last proof b}
		\end{subfigure}
		\hfil
		\begin{subfigure}[t]{.23\textwidth}\centering
			\begin{tikzpicture}[scale=.3]

    	% level 3 
    	
    	{ %P1
    		\def\i{0}\def\j{1}\def\k{-1}\def\l{1}\def\m{1}\def\n{1}
    		\draw ({2*cos(120*\i+30*\j)},{2*sin(120*\i+30*\j)}) 
    			-- ({3*cos(120*\i+30*\j+15*\k)},{3*sin(120*\i+30*\j+15*\k )})
    			-- ({4*cos(120*\i+30*\j+15*\k)},{4*sin(120*\i+30*\j+15*\k )}); 
    			 
    		\filldraw[fill=black] ({3*cos(120*\i+30*\j+15*\k)},{3*sin(120*\i+30*\j+15*\k)}) circle (6pt);  
    		
    		\filldraw[fill=white] ({4*cos(120*\i+30*\j+15*\k)},{4*sin(120*\i+30*\j+15*\k)}) circle (6pt); 
    	} 
    	
    	{ %P2
    		\def\i{1}\def\j{1}\def\k{1}\def\l{1}\def\m{1}\def\n{1}
    		
    		\draw ({2*cos(120*\i+30*\j)},{2*sin(120*\i+30*\j)}) 
    			-- ({3*cos(120*\i+30*\j+15*\k)},{3*sin(120*\i+30*\j+15*\k )})
    			-- ({4*cos(120*\i+30*\j+15*\k)},{4*sin(120*\i+30*\j+15*\k )}); 
    			
    		\filldraw[fill=black] ({3*cos(120*\i+30*\j+15*\k)},
    		{3*sin(120*\i+30*\j+15*\k)}) circle (6pt); 
    		
    		\filldraw[fill=white] ({4*cos(120*\i+30*\j+15*\k)},
    		{4*sin(120*\i+30*\j+15*\k)}) circle (6pt); 
    	} 
    	
    	{ %P3
    		\def\i{2}\def\j{1}\def\k{1}\def\l{1}\def\m{1}\def\n{1}
    		
    		\draw ({2*cos(120*\i+30*\j)},{2*sin(120*\i+30*\j)}) 
    			-- ({3*cos(120*\i+30*\j+15*\k)},{3*sin(120*\i+30*\j+15*\k )})
    			-- ({4*cos(120*\i+30*\j+15*\k)},{4*sin(120*\i+30*\j+15*\k )}); 
    			
    		\filldraw[fill=black] ({3*cos(120*\i+30*\j+15*\k)},{3*sin(120*\i+30*\j+15*\k)}) circle (6pt); 
    		\filldraw[fill=white] ({4*cos(120*\i+30*\j+15*\k)},{4*sin(120*\i+30*\j+15*\k)}) circle (6pt); 
    	} 
    	
    	{ %P4
    		\def\i{2}\def\j{1}\def\k{-1}\def\l{1}\def\m{1}\def\n{1}
    		\draw ({2*cos(120*\i+30*\j)},{2*sin(120*\i+30*\j)}) 
    			-- ({3*cos(120*\i+30*\j+15*\k)},{3*sin(120*\i+30*\j+15*\k )})
    			-- ({4*cos(120*\i+30*\j+15*\k)},{4*sin(120*\i+30*\j+15*\k )}); 
    		\filldraw[fill=black] ({3*cos(120*\i+30*\j+15*\k)},{3*sin(120*\i+30*\j+15*\k)}) circle (6pt); 
    		\filldraw[fill=white] ({4*cos(120*\i+30*\j+15*\k)},{4*sin(120*\i+30*\j+15*\k)}) circle (6pt); 
    	} 
    	
    	{ %P5
    		\def\i{2}\def\j{-1}\def\k{-1}\def\l{1}\def\m{1}\def\n{1}
    		\draw ({2*cos(120*\i+30*\j)},{2*sin(120*\i+30*\j)}) 
    			-- ({3*cos(120*\i+30*\j+15*\k)},{3*sin(120*\i+30*\j+15*\k )})
    			-- ({4*cos(120*\i+30*\j+15*\k)},{4*sin(120*\i+30*\j+15*\k )}); 
    		\filldraw[fill=black] ({3*cos(120*\i+30*\j+15*\k)},{3*sin(120*\i+30*\j+15*\k)}) circle (6pt); 
    		\filldraw[fill=white] ({4*cos(120*\i+30*\j+15*\k)},{4*sin(120*\i+30*\j+15*\k)}) circle (6pt); 
    	} 
    	
    	{ %P6
    		\def\i{0}\def\j{-1}\def\k{-1}\def\l{1}\def\m{1}\def\n{1}
    		\draw ({2*cos(120*\i+30*\j)},{2*sin(120*\i+30*\j)}) 
    			-- ({3*cos(120*\i+30*\j+15*\k)},{3*sin(120*\i+30*\j+15*\k )})
    			-- ({4*cos(120*\i+30*\j+15*\k)},{4*sin(120*\i+30*\j+15*\k )}); 
    		\filldraw[fill=black] ({3*cos(120*\i+30*\j+15*\k)},{3*sin(120*\i+30*\j+15*\k)}) circle (6pt); 
    		\filldraw[fill=white] ({4*cos(120*\i+30*\j+15*\k)},{4*sin(120*\i+30*\j+15*\k)}) circle (6pt); 
    	} 
    	
    	{ %P7
    		\def\i{0}\def\j{-1}\def\k{1}\def\l{1}\def\m{1}\def\n{1}
    		\draw ({2*cos(120*\i+30*\j)},{2*sin(120*\i+30*\j)}) 
    			-- ({3*cos(120*\i+30*\j+15*\k)},{3*sin(120*\i+30*\j+15*\k )})
    			-- ({4*cos(120*\i+30*\j+15*\k)},{4*sin(120*\i+30*\j+15*\k )}); 
    		\filldraw[fill=black] ({3*cos(120*\i+30*\j+15*\k)},{3*sin(120*\i+30*\j+15*\k)}) circle (6pt); 
    		\filldraw[fill=white] ({4*cos(120*\i+30*\j+15*\k)},{4*sin(120*\i+30*\j+15*\k)}) circle (6pt); 
    	} 
    	
    	% level 2
    	\foreach \i in {0,2} 
    		\foreach \j in {-1,1}
    		{
    			\draw ({cos(120*\i)},{sin(120*\i)}) -- ({2*cos(120*\i+30*\j)},{2*sin(120*\i+30*\j)}); 
    			\filldraw[fill=red] ({2*cos(120*\i+30*\j)},{2*sin(120*\i+30*\j)}) circle (7pt);   
    		}   
    	
    	% level 2
    	\foreach \i in {1} 
    		\foreach \j in {1}
    		{
    			\draw ({cos(120*\i)},{sin(120*\i)}) -- ({2*cos(120*\i+30*\j)},{2*sin(120*\i+30*\j)}); 
    			\filldraw[fill=red] ({2*cos(120*\i+30*\j)},{2*sin(120*\i+30*\j)}) circle (7pt);   
    		}   
		
    	% level 1
    	\foreach \i in {0,...,3}
    	{	  
    		\draw (0,0) -- ({cos(120*\i)},{sin(120*\i)});
    		\filldraw[fill=red] ({cos(120*\i)},{sin(120*\i)}) circle (8pt);  
    	}
    	
    	% level 0 
    	\filldraw[fill=red] (0,0) circle (10pt);   
	
     \end{tikzpicture}
		\caption{The subgraph of \( T \) induced by \( V_{2,1} \), and the set \( S_1 \) in black and red.}
		\label{figure: last fig 3}
		\end{subfigure}
		\hfil
	\begin{subfigure}[t]{.23\textwidth}\centering
			\begin{tikzpicture}[scale=.3]

    	% level 1
    	\foreach \i in {0,...,6}
    	{	  
    		\draw (0,0) -- ({2.4*cos(51.43*\i)},{2.4*sin(51.43*\i)});
    		\filldraw[fill=black] ({1.2*cos(51.43*\i)},{1.2*sin(51.43*\i)}) circle (6pt);  
    		\filldraw[fill=white] ({2.4*cos(51.43*\i)},{2.4*sin(51.43*\i)}) circle (6pt);  
    	}
    	
    	% level 0 
    	\filldraw[fill=red] (0,0) circle (10pt);   
	
     \end{tikzpicture}
		\caption{\( T_7 \) and the set \( S_0 \) in black and red.}\label{figure: new examples of subtrees d}\label{fig: d}
		\end{subfigure}
		\caption{In the figures above, for a tree \( T, \) we draw the related trees and paths that appear in the proof of Theorem~\ref{theorem: p close to zero}. Note that in this case we have \( k(7) = 2. \)}\label{figure: new examples of subtrees}
	\end{figure}
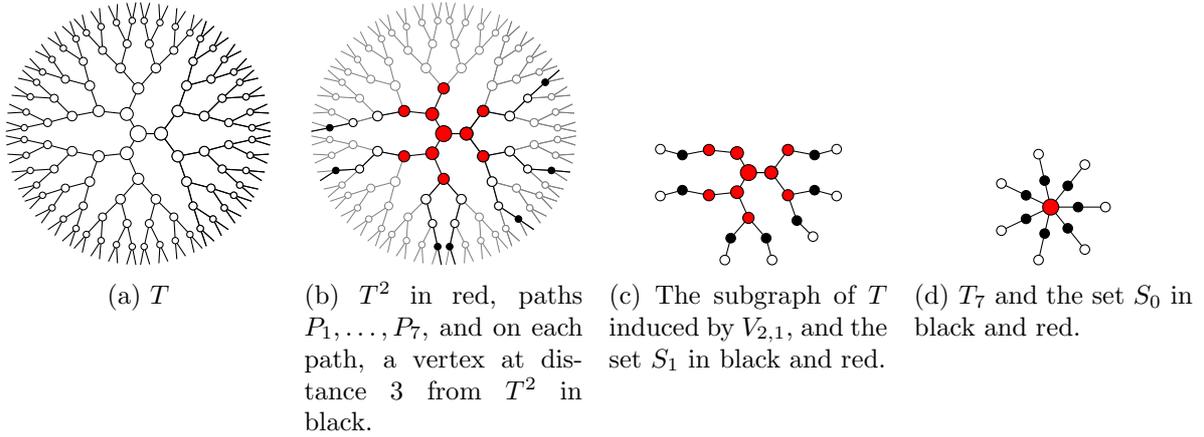

	Let \( T_m^{\leq 2} \) denote the subgraph of \( T_m \) induced by all vertices at distance at most two from the origin.
	Let \( p^0 \in (0,1), \) and let \( X^0 \) be the tree-indexed Markov chain on \( T_m^{\leq 2} \) with parameters \( (r,p^0) \) (see~Figure~\ref{fig: d}), and let \( \nu_0 \) be the corresponding signed measure. 
	Further, let \( S_0 \coloneqq \{ o \} \cup N(\{ o \}). \) 
	Since \( \Li_{1-m}(-(1-r)/r) > 0\), using a Taylor expansion, it follows from Proposition~\ref{proposition: new negative lemma for p close to one} implies that that \( \nu_0(S_0)<0 \) whenever \( p^0 \) is sufficiently close to one.
	By continuity, for all sufficiently small \( \varepsilon > 0 \) there is a non-empty interval \( I_0 \subseteq (0,1) \) such that \( {\nu_{0}(S_0) < -\varepsilon} \) whenever \( p^0 \in I_0. \)

	Next, let \( T^{k,1} \) be the subgraph of \( T \) induced by \( V_{k,1}. \)
	 Let \( p^1 \in [0,1), \) and let \( X^1 \) be the tree-indexed Markov chain on \( T^{k,1} \) with parameters \( (r,(p_e)), \) where \( p_e = p^1 \) for all \( e \in  E(T^k) \) and  \(p_e = p^0 \) for all \( e \in E(T^{k,1})\smallsetminus E(T^k).\)  
	Let \( \nu_1 \) be the corresponding signed measure. 
	Let \( S_1 \coloneqq V_k \cup \mathcal{B}_{T_{k,1}}^+ \!\!\!(V_k). \)
	Note that if \( p^1 = 0\) then \( X^1(\{ o \} \cup (V_{k,1}\smallsetminus V_k )) \) is equal in distribution to \( X^0, \) and hence, by Lemma~\ref{lemma: finite to infinite}, we have \( \nu_1(S_1) = \nu_0(S_0) < 0 \) when \( p^1 = 0. \)
	Since \( \nu_1(S_1) \) is a continuous function in \( p^0 \) and \( p^1 \) and \( V_k \) is finite, it follows that there is \( \delta > 0 \) and a non-empty interval \( I_1 \subseteq I_0 \) such that \( \nu_1(S_1)<-\varepsilon/2 \) whenever \( p^1 \in (0,\delta) \) and \( p^0 \in I_1. \) Hence, for such \( p^0 \) and \( p^1, \) we have \( X^1 \notin \mathcal{R}. \)

	Now let \( X \) be a tree-indexed Markov chain on \( T \) with parameters \( (r,p). \) Assume that \( p \) is small enough to ensure \( p = p_1 \in (0,\delta) \) and that \( p_0=1-(1-p)^j \in I_1 \) for some \( j \geq 1. \) Fix any \( j \geq 1 \) such that this holds.
	Then \( X(V_{k,j}) \) is equal in distribution to \( X^1(S_1). \) Since \( X^1 \notin \mathcal{R} \), using Lemma~\ref{lemma: finite to infinite}, it follows that \( X \notin\mathcal{R}. \) This concludes the proof.
	\end{proof}

	%Let \(  T^{k,1} \) be the subgraph of \( T \) induced by \( V_{k,1} \) (see Figure~\ref{figure: last fig 3}). 
	%
	%
	%Note that the tree-indexed Markov chain on \(  T^{k,1} \) with resampling probabilities \( p \) on \( E(T^k) \) and \( 1-(1-p)^j \) on the remaining edges is equal in distribution to \( X ( V_{k,j} ) .\) 


\begin{thebibliography}{99}


\bibitem{fgs} Forsstr\"om, M. P., Gantert, N., Steif, J., Poisson Representable Processes, Probab. Theory Relat. Fields (2025).

\bibitem{fs} Forsstr\"om, M. P., Bethuelsen, S. A., Mixing for Poisson representable processes and consequences for the Ising model and the contact process, preprint (2025), available as arXiv:2501.14445

\bibitem{gbook} Georgii, H-O. Gibbs Measures and Phase Transitions, Berlin, New York: De Gruyter, (2011).

\bibitem{mp2003} Mossel, M., Peres, Y., Information flow on trees, The Annals of Applied Probability, Vol. 13, No. 3, 817--844 (2003).
 

\bibitem{st2019} Steif, J. E., Tykesson, J., Generalized Divide and Color Models, ALEA, Lat. Am. J. Probab. Math. Stat. 16, 1--57 (2019).


%\bibitem{sv2001} Subbarao, M. V., Verma, A., Some remarks on a product expansion: An unexplored partition function. In F. G. Garvan and M. E. H. Ismail, editors, Symbolic Computation, Number Theory, Special Functions, Physics and Combinatorics, Gainsville, Florida, 1999, pages 267--283. Kluwer, Dordrecht, 2001.



 

\end{thebibliography}
\end{document}